\documentclass[11pt]{article}
\usepackage{amsfonts,mathrsfs,amssymb,amsthm,mathptm}
\usepackage{amsmath,amscd}
\usepackage{mathptm,pslatex}
\usepackage{color}
\usepackage[dvips]{graphicx}
\usepackage[all]{xy}
\usepackage{graphicx}
\usepackage{fancyhdr}

\begin{document}
\newtheorem{Def}{Definition}[section]
\newtheorem{Bsp}[Def]{Example}
\newtheorem{Prop}[Def]{Proposition}
\newtheorem{Theo}[Def]{Theorem}
\newtheorem{Lem}[Def]{Lemma}
\newtheorem{Koro}[Def]{Corollary}
\theoremstyle{definition}
\newtheorem{Rem}[Def]{Remark}

\newcommand{\add}{{\rm add}}
\newcommand{\gd}{{\rm gl.dim}}
\newcommand{\dm}{{\rm dom.dim}}
\newcommand{\E}{{\rm E}}
\newcommand{\Mor}{{\rm Morph}}
\newcommand{\End}{{\rm End}}
\newcommand{\ind}{{\rm ind}}
\newcommand{\rsd}{{\rm res.dim}}
\newcommand{\rd} {{\rm rep.dim}}
\newcommand{\ol}{\overline}
\newcommand{\overpr}{$\hfill\square$}
\newcommand{\rad}{{\rm rad}}
\newcommand{\soc}{{\rm soc}}
\renewcommand{\top}{{\rm top}}
\newcommand{\pd}{{\rm proj.dim}}
\newcommand{\id}{{\rm inj.dim}}
\newcommand{\fld}{{\rm flat.dim}}
\newcommand{\Fac}{{\rm Fac}}
\newcommand{\Gen}{{\rm Gen}}
\newcommand{\fd} {{\rm fin.dim}}
\newcommand{\DTr}{{\rm DTr}}
\newcommand{\cpx}[1]{#1^{\bullet}}
\newcommand{\D}[1]{{\mathscr D}(#1)}
\newcommand{\Dz}[1]{{\mathscr D}^+(#1)}
\newcommand{\Df}[1]{{\mathscr D}^-(#1)}
\newcommand{\Db}[1]{{\mathscr D}^b(#1)}
\newcommand{\C}[1]{{\mathscr C}(#1)}
\newcommand{\Cz}[1]{{\mathscr C}^+(#1)}
\newcommand{\Cf}[1]{{\mathscr C}^-(#1)}
\newcommand{\Cb}[1]{{\mathscr C}^b(#1)}
\newcommand{\K}[1]{{\mathscr K}(#1)}
\newcommand{\Kz}[1]{{\mathscr K}^+(#1)}
\newcommand{\Kf}[1]{{\mathscr  K}^-(#1)}
\newcommand{\Kb}[1]{{\mathscr K}^b(#1)}
\newcommand{\modcat}{\ensuremath{\mbox{{\rm -mod}}}}
\newcommand{\Modcat}{\ensuremath{\mbox{{\rm -Mod}}}}

\newcommand{\stmodcat}[1]{#1\mbox{{\rm -{\underline{mod}}}}}
\newcommand{\pmodcat}[1]{#1\mbox{{\rm -proj}}}
\newcommand{\imodcat}[1]{#1\mbox{{\rm -inj}}}
\newcommand{\Pmodcat}[1]{#1\mbox{{\rm -Proj}}}
\newcommand{\Imodcat}[1]{#1\mbox{{\rm -Inj}}}
\newcommand{\opp}{^{\rm op}}
\newcommand{\otimesL}{\otimes^{\rm\mathbb L}}
\newcommand{\rHom}{{\rm\mathbb R}{\rm Hom}\,}
\newcommand{\projdim}{\pd}
\newcommand{\Hom}{{\rm Hom}}
\newcommand{\Coker}{{\rm coker}}
\newcommand{ \Ker  }{{\rm Ker}}
\newcommand{ \Img  }{{\rm Im}}
\newcommand{\Ext}{{\rm Ext}}
\newcommand{\StHom}{{\rm \underline{Hom}}}

\newcommand{\gm}{{\rm _{\Gamma_M}}}
\newcommand{\gmr}{{\rm _{\Gamma_M^R}}}

\def\vez{\varepsilon}\def\bz{\bigoplus}  \def\sz {\oplus}
\def\epa{\xrightarrow} \def\inja{\hookrightarrow}

\newcommand{\lra}{\longrightarrow}
\newcommand{\lraf}[1]{\stackrel{#1}{\lra}}
\newcommand{\ra}{\rightarrow}
\newcommand{\dk}{{\rm dim_{_{k}}}}

\newcommand{\colim}{{\rm colim\, }}
\newcommand{\limt}{{\rm lim\, }}
\newcommand{\Add}{{\rm Add }}
\newcommand{\Tor}{{\rm Tor}}
\newcommand{\Cogen}{{\rm Cogen}}

{\Large \bf
\begin{center}
Good tilting modules and recollements of derived module categories
\end{center}}
\medskip

\centerline{{\bf Hongxing Chen} and {\bf Changchang Xi$^*$}}
\begin{center} School of Mathematical Sciences, Beijing Normal University, \\
Laboratory of Mathematics and Complex Systems, \\
100875 Beijing, People's Republic of  China \\ E-mail: xicc@bnu.edu.cn \quad chx19830818@163.com\\
\end{center}
\bigskip

\renewcommand{\thefootnote}{\alph{footnote}}
\setcounter{footnote}{-1} \footnote{ $^*$ Corresponding author.
Email: xicc@bnu.edu.cn; Fax: 0086 10 58802136; Tel.: 0086 10
58808877.}
\renewcommand{\thefootnote}{\alph{footnote}}
\setcounter{footnote}{-1} 
\renewcommand{\thefootnote}{\alph{footnote}}
\setcounter{footnote}{-1} 

\begin{abstract}
Let $T$ be an  infinitely generated tilting  module of projective
dimension at most one over an arbitrary associative ring $A$, and
let $B$ be the endomorphism ring of $T$. In this paper, we prove
that if $T$ is good then there exists a ring $C$, a homological ring
epimorphism $B\ra C$ and a recollement among the (unbounded) derived
module categories $\D{C}$ of $C$, $\D{B}$ of $B$,  and $\D{A}$ of
$A$. In particular, the kernel of the total left derived functor
$T\otimes_B^{\mathbb L}-$ is triangle equivalent to the derived
module category $\D{C}$. Conversely, if the functor
$T\otimes_B^{\mathbb L}-$ admits a fully faithful left adjoint
functor, then $T$ is a good tilting module. We apply our result to
tilting modules arising from ring epimorphisms, and can then
describe the rings $C$ as coproducts of two relevant rings. Further,
in case of commutative rings, we can weaken the condition of being
tilting modules, strengthen the rings $C$ as tensor products of two
commutative rings, and get similar recollements. Consequently, we
can produce examples (from commutative algebra and $p$-adic number
theory, or Kronecker algebra) to show that two different
stratifications of the derived module category of a ring by derived
module categories of rings may have completely different derived
composition factors (even up to ordering and up to derived
equivalence),or different lengths. This shows that the
Jordan-H\"older theorem fails even for stratifications by derived
module categories, and also answers negatively an open problem by
Angeleri-H\"ugel, K\"onig and Liu.
\end{abstract}

\tableofcontents

\section{Introduction}
The theory of finitely generated tilting modules has been
successfully applied, in the representation theory of algebras and
groups, to understanding different aspects of algebraic structure
and homological features of (algebraic) groups, algebras and modules
(for instance, see \cite{BI, BB, CPS, donkin},
\cite{Happel}-\cite{Keller}). Recently, infinitely generated tilting
modules over arbitrary associated rings have become of interest in
and attracted increasingly attentions toward to understanding
derived categories and equivalences of general rings
(\cite{HA}-\cite{HJ2}, \cite{Bz}-\cite{Bz2}, \cite{ct, CTT},
\cite{R}-\cite{TP}). In this general situation, many classical
results in the tilting theory appear in a very different new
fashion. For example, Happel's Theorem (see also \cite{CPS}) on
derived equivalences induced by infinitely generated tilting modules
comes up with a new formulation in which quotient categories are
involved (see \cite{Bz}). This more general context of tilting
theory not only renews our view on features of finitely generated
tilting modules, but also provides us completely different
information about the whole tilting theory. Let us recall the
definition of tilting modules over an arbitrary ring from \cite{ct}.

Let $A$ be a ring with identity, and let $T$ be a left $A$-module
which may be infinitely generated. The module $T$ is called a
tilting module (of projective dimension at most $1$) provided that

$(T1)$ $T$ has projective dimension at most one,

$(T2)$ $\Ext^i_A(T,T^{({\alpha})})=0$ for each $i\geq 1$ and each
cardinal $\alpha$, and

$(T3)$ there exists  an exact sequence $0\rightarrow A\rightarrow
T_0\rightarrow T_1\rightarrow 0$ of left $A$-modules, where  $T_0$
and $T_1$ are isomorphic to direct summand of arbitrary direct sums
of copies of $T$.

\medskip
If, in addition, $T$ is finitely presented, then we  say that $T$ is
a classical tilting module. If the modules $T_0$ and $T_1$ in $(T3)$
are isomorphic to direct summands of finite direct sums of copies of
$T$, then we say that $T$ is a good tilting module, following
\cite{Bz2}.  Actually, each classical tilting module is good,
furthermore, it is proved in \cite{Bz2} that, for an arbitrary
tilting $A$-module $T$, there exists a good tilting $A$-module $T'$
which is equivalent to $T$, that is, $T$ and $T'$ generate the same
full subcategories in the category of all left $A$-modules.

One of the realizations of tilting modules is universal
localizations. It is shown in \cite{HA} that every tilting module
over a ring is associated in a canonical manner with a ring
epimorphism which can be interpreted as a universal localization at
a set of homomorphisms between finitely presented modules of
projective dimension at most one.

As in the theory of classical tilting modules, a natural context for
studying infinitely generated tilting modules is the relationship of
derived categories and equivalences induced by infinitely generated
tilting modules. In fact, if $T$ is a good tilting module over  a
ring $A$, and if $B$ is the endomorphism ring of $T$, then Bazzoni
 proves in \cite{Bz} that the total right derived functor $\rHom_A(T,-)$
induces an equivalence between the (unbounded) derived category
$\D{A}$ of $A$ and the quotient category of the derived category
$\D{B}$ of $B$ modulo the full triangulated subcategory
$\Ker(T\otimesL_{B}-)$ which is the kernel of the total left derived
functor $T\otimesL_{B}-$. Thus, in general, the total right derived
functor $\rHom_A(T,-)$ does not define a derived equivalence between
$A$ and $B$. This is a contrary phenomenon to the classical
situation (see \cite{CPS}). The condition for $A$ and $B$ to be
derived-equivalent depends on the vanishing of
$\Ker(T\otimesL_{B}-)$. It is shown in \cite{Bz} that
$\Ker(T\otimesL_{B}-)$ vanishes if and only if $T$ is a classical
tilting module. From this point of view, the triangulated category
$\Ker(T\otimesL_{B}-)$ measures how far a good tilting module is
from being classical, in other words, the difference between the two
derived categories $\D{A}$ and $\D{B}$. It is certainly of interest
to have a little bit knowledge about the categories
$\Ker(T\otimesL_{B}-)$ for infinitely generated tilting modules $T$.
This might help us to understand some new aspects of the tilting
theory of infinitely generated tilting modules.

The main purpose of this paper is to give a characterization of the
triangulated categories $\Ker(T\otimesL_{B}-)$ for infinitely
generated tilting modules $T$, namely, we show that if the tilting
module $T$ is good then the triangulated category
$\Ker(T\otimesL_{B}-)$ is equivalent to the derived category of a
ring $C$, and therfore, there is a recollemment among the derived
categories of rings $A$, $B$ and $C$. Conversely, the existence of
such a recollement implies that the given tilting module $T$ is
good. More precisely, our result can be stated as follows:

\begin{Theo}
Let $A$ be a ring, $T$ a  tilting $A$-module of projective dimension
at most $1$ and $B$ the endomorphism ring of $T$.

$(1)$ If $T$ is good, then  there is a ring $C$, a homological ring
epimorphism $\lambda: B\to C$ and  a recollement among the unbounded
derived categories of the rings $A$, $B$ and $C$:
$$\xymatrix{\D{C}\ar[r]&\D{B}\ar^{j^!}[r]\ar@/^1.2pc/[l]\ar@/_1.2pc/[l]
&\D{A}\ar@/^1.2pc/[l]\ar@/_1.2pc/[l]}.$$ such that the triangle
functor $j^!$ is isomorphic to the total left derived functor
${_A}T\otimesL_B-$. In this case, the kernel of the functor
$T\otimesL_B-$ is equivalent to the unbounded derived category
$\D{C}$ of $C$ as triangulated categories.

$(2)$ If the triangle functor $T\otimesL_B-: \D{B}\ra \D{A}$ admits
a fully faithful left adjoint $j_!:\D{A}\to\D{B}$, then the given
tilting module $T$ is good. \label{th01}
\end{Theo}

Let us remark that a noteworthy difference of Theorem \ref{th01}(1)
from the result \cite[Proposition 1.7]{HKL} is that our recollement
is  over derived module categories of precisely determined rings,
while the recollement in \cite[Proposition 1.7]{HKL} involves a
triangulated category. Theorem \ref{th01}(1) realizes this abstract
triangulated category by a derived module category via describing
the kernel of the functor $T\otimesL_{B}-$. Our result also
distinguishes itself from the one in \cite{DY} where $C$ is a
differential graded ring instead of a usual ring, and where the
consideration is restricted to ground ring being a field.

If we apply Theorem \ref{th01} to tilting modules arising from ring
epimorphisms, then we can see that, in most cases, the recollements
given in Theorem \ref{th01} are different from the usual ones
induced from the structure of triangular matrix rings. The following
corollary is a consequence of the proof of Theorem \ref{th01}.

\begin{Koro} $(1)$ Let $R\ra S$ be an injective
ring epimorphism such that $\Tor_1^R(S,S)=0$ and that $_RS$ has
projective dimension at most one. Then there is a recollement of
derived module categories:

$$\xymatrix@C=1.2cm{\D{S\sqcup_RS'}\ar[r]&\D{\End_R(S\oplus S/R)}\ar[r]
\ar@/^1.2pc/[l]\ar@/_1.2pc/[l]
&\D{R}\ar@/^1.2pc/[l]\ar@/_1.2pc/[l]},\vspace{0.3cm}$$ where $S'$ is
the endomorphism ring of the $R$-module $S/R$, and $S\sqcup_RS'$ is
the coproduct of $S$ and $S'$ over $R$.

$(2)$ Suppose that $\lambda: R\ra S$ is an injective homological
ring epimorphism between commutative rings $R$ and $S$. Then there
is a recollement of derived module categories:
$$\xymatrix@C=1.2cm{\D{S\otimes_RS'}\ar[r]&\D{\End_R(S\oplus S/R)}\ar[r]
\ar@/^1.2pc/[l]\ar@/_1.2pc/[l]
&\D{R}\ar@/^1.2pc/[l]\ar@/_1.2pc/[l]},\vspace{0.3cm}$$ where
$S':=\End_R(S/R)$ is a commutative ring, and $S\otimes_RS'$ is the
tensor product of $S$ and $S'$ over $R$.

$(3)$ For every prime number $p\ge 2$, the derived category of the
ring $\left(\begin{array}{cc} {\mathbb Q} & {\mathbb Q}_p
\\ 0 & {\mathbb Z}_p\end{array}\right)$ admits two stratifications,
one of which clearly has composition factors $\mathbb Q$ and
${\mathbb Z}_p$, and the other has composition factors ${\mathbb
Q}_{(p)}$ and ${\mathbb Q}_p$, where ${\mathbb Q}_{(p)}$, $\mathbb
Q$, ${\mathbb Z}_p$ and ${\mathbb Q}_p$ denote the rings of
$p$-integers, rational numbers, $p$-adic integers and $p$-adic
numbers, respectively. \label{1.2}
\end{Koro}

As pointed out in \cite{akl}), the Jordan-H\"older theorem fails for
stratifications of derived module categories by triangulated
categories. Our Corollary \ref{1.2}(3) (see also the example in
Section \ref{sect8} below) shows that the Jordan-H\"older theorem
fails even for stratifications of derived module categories by
derived module categories, and therefore the problem posed in
\cite{akl} gets a negative answer.

The paper is organized as follows: In Section \ref{sect2}, we recall
some definitions, notations and useful results which are needed for
our proofs. In Section \ref{sect3}, we shall first establish a
connection between universal localizations and recollements of
triangulated categories, and then prove Proposition \ref{th2} which
is crucial for the proof of the main result. In Section \ref{sect4},
we discuss some homological properties of good tilting modules, and
establish another crucial result, Proposition \ref{prop1}, for the
proof of the main result Theorem \ref{th01}. After these
preparations, we apply the results obtained in Section \ref{sect3}
to prove Theorem \ref{th01}(1). In Section \ref{sect5}, we prove the
second part of Theorem \ref{th01}. This may be regarded as a
converse statement of the first part. In Section \ref{6}, we apply
Theorem \ref{th01} to good tilting modules arising from ring
epimorphisms, and prove Corollary \ref{1.2}(1). In these cases the
universal localization rings in Theorem \ref{th01} can be given by
coproducts of rings. Our discussion in this section is actually
carried out under the general assumption of injective homological
ring epimorphisms. In Section \ref{sect7}, we first consider the
existence of the recollements in Theorem \ref{th01} for commutative
rings without assumption that the involved modules are tilting
modules, and then make special consideration of localizations of
commutative one-Gorestein rings. In particular, we prove Corollary
\ref{1.2}(2) and Corollary \ref{1.2}(3). It turns out that many
derived module categories of rings possess stratifications by
derived module categories of rings, such that, even up to ordering
and up to derived equivalence, not all of their composition factors
are the same; for instance, the derived category of the endomorphism
ring of the abelian group ${\mathbb Q}\oplus {\mathbb Q}/{\mathbb
Z}$ (or its variation ${\mathbb Q}\oplus {\mathbb Q}/{\mathbb
Q}_{(p)}$). Note that, in the examples presented in this section,
the two stratifications all have the same lengths. In Section
\ref{sect8}, we give an example of a non-commutative algebra over
which the derived category of the endomorphism ring of a tilting
module has two stratifications of different finite lengths. This,
together with the examples in Section \ref{sect7}, gives a complete
answer to an open problem in \cite{akl} negatively.

\medskip
The research work of the corresponding author C.C.Xi is partially
supported by the Fundamental Research Funds for the Central
Universities (2009SD-17), while the author H.X.Chen is supported by
the Doctor Funds of the Beijing Normal University. Also, C.C.Xi
thanks Lidia Angeleri-H\"ugel for some discussions on localizations
of commutative rings.

\section{Preliminaries\label{sect2}}
In this section, we shall recall some definitions, notations and
basic results which are related to our proofs. In particular, we
recall the notions of recollements and TTF triples as well as their
relationship.

\subsection{Some conventions}
In this  subsection, we recall some standard notations which will be
used throughout this paper.

All rings considered in this paper are assumed to be associative and
with identity, and all ring homomorphisms preserve identity.

Let $A$ be a ring. We denote by $A\Modcat$ the category of all
unitary left $A$-modules. For an $A$-module $M$, we denote by
$\add(M)$ (respectively, $\Add(M)$) the full subcategory of
$A\Modcat$ consisting of all direct summands of finite
(respectively, arbitrary) direct sums of copies of $M$. In many
circumstances, we shall write $A\pmodcat$ and $A\Pmodcat$ for
$\add(_AA)$ and $\Add(_AA)$, respectively. If $I$ is an index, we
denote by $M^{(I)}$ the direct sum of $I$ copies of $M$. If there is
a surjective homomorphism from $M^{(I)}$ to an $A$-module $X$, we
say that $X$ is generated by $M$, or $M$ generates $X$. By $\Gen(M)$
we denote the full subcategory of $A\Modcat$ generated by $M$.

If $f: M\ra N$ is a homomorphism of $A$-modules, then the image of
$x\in M$ under $f$ is denoted by $(x)f$ instead of $f(x)$. Also, for
any $A$-module $X$, the induced morphisms $\Hom_A(X,f):
\Hom_A(X,M)\ra \Hom_A(X,N)$ and $\Hom_A(f,X): \Hom_A(N, X)\ra
\Hom_A(M, X)$ is denoted by $f^*$ and $f_*$, respectively.

Let $\mathcal C$ be an additive category.

Given two morphisms $f: X\to Y$ and $g: Y\to Z$ in $\mathcal C$, we
denote the composition of $f$ and $g$ by $fg$ which is a morphism
from $X$ to $Z$, while we denote the composition of a functor
$F:\mathcal {C}\to \mathcal{D}$ between categories $\mathcal C$ and
$\mathcal D$ with a functor $G: \mathcal{D}\to \mathcal{E}$ between
categories $\mathcal D$ and $\mathcal E$ by $GF$ which is a functor
from $\mathcal C$ to $\mathcal E$. The image of the functor $F$ is
denoted by Im$(F)$ which is a full subcategory of $\mathcal D$.

Throughout the paper, a full subcategory $\mathcal D$ of $\mathcal
C$ is always assumed to be closed under isomorphisms, that is, if
$X$ and $Y$ are objects in $\cal C$, then $Y\in{\mathcal D}$
whenever $Y\simeq X$ with $X\in  {\mathcal D}$.

Let $\mathcal{Y}$ be a full subcategory of $\mathcal{C}$. By
$\Ker(\Hom_{\mathcal{C}}(-,\mathcal{Y}))$ we denote the left
orthogonal subcategory with respect to $\mathcal{Y}$, that is, the
full subcategory of $\mathcal{C}$ consisting of the objects $X$ such
that $\Hom_{\mathcal{C}}(X,Y)=0$ for all objects $Y$ in
$\mathcal{Y}$. Similarly, $\Ker(\Hom_{\mathcal{C}}(\mathcal{Y},-))$
stands for the right orthogonal subcategory of $\cal C$ with respect
to $\mathcal{Y}$.

By a complex $\cpx{X}$ over $\mathcal{C}$ we mean a sequence of
morphisms $d_X^i$ between objects $X^i$ in $\mathcal{C}:\; \cdots\to
X^i\epa{d_X^i } X^{i+1}\epa{d_X^{i+1}} X^{i+2}\to\cdots$, such that
$d_X^id_X^{i+1}=0$ for all $i\in\mathbb{Z}$. In this case, we write
$\cpx{X}=(X^i, d_X^i)_{i\in\mathbb{Z}}$, and call $d_X^i$ a
differential of $\cpx{X}$. Sometimes, for simplicity, we write
$(X^i)_{i\in\mathbb{Z}}$ for $\cpx{X}$ without mentioning the
morphisms $d^i_X$. For a fixed integer $n$, we denote by
$\cpx{X}[n]$ the complex obtained from $\cpx{X}$ by shifting $n$
degrees, that is, $(\cpx{X}[n])^0=X^n$, and $H^n(\cpx{X})$ the
cohomology of $\cpx{X}$ in degree $n$.

Let $\C{\mathcal{C}}$ be the category of all complexes over
$\mathcal{C}$ with chain maps, and $\K{\mathcal{C}}$ the homotopy
category of $\C{\mathcal{C}}$. We denote by $\Cb{\mathcal{C}}$ and
$\Kb{\mathcal{C}}$ the full subcategories of $\C{\mathcal{C}}$ and
$\K{\mathcal{C}}$ consisting of bounded complexes over
$\mathcal{C}$, respectively. When $\mathcal{C}$ is abelian, the
derived category of $\mathcal{C}$ is denoted by $\D{\mathcal{C}}$,
which is the localization of $\K{\mathcal C}$ at all
quasi-isomorphisms. The full subcategory of $\D{\mathcal{C}}$
consisting of bounded complexes over $\mathcal{C}$ is denoted by
$\Db{\mathcal{C}}$. As usual, for a ring $A$, we simply write
$\C{A}$ for $\C{A\Modcat}$, $\K{A}$ for $\K{A\Modcat}$, $\Cb{A}$ for
$\Cb{A\Modcat}$, and $\Kb{A}$ for $\Kb{A\Modcat}$. Similarly, we
write $\D{A}$ and $\Db{A}$ for $\D{A\Modcat}$ and $\Db{A\Modcat}$,
respectively. Furthermore, we always identify $A\Modcat$ with the
full subcategory of $\D{A}$ consisting of all stalk complexes
concentrated on degree zero.

Now we recall some  basic facts about derived functors defined on
derived module categories. We refer to \cite{bn} for details and
proofs.

Let $R$ and $S$ be rings, and let $H$ be an additive functor from
$R$-Mod to $S$-Mod.

(1) For each complex $\cpx{X}$ in $\D{R}$, there is a complex
$\cpx{I}\in \C{\Imodcat{R}}$ such that $\cpx{X}$ is quasi-isomorphic
to $\cpx{I}$, where $\Imodcat{R}$ is the full subcategory of $R$-Mod
consisting of all injective $R$-modules. Dually, for each complex
$\cpx{Y}$ in $\D{R}$, there is a complex $\cpx{P}\in
\C{\Pmodcat{R}}$ such that $\cpx{P}$ is quasi-isomorphic to
$\cpx{Y}$.

(2) There is a total right derived functor ${\mathbb R}H$ and a
total left derived functor ${\mathbb L}H$ defined on $\D{R}$. If
$\cpx{X}, \cpx{Y}\in \D{R}$,  then ${\mathbb R}H(\cpx{X})=
H(\cpx{I})$ and ${\mathbb L}H(\cpx{Y})=H(\cpx{P})$, where $\cpx{I}$
and $\cpx{P}$ are chosen as in (1). Here we think of $H$ as an
induced functor between homotopy categories, and if
$\cpx{X}=(X^i,d_X^i)_{i\in \mathbb Z}$ then $H(\cpx{X})
:=\big(H(X^i),H(d^i_X)\big)_{i\in \mathbb Z}$.

In case $T$ is an $R$-$S$-bimodule, the total right derived functor
of $\Hom_R(T,-)$ is denoted by ${\mathbb R}\Hom_A(T,-)$, and the
total left derived functor of $T\otimes_B-$ is denoted by
$T\otimesL_B-$.

(3) Any adjoint pair of functors $(G,H)$ between $R$-Mod and $S$-Mod
induces an adjoint pair $({\mathbb L}G,{\mathbb R}H)$ between the
unbounded derived categories of $R$ and $S$.

\subsection{Homological ring epimorphisms\label{epimorphism}}

Let $R$ and $S$ be rings. Recall that a homomorphism $\lambda: R\ra
S$ of rings is called a ring epimorphism if, for any two
homomorphisms $f_1,f_2: S\ra T$ of rings, the equality $\lambda
f_1=\lambda f_2$ implies that $f_1=f_2$. It is known that $\lambda$
is a ring epimorphism if and only if the multiplication map
$S\otimes_RS\ra S$ is an isomorphism as $S$-$S$-bimodules if and
only if $x\otimes 1=1\otimes x$ in $S\otimes_RS$ for any $x\in S$.
It follows that, for a ring epimorphism, we have $X\otimes_SY\simeq
X\otimes_RY$ for any $S$-modules $X_S$ and $_SY$. An example of ring
epimorphisms is the inclusion ${\mathbb Z}\hookrightarrow {\mathbb
Q}$. Note that $\mathbb Q$ is an injective and a flat $\mathbb
Z$-module.

Given a ring epimorphism $\lambda: R\ra S$ between two rings $R$ and
$S$, we can regard $S$-Mod as a full subcategory of $R$-Mod via
$\lambda$. This means that $\Hom_S(X,Y)\simeq \Hom_R(X,Y)$ for all
$S$-modules $X$ and $Y$.

Two ring epimorphisms $\lambda: R\ra S$ and $\lambda': R\ra S'$ are
said to be equivalent if there is a ring isomorphism $\psi:S\to S'$
such that $\lambda'=\lambda\psi$. This defines an equivalence
relation on the class of ring epimorphisms $R\ra S$ with $R$ fixed.
The equivalence classes with respect to this equivalence relation
are called the epiclasses of $R$. This notion is associated with
bireflective subcategories of module categories.

Recall that a full subcategory $\mathcal D$ of $R\Modcat$ is said to
be reflective if every $R$-module $X$ admits a
$\mathcal{D}$-reflection, that is, there exists an $R$-module
$D'\in\mathcal{D}$ and a homomorphism $f:D'\to X$ of $R$-modules
such that $\Hom_R(D,f):\Hom_R(D,D')\to\Hom_R(D,X)$ is an isomorphism
as abelian groups for any module $D\in\mathcal{D}$. Dually, one
defines the notion of  coreflective subcategories of $R\Modcat$. The
full subcategory $\mathcal D$ of $R\Modcat$ is called bireflective
if it is both reflective and coreflective.

Ring epimorphisms are related to bireflective subcategories in the
following way.

\begin{Lem}{\rm  \cite[Theorm 1.4]{HA}}
For a full subcategory $\mathcal{D}$ of $R\Modcat$, the following
statements are equivalent.

$(1)$ There is  a ring epimorphism $\lambda: R\to S$ such that the
category  $\mathcal{D}$ is the image of the restriction functor
\indent\indent $\lambda_*:S\Modcat\to R\Modcat$.

$(2)$ $\mathcal{D}$ is a bireflective subcategory of $R\Modcat$.

$(3)$ $\mathcal{D}$ is closed under direct sums, products, kernels
and cokernels.

\noindent Thus, there is a bijection between the epiclasses of $R$
and the bireflective subcategories of $R$-Mod. Furthermore, the map
$\lambda: R\to S$ in $(1)$, viewed as a homomorphism of $R$-modules,
is a $\mathcal{D}$-reflection of $R$. \label{2.0}
\end{Lem}

Following Geigle and Lenzing \cite{GL},  we say that a ring
epimorphism $\lambda: R\ra S$ is homological if $\Tor^R_i(S, S)=0$
for all $i>0$. This is equivalent to saying that the restriction
functor $\lambda_*:\D{S}\to\D{R}$ induced by $\lambda$ is fully
faithful. In \cite[Theorem 4.4]{GL}, the following lemma is proved.

\begin{Lem} {\rm \cite[Theorem 4.4]{GL}}
For a homomorphism $\lambda: R\ra S$ of rings, the following
assertions are equivalent:

$(1)$ $\lambda$ is homological,

$(2)$ For all right $S$-modules $X$ and all left $S$-modules $Y$,
the natural map $\Tor_i^R(X,Y)\ra \Tor^S_i(X,Y)$ is an isomorphism
for all $i\ge 0$.

$(3)$ For all $S$-modules $X$ and $Y$, the natural map
$\Ext_S^i(X,Y)\ra \Ext^i_R(X,Y)$ is an isomorphism for all $i\ge 0$.
\label{2.1}
\end{Lem}

Note that the condition (3) in Lemma \ref{2.1} can be replaced by
the corresponding version of right modules. For more details, one
may look at \cite{GL} and \cite[Section 5.3]{NS2}.

On ring epimorphisms, we have the following property which will be
used in Section \ref{sect7}.

\medskip
\begin{Lem}
Let $g:\Lambda\to\Gamma$ and $h:\Gamma\to\Delta$ be ring
homomorphisms such that $gh:\Lambda\to\Delta$ is a ring epimorphism.
Then $h$ is a ring epimorphism. Suppose further that $h$ is
injective. If $\Gamma_\Lambda$ and $_\Lambda\Delta$ (respectively,
$_\Lambda\Gamma$ and $\Delta_\Lambda$) are flat, then both $g$ and
$h$ are homological ring epimorphisms. \label{prop6.7}
\end{Lem}

{\it Proof.} By the  definition of ring epimorphisms, we can readily
show that $h$ is a ring epimorphism. Note that we always have the
following commutative diagram:
$$\xymatrix{\Gamma\otimes_\Lambda\Gamma\ar[r]^-{h\otimes_\Lambda
h}\ar[d]_-{\mu_1}&\Delta\otimes_\Lambda\Delta\ar[d]_-{\mu_2}\\
\Gamma\ar[r]^-{h}&\Delta,}$$ where $\mu_1$ and $\mu_2$ are the
canonical multiplication maps.  Suppose that $h$ is injective. If
$\Gamma_\Lambda$ and $_\Lambda\Delta$ are flat, then the map
$h\otimes_\Lambda h$ is injective. Since $gh:\Lambda\to\Delta$ is a
ring epimorphism, the map $\mu_2$ is an isomorphism. It follows that
$\mu_1$ is injective, and therefore it is an isomorphism. This means
that $g:\Lambda\to\Gamma$ is a ring epimorphism. Note that
$\Gamma_\Lambda$ is flat. Thus $g$ is a homological ring
epimorphism. To prove that $h$ also is a homological ring
epimorphism, we claim that $_\Gamma\Delta$ is flat. In fact, this
follows from Lemma \ref{2.1} because $g$ is a homological ring
epimorphism and because $\Delta$ is flat as a $\Lambda$-module.
Similarly, we can prove that if $_\Lambda\Gamma$ and
$\Delta_\Lambda$ are flat, then both $g$ and $h$ are homological
ring epimorphisms. $\square$

\subsection{Recollements and TTF triples}\label{sect2.2}

In this subsection, we first recall the definitions of recollements
and TTF triples, and then state a correspondence between them.

From now on, $\mathcal{D}$ denotes a triangulated category with
small coproducts (that is, coproducts indexed over a set), and $[1]$
the shift functor of  $\mathcal{D}$.

The notion of recollements was first defined by Beilinson, Bernstein
and Deligne in \cite{BBD} to study ``exact sequences" of derived
categories of coherent sheaves over geometric objects.

\begin{Def}\rm
Let $\mathcal{D'}$ and $\mathcal{D''}$ be triangulated categories.
We say that $\mathcal{D}$ is a recollement of $\mathcal{D'}$ and
$\mathcal{D''}$ if there are six triangle functors as in the
following diagram
$$\xymatrix{\mathcal{D''}\ar^-{i_*=i_!}[r]&\mathcal{D}\ar^-{j^!=j^*}[r]
\ar^-{i^!}@/^1.2pc/[l]\ar_-{i^*}@/_1.2pc/[l]
&\mathcal{D'}\ar^-{j_*}@/^1.2pc/[l]\ar_-{j_!}@/_1.2pc/[l]}$$ such
that

$(1)$ $(i^*,i_*),(i_!,i^!),(j_!,j^!)$ and $(j^*,j_*)$ are adjoint
pairs;

$(2)$ $i_*,j_*$ and $j_!$ are fully faithful functors;

$(3)$ $i^!j_*=0$ (and thus also $j^! i_!=0$ and $i^*j_!=0$); and

$(4)$ for each object $C\in\mathcal{D}$, there are two triangles in
$\mathcal D$:
$$
i_!i^!(C)\lra C\lra j_*j^*(C)\lra i_!i^!(C)[1],
$$
$$
j_!j^!(C)\lra C\lra i_*i^*(C)\lra j_!j^!(C)[1].
$$
\label{def01}
\end{Def}

Recollements are closely related to TTF triples which are defined in
terms of torsion pairs. So, let us first recall the notion of
torsion pairs in triangulated categories.

\begin{Def}\rm \cite{BI} A torsion pair in $\mathcal{D}$ is a pair
$(\mathcal{X},\mathcal{Y})$ of full subcategories $\mathcal{X}$ and
$\mathcal{Y}$ of $\mathcal{D}$ satisfying the following conditions:

$(1)$ $\Hom_{\mathcal{D}}(\mathcal{X},\mathcal{Y})=0$;

$(2)$ $\mathcal{X}[1]\subseteq\mathcal{X}$ and
$\mathcal{Y}[-1]\subseteq\mathcal{Y}$; and

$(3)$ for each object $C\in\mathcal{D}$, there is a triangle
$$
X_C\lra C\lra Y^C\lra X_C[1]
$$
in $\mathcal{D}$ such that $X_C\in\mathcal{X}$ and
$Y^C\in\mathcal{Y}$. In this case, $\mathcal{X}$ is called a torsion
class and $\mathcal{Y}$ is called a torsion-free class. If, in
addition, $\mathcal{X}$ is a triangulated subcategory  of
$\mathcal{D}$ (or equivalently, $\mathcal{Y}$ is a triangulated
subcategory of $\mathcal{D}$), then  the torsion pair
$(\mathcal{X},\mathcal{Y})$ is said to be  hereditary (see
\cite[Chapter I, Proposition 2.6]{BI}). \label{def02}
\end{Def}

Note that, if $(\mathcal{X},\mathcal{Y})$ is a torsion pair in
$\mathcal{D}$, then
$\mathcal{X}=\Ker(\Hom_{\mathcal{C}}(-,\mathcal{Y}))$ which is
closed under  small coproducts, and
$\mathcal{Y}=\Ker(\Hom_{\mathcal{C}}(\mathcal{X},-))$ which is
closed under small products.

\begin{Def}\rm \cite{BI} A torsion torsionfree triple, or {\rm TTF} triple for short,  in
$\mathcal{D}$ is a triple $(\mathcal{X},\mathcal{Y},\mathcal{Z})$ of
full subcategories $\mathcal{X},\mathcal{Y}$ and $\mathcal{Z}$ of
$\mathcal{D}$ such that both $(\mathcal{X},\mathcal{Y})$ and
$(\mathcal{Y},\mathcal{Z})$ are torsion pairs. In this case,
$\mathcal{X}$  is said to be a smashing subcategory of
$\mathcal{D}$. \label{def03}
\end{Def}

It follows from \cite[Chapter I.2.]{BI} that, associated with a TTF
triple $(\mathcal{X},\mathcal{Y},\mathcal{Z})$ in $\mathcal{D}$,
there are seven triangle functors demonstrated in the following
diagram
$$\xymatrix{\mathcal{X}\ar^-{\bf{i}}@/^0.8pc/[r]&\mathcal{D}
\ar^-{\bf{R}}@/^0.8pc/[l]\ar^-{\bf{L}}@/^0.8pc/[r]
&\mathcal{Y}\ar^-{\bf{j}}@/^0.8pc/[l]
&\mathcal{D}\ar_-{\bf{V}}@/^0.8pc/[l]\ar^-{\bf{U}}@/^0.8pc/[r]
&\mathcal{Z}\ar^-{\bf{k}}@/^0.8pc/[l]}$$ such that

$(1)$ $\bf{i},\bf{j}$ and $\bf{k}$ are canonical inclusions; and

$(2)$ $(\bf{i},\bf{R}),(\bf{L},\bf{j}),(\bf{j},\bf{V})$ and
$(\bf{U},\bf{k})$ are adjoint pairs; and

$(3)$ the  composition functor $\bf{Ui}:\mathcal{X}\to\mathcal{Z}$
of th functors $\bf{i}$ and $\bf{U}$ is a triangle equivalence with
the quasi-inverse functor $\bf{Rk}$ which is the composition of the
functors $\bf{k}$ and $\bf{R}$.

\medskip
Note that if $(\mathcal{X},\mathcal{Y},\mathcal{Z})$ is a TTF triple
in $\mathcal{D}$, then it is easy to check that $\mathcal X$,
$\mathcal Y$ and $\mathcal Z$ are automatically triangulated
subcategories of $\mathcal D$.

Observe also that the existence of the  functors $\bf{R}$ and
$\bf{L}$ in the above diagram follows from the fact that
$(\mathcal{X},\mathcal{Y})$ is a torsion pair in $\mathcal{D}$ (see
\cite[Chapter I, Proposition 2.3]{BI} for details). Furthermore,
$\mathcal{Y}$ is closed under small coproducts and products.

\medskip
Now, we state a correspondence between recollements and TTF triples
given in \cite[Section 9.2]{AN2} and \cite[Section 4.2]{NS2}. For
more details, we refer the reader to these papers.

\begin{Lem}
$(1)$  If $\mathcal{D}$ is a recollement of $\mathcal{D'}$ and
$\mathcal{D''}$ in Definition \ref{def01}, then
$\big(j_!(\mathcal{D}'),i_*(\mathcal{D}''),j_*(\mathcal{D}')\big)$
is a {\rm{TTF}} triple in $\mathcal{D}$.

$(2)$ If $(\mathcal{X},\mathcal{Y},\mathcal{Z})$ is a {\rm{TTF}}
triple in $\mathcal{D}$, then $\mathcal{D}$ is a recollement of
$\mathcal{X}$ and $\mathcal{Y}$ as follows:
$$\xymatrix{\mathcal{Y}\ar^-{\bf{j}}[r]&\mathcal{D}\ar^-{\bf{R}}[r]
\ar^-{\bf{V}}@/^1.2pc/[l]\ar_-{\bf{L}}@/_1.2pc/[l]
&\mathcal{X}\ar^-{\bf{kUi}}@/^1.2pc/[l]\ar_-{\bf{i}}@/_1.2pc/[l]}.$$
\label{lem01}
\end{Lem}

\subsection{Generators and compact objects}
In this subsection, we shall recall some definitions and facts on
generators in a triangulated category.

Given a class of objects $\mathcal{U}$ in $\mathcal{D}$, we denote
by $\rm{Tria}(\mathcal{U})$  the smallest full triangulated
subcategory of $\mathcal{D}$ which contains $\mathcal{U}$ and is
closed under small coproducts. If $\mathcal{U}$ consists of only one
single object $U$, then we simply write $\rm{Tria}(U)$ for
$\rm{Tria}(\{U\})$.

\begin{Def}\rm
A class $\mathcal{U}$ of objects in $\mathcal{D}$ is called a class
of generators of $\mathcal{D}$ if an object $D$ in $\mathcal{D}$ is
zero whenever $\Hom_{\mathcal{D}}(U[n],D)=0$ for every object $U$ of
$\mathcal{U}$ and every $n$ in $\mathbb{Z}$.

An object $P$ in $\mathcal{D}$ is called compact if the functor
$\Hom_{\mathcal{D}}(P,-)$ preserves small coproducts, that is,
$\Hom_{\mathcal D}(P, {}\oplus_{i\in I}X_i)\simeq \oplus_{i\in
I}\Hom_{\mathcal D}(P,X_i)$, where $I$ is a set;  and exceptional if
$\Hom_{\mathcal D}(P,P[i])=0$ for all $i\ne 0$. The object $P$ is
called a tilting object if $P$ is compact, exceptional and a
generator of $\mathcal D$. Note that, for a compact generator $P$,
we have Tria$(P) = \mathcal D$ (see \cite{NS2}, for instance).

The category $\mathcal{D}$ is said to be compactly generated if
$\mathcal{D}$ admits a set $\mathcal{V}$ of compact generators. In
this case, $\mathcal{D}=\rm{Tria}(\mathcal{V})$, and we say that
$\mathcal{D}$ is compactly generated by $\mathcal{V}$.
 \label{def04}
\end{Def}

It is well-known that, for a ring $A$, the unbounded derived
category $\D{A}$ is a compactly generated triangulated category, and
one of its compact generators is $_AA$. Moreover, a complex
$\cpx{P}\in\D{A}$ is compact if and only if it is quasi-isomorphic
to a bounded complex of finitely generated projective $A$-modules.
The importance of compact objects can be seen from the following
lemma, due to Keller in \cite[Corollary 8.4, Theorem 8.5]{Keller}.

\begin{Lem}  Let $A$ be a ring. If $\cpx{P}$ is a compact exceptional object in
$\D{A}$, then ${\rm Tria}(\cpx{P})$ is equivalent to
$\D{\End_{\D{R}}(\cpx{P})}$ as triangulated categories.
\label{compact-deriv}
\end{Lem}

The following result is proved in \cite[Proposition 5.14]{Be1},
which shows that, under certain  natural assumptions, torsion pairs
in compactly generated triangulated categories can be lifted to TTF
triples.

\begin{Lem} Let $\mathcal{C}$ be a compactly generated triangulated
category which admits all small coproducts and products. Suppose
that $(\mathcal{Y},\mathcal{Z})$ is a hereditary torsion pair in
$\mathcal{C}$. Then we have the following.

$(1)$ If $\mathcal{Y}$ is closed under all small products, then
there exists a {\rm{TTF}} triple
$(\mathcal{X},\mathcal{Y},\mathcal{Z})$ in $\mathcal{C}$. In this
case, $\mathcal{Y}$ is compactly generated.

$(2)$ If $\mathcal{Z}$ is closed under all small coproducts, then
there exists a {\rm{TTF}} triple
$(\mathcal{Y},\mathcal{Z},\mathcal{W})$ in $\mathcal{C}$. In this
case, $\mathcal{Z}$ is compactly generated. \label{lem02}
\end{Lem}

The relationship between compact objects and TTF triples is
explained in the next result, which states that any set of compact
objects in a triangulated category with small coproducts gives rise
to a TTF triple. For more details, we refer the reader to
\cite[Chapter III, Theorem 2.3; Chapter IV, Proposition 1.1]{BI}.

\begin{Lem}
Let $\mathcal{C}$ be a triangulated category which admits all small
coproducts.  Suppose that  $\mathcal{P}$ is  a set of compact
objects in $\mathcal{C}$. Set
$\mathcal{X}:={\rm{Tria}}(\mathcal{P}),
\mathcal{Y}:=\Ker(\Hom_{\mathcal{C}}(\mathcal{X},-))$ and
$\mathcal{Z}:=\Ker(\Hom_{\mathcal{C}}(\mathcal{Y},-))$. Then
$(\mathcal{X},\mathcal{Y},\mathcal{Z})$ is a {\rm{TTF}} triple in
$\mathcal{C}$. Moreover,  $\mathcal{Y}$ coincides with the full
subcategory of $\mathcal{C}$ consisting of the objects $Y$ such that
$\Hom_{\mathcal{C}}(P[n],Y)=0$ for every $P\in\mathcal{P}$ and
$n\in\mathbb{Z}$. \label{lem03}
\end{Lem}

\section{Universal localizations and recollements\label{sect3}}

In this section, we discuss the connection between universal
localizations and recollements of triangulated categories. In our
considerations, homological ring epimorphisms and perpendicular
categories will play a role.

Now, we fix a ring $R$, and suppose that $\Sigma$ is a set of
 homomorphisms between finitely generated projective $R$-modules.
For each $f:P^{-1}\to P^0$ in $\Sigma$, we  denote by $\cpx{P_f}$
the following complex of $R$-modules:
$$
\cdots\lra 0\lra P^{-1}\lraf{f} P^{0}\lra 0\lra \cdots,
$$
where $P^{-1}$ and $P^0$ are concentrated in the degrees $-1$ and
$0$, respectively.

Set $$ \cpx{\Sigma}:=\{\cpx{P_f}~|~f\in\Sigma\},
$$
$$
{\Sigma^{\perp}}:=\{X\in R\Modcat~|~ \Hom_{\D{R}}(\cpx{P}, X[i])=0
\mbox{ for\, all \,} \cpx{P}\in\cpx{\Sigma}\mbox{ and\,all
\,}{i\in\mathbb{Z}}\},
$$
$$ \mathscr{D}(R)_{\Sigma^{\perp}}:=\{\cpx{Y}\in\D{R}~|~H^n(\cpx{Y})\in\Sigma^{\perp}
\mbox{ for\, all \,}n\in\mathbb{Z}\},
$$
where $H^n(\cpx{Y})$ is the $n$-th cohomology of the complex
$\cpx{Y}$. Note that some special cases of $\Sigma^{\perp}$ have
been discussed in literature (see, for example, \cite{HA, HKL, CTT,
GL}). For example, the set $\Sigma$ consists of injective
homomorphisms or only one single homomorphism. In those papers, such
a category $\Sigma^{\perp}$ is called the perpendicular category of
$\Sigma$.

Universal localizations were pioneered by Ore and Cohn, in order to
study embeddings of noncommutative rings in skewfields.

Before recalling the definition of universal localizations, we
mention the following result, due initially to Cohen (see also
\cite{Sch}), which explains how universal localizations arise.

\begin{Lem}{\rm \cite{cohenbook1}} Let $R$ and $\Sigma$ be as above.
Then there is a ring $R_{\Sigma}$ and a homomorphism  $\lambda:R\to
R_{\Sigma}$ of rings  with the following properties:

$(1)$ $\lambda$ is $\Sigma$-inverting, that is, if $\alpha:P\to Q$
belongs to $\Sigma$, then
$R_{\Sigma}\otimes_R\alpha:R_{\Sigma}\otimes_RP\to
R_{\Sigma}\otimes_R Q$ is an isomorphism of $R_{\Sigma}$-modules,
and

$(2)$ $\lambda$ is universal $\Sigma$-inverting, that is, if $S$ is
a ring such that there exists a $\Sigma$-inverting homomorphism
$\varphi:R\to S$, then there exists a unique homomorphism
$\psi:R_{\Sigma}\to S$ of rings such that $\varphi=\lambda\psi$.
\label{lem32}
\end{Lem}

The homomorphism $\lambda:R\to R_{\Sigma}$ in Lemma \ref{lem32} is a
ring epimorphism with $\Tor^R_1(R_{\Sigma}, R_{\Sigma})=0.$ It is
called the universal localization of $R$ at $\Sigma$. It is easy to
see that if $R$ has weak dimension at most $1$, then the
localization $\lambda: R\ra R_{\Sigma}$ of $R$ at any set $\Sigma$
is homological, and moreover, the weak dimension of $R_{\Sigma}$ is
also at most $1$ by Lemma \ref{2.1}.

If $\Sigma$ is a finite set, then we may assume that $\Sigma$
contains only one homomorphism since the universal localization at
$\Sigma$ is the same as the universal localization at the direct sum
of the homomorphisms in $\Sigma$.

The following result is a general formulation of the case discussed
in \cite{HA} and \cite{HKL}. Nevertheless, many arguments of the
proof there work in this general situation. We outline here a
modified proof.

\begin{Prop}
$(1)$ $\Sigma^{\perp}$ is closed under isomorphic images,
extensions, kernels, cokernels, direct sums and products.

$(2)$ $\Sigma^{\perp}$ coincides with the image of the restriction
functor $\lambda_*: {R_{\Sigma}}\Modcat\to R\Modcat$ induced by the
ring homomorphism $\lambda$ defined in Lemma \ref{lem32}. In this
sense, we can identify $ \Sigma^{\perp}$ with ${R_{\Sigma}}\Modcat$
via the homomorphism $\lambda$.

$(3)$ $\mathscr{D}(R)_{\Sigma^{\perp}}=
\Ker(\Hom_{\D{R}}({\rm{Tria}}(\cpx{\Sigma}), -))$. \label{prop31}
\end{Prop}
In order to prove Proposition \ref{prop31}, we need the following
known homological result.
\begin{Lem}
Suppose that $\cpx{W}=(W^i)_{i\in\mathbb{Z}}$ is a  complex in
${\mathscr C}(\Pmodcat{A})$ such that $W^n=0$ for all
$n\in\mathbb{Z}\backslash\{-1,0\}$. Then, for each $\cpx{X}\in\D{R}$
and $n\in\mathbb{Z}$, there is an exact sequence of abelian groups:
$$
0\lra\Hom_{\D{R}}(\cpx{W},H^{n-1}(\cpx{X})[1])\lra
\Hom_{\D{R}}(\cpx{W},\cpx{X}[n])\lra\Hom_{\D{R}}(\cpx{W},H^{n}(\cpx{X}))\lra
0.$$ \label{lem31}
\end{Lem}

{\it Proof.} It is sufficient to show the statement for $n=0$. In
this case, it follows from the triangle $ W^{-1}\ra W^0\ra
\cpx{W}\ra W^{-1}[1]$ that the following diagram is commutative and
exact: {\tiny
$$\begin{CD}
\Hom_{\K{R}}(W^0[1],\cpx{X})@>>>\Hom_{\K{R}}(W^{-1}[1],\cpx{X})@>>>\Hom_{\K{R}}(\cpx{W},\cpx{X})@>>>\Hom_{\K{R}}(W^0,\cpx{X})@>>>\Hom_{\K{R}}(W^{-1},\cpx{X})\\
@VV{\simeq}V @VV{\simeq}V @. @V{\simeq}VV @V{\simeq}VV \\
\Hom_{\K{R}}\big(W^0,
H^{-1}(\cpx{X})\big)@>>>\Hom_{\K{R}}\big(W^{-1},H^{-1}(\cpx{X})\big)@.
@.\Hom_{\K{R}}\big(W^0,H^{0}(\cpx{X})\big)@>>>\Hom_{\K{R}}\big(W^{-1},H^{0}(\cpx{X})\big).
\end{CD} $$}
Here we use the fact that $\Hom_{\D{R}}(P,
\cpx{X}[n])=\Hom_{\K{R}}(P, \cpx{X}[n])\simeq
\Hom_R(P,H^n(\cpx{X}))$ for every projective module $P$ and $n\in
{\mathbb Z}$. Thus Lemma \ref{lem31} follows.

\medskip
{\bf Proof of Proposition \ref{prop31}}. $(1)$ Clearly,
$\Sigma^{\perp}$ is closed under isomorphic images and extensions.
In the following, we shall prove that $\Sigma^{\perp}$ is closed
under kernels and cokernels. Recall that $\Sigma^{\perp}$ is defined
to be the full subcategory of $R\Modcat$ consisting of those
$R$-modules $X$ that
$\Hom_{\D{R}}(\cpx{U},X)=\Hom_{\D{R}}(\cpx{U},X[1])=0$ for all
$\cpx{U}\in\cpx{\Sigma}$. Suppose that  $f:Y\to Z$ is  a
homomorphism between two modules $Y$ and $Z$ in $\Sigma^{\perp}$.
Set $K:=\Ker(f),I:=\Img(f)$ and $C:=\Coker(f)$. Then we have two
exact sequences of $R$-modules:
$$0\to K\to Y\to I\to 0\quad \mbox{ and }\quad 0\to I\to Z\to C\to 0.$$
Since any short exact sequence in $R\Modcat$ can canonically be
extended to a triangle in $\D{R}$, we get two triangles in $\D{R}$:
$$ K\to Y\to I\to K[1]\quad \mbox{ and }\quad  I\to Z\to C\to I[1].$$
For convenience, we will write ${_{\D{R}}}(\cpx{X},\cpx{Y})$  for
the Hom-set $\Hom_{\D{R}}(\cpx{X},\cpx{Y})$ for
$\cpx{X},\cpx{Y}\in\D{R}$. Let $\cpx{P}\in\cpx{\Sigma}$. Then, by
applying ${_{\D{R}}}(\cpx{P},-)$ to these triangles, we obtain two
long exact sequences of abelian groups
$$
0\to{_{\D{R}}}(\cpx{P},K)\to{_{\D{R}}}(\cpx{P},Y)\to{_{\D{R}}}(\cpx{P},I)\to
{_{\D{R}}}(\cpx{P},K[1])\to {_{\D{R}}}(\cpx{P},Y[1]) \to
{_{\D{R}}}(\cpx{P},I[1]\to 0;
$$
$$
0\to{_{\D{R}}}(\cpx{P},I)\to{_{\D{R}}}(\cpx{P},Z)\to{_{\D{R}}}(\cpx{P},C)\to
{_{\D{R}}}(\cpx{P},I[1])\to {_{\D{R}}}(\cpx{P},Z[1]) \to
{_{\D{R}}}(\cpx{P},C[1]\to 0.
$$
Since $Y$ and $Z$ lie in $\Sigma^{\perp}$, we know
${_{\D{R}}}(\cpx{P},Y)={_{\D{R}}}(\cpx{P},Z)={_{\D{R}}}(\cpx{P},Y[1])={_{\D{R}}}(\cpx{P},Z[1])=0.$
It follows that ${_{\D{R}}}(\cpx{P},K)={_{\D{R}}}(\cpx{P},I)=0$, and
so ${_{\D{R}}}(\cpx{P},K[1])=0$. This implies $K\in\Sigma^{\perp}$.
Similarly, we can conclude that $I$ and $C$ belong to
$\Sigma^{\perp}$. Hence  $\Sigma^{\perp}$ is closed under kernels,
images and cokernels. By the definition of $\Sigma^{\perp}$ and the
fact that Hom-funcors commute with products, we infer that
$\Sigma^{\perp}$ is closed under products. Since $\cpx{\Sigma}$ is a
set of bounded complexes over finitely generated projective
$R$-modules, these complexes are compact, and therefore
$\Sigma^{\perp}$ is closed under direct sums.

$(2)$ Observe that, for each element $f:P^{-1}\to P^0$ in $\Sigma$,
there is a canonical triangle in $\D{R}$:
$$
(\ast)\quad P^{-1}\lraf{f} P^0\lra\cpx{P_f}\lra P^{-1}[1].
$$
If, in addition, $f$ is injective, then we have a short exact
sequence of $R$-modules:
$$
(\ast\ast)\quad 0\lra P^{-1}\lraf{f} P^0\lra\Coker(f)\lra 0.
$$
In this case, we get $\cpx{P_f}\simeq\Coker(f)$ in $\D{R}$. Note
that the same statement as $(2)$ is obtained in \cite[Lemma 1.6,
Proposition 1.7]{HA} under the extra assumption that each element in
$\Sigma$ is injective, where the sequence $(\ast\ast)$ is used. In
fact, this assumption is not necessary since we can replace ($**$)
by $(\ast)$ and modify the proof there to show the general case. For
more details, we refer the reader to \cite{HA}.

$(3)$ This follows directly from Lemma \ref{lem03} and Lemma
\ref{lem31}. $\square$

\medskip
Combining Lemma \ref{2.0} with Proposition \ref{prop31}, we have the
following result, which says that, in some sense, Morita equivalence
preserves universal localization.

\begin{Koro} Let $\lambda:R\to R_{\Sigma}$ be the universal
localization of the ring $R$ at the set $\Sigma$. Suppose that  $P$
is a finitely generated projective generator for $R\Modcat$. Set
$\Delta:=\{\Hom_R(P,f)\mid f\in\Sigma\}$. Then  the ring
homomorphism
$\mu:\End_R(P)\to\End_{R_{\Sigma}}(R_{\Sigma}\otimes_RP)$, defined
by $g\mapsto R_{\Sigma}\otimes_Rg$ for any $g\in\End_R(P)$, is the
universal localization of the ring $\End_R(P)$  at the set $\Delta$.
\label{Morita}
\end{Koro}

{\it Proof.} Let $S:=\End_R(P)$. Since $_RP$ is a finitely generated
projective generator for $R\Modcat$, the Hom-functor
$\Hom_R(P,-):R\Modcat\to S\Modcat$ is an equivalence, which extends
to a triangle equivalence  between $\D{R}$ and $\D{S}$. By the
definitions of ${\Sigma^{\perp}}$ and ${\Delta^{\perp}}$, the
restriction of $\Hom_R(P,-)$ induces an equivalence from
${\Sigma^{\perp}}$ to ${\Delta^{\perp}}$. Note that
$R_{\Sigma}\otimes_RP$ is a finitely generated projective generator
for $R_{\Sigma}\Modcat$.  Since the functor $\lambda_*:
R_{\Sigma}\Modcat\ra R\Modcat$ is fully faithful and since the image
of $\lambda_*$ coincides with ${\Sigma^{\perp}}$ by Proposition
\ref{prop31}(2), it follows from the following commutative diagram
of functors:
$$\xymatrix{
R_{\Sigma}\Modcat\ar^-{\Hom_{R_{\Sigma}}(R_{\Sigma}\otimes_RP,-)}_-{\simeq}[rrr]\ar_-{\lambda_*}[d]
&&&\End_{R_{\Sigma}}(R_{\Sigma}\otimes_RP)\Modcat\ar^{\mu_{*}}[d] \\
R\Modcat\ar^-{\Hom_R(P,-)}_-{\simeq}[rrr] &&&S\Modcat}
$$
that $\mu_*$ is fully faithful, and that the image of $\mu_*$
coincides with ${\Delta^{\perp}}$. This implies also that $\mu$ is a
ring epimorphism.  Note that, under our conventions, full
subcategories are always closed under isomorphic images.

On the other hand, if $\varphi: S\ra S_{\Delta}$ is the universal
localization of $S$ at $\Delta$, then, by Proposition
\ref{prop31}(2), the image of $\varphi_*$ coincides with
$\Delta^{\perp}$. Thus the two ring epimorphisms $\mu$ and $\varphi$
are equivalent by Lemma \ref{2.0}. This means that the two rings
$S_{\Delta}$ and $\End_{R_{\Sigma}}(R_{\Sigma}\otimes_RP)$ are
isomorphic. Thus $\mu$ is the universal localization of $S$ at
$\Delta$. $\square$

\medskip
Motivated by \cite[Theorem 10.8]{nr}, see also \cite[Thereom 4.8
(3)]{HKL}, we shall establish the following connection between
universal localizations and recollements of triangulated categories.
The last condition (5) of Proposition \ref{th2} below seems to
appear for the first time in the work, and will be used in Section
\ref{sect4} to prove Theorem \ref{th01}(1).

\begin{Prop} $(a)$ Let ${\bf{j}}$ be the canonical embedding of $\mathscr{D}(R)_{\Sigma^{\perp}}$ into $\mathscr{D}(R)$. Then there is a
recollement
$$\xymatrix@C=1.2cm{\mathscr{D}(R){_{\Sigma^{\perp}}}\ar[r]^-{{\bf{j}}}
&\D{R}\ar[r]\ar@/^1.2pc/[l]\ar_-{\bf{L}}@/_1.2pc/[l]
&{\rm{Tria}}{(\cpx{\Sigma})} \ar@/^1.2pc/[l]\ar@/_1.2pc/[l]}$$

\medskip
\noindent such that $\bf{L}$ is the left adjoint of ${\bf j}$ and
$\cpx{T}:={\bf{L}}(R)$ is a compact generator of
$\mathscr{D}(R){_{\Sigma^{\perp}}}$.

$(b)$ The following statements are equivalent:

\;\quad $(1)$ $\lambda: R\to  R_{\Sigma}$ is a homological
epimorphism of rings;

\quad $(2)$ $\lambda_{*}: \D{R_{\Sigma}}\xrightarrow {\sim}
\mathscr{D}(R){_{\Sigma^{\perp}}}$;

\quad $(3)$ the complex $\cpx{T}$ in (a) is a tilting object in
$\mathscr{D}(R){_{\Sigma^{\perp}}}$;

\quad $(4)$ the complex $\cpx{T}$ in (a) is isomorphic to
$R_{\Sigma}$ in $\D{R}$;

\quad $(5)$ the complex $\cpx{T}$ in (a) is isomorphic in $\D{R}$ to
a complex $\cpx{X}:=(X^i)_{i\in\mathbb{Z}}$ such that $X^i\in
\Sigma^{\perp}$ for all $i\in\mathbb{Z}$.

\label{th2}
\end{Prop}

{\bf Proof.}  The existence of the above recollement is an immediate
consequence of  Lemmata \ref{lem01}(2), \ref{lem03} and Proposition
\ref{prop31}. The property in $(a)$ follows from the proof in
\cite[Chapter IV, Proposition 1.1]{BI}. As to the property $(b)$, we
notice that the equivalences among the first four statements in (b)
can be deduced from \cite[Proposition 1.7, Lemma 4.6]{HKL}. Clearly,
the statement $(4)$ implies the statement $(5)$. We shall show that
(5) implies (4).

Let $\lambda:R\to R_{\Sigma}$ be the universal localization of $R$
at $\Sigma$. In what follows,  we always  identify $\Sigma^{\perp}$
with $R_{\Sigma}\Modcat$ via $\lambda$. This is due to Proposition
\ref{prop31}(2).

Suppose that $\cpx{T}\simeq \cpx{X}:=(X^i)_{i\in\mathbb{Z}}$ in
$\D{R}$ such that $X^i\in R_{\Sigma}\Modcat$ for all
$i\in\mathbb{Z}$. Since $\lambda$ is a ring epimorphism, we get
$\Hom_{R_{\Sigma}}(X,Y)\simeq\Hom_R(X,Y)$ for all $X,Y\in
{R_{\Sigma}}\Modcat$. Thus $\cpx{X}$ can be considered as a complex
over $R_{\Sigma}\Modcat$, that is, $\cpx{X}\in\C{R_{\Sigma}}$. Let
$\lambda_1$ be the map $\Hom_{\D{R}}(\lambda, \cpx{X}):
\Hom_{\D{R}}(R_{\Sigma}, \cpx{X})\to\Hom_{\D{R}}(R,\cpx{X})$. We
claim that $\lambda_1$ is surjective. In fact, there is a
commutative diagram:
$$\xymatrix{\Hom{_{\K{R}}}(R_{\Sigma}, \cpx{X})\ar[r]^-{q_1}\ar[d]^-{\lambda_2}&
\Hom_{\D{R}}(R_{\Sigma}, \cpx{X})\ar[d]^-{\lambda_1}\\
\Hom_{\K{R}}(R, \cpx{X})\ar[r]^-{q_2}&\Hom_{\D{R}}(R, \cpx{X}),}$$
where $\lambda_2=\Hom_{\K{R}}(\lambda, \cpx{X})$, and $q_1,q_2$ are
induced by the localization functor $q:\K{R}\to\D{R}$. Clearly,
$q_2$ is a bijection. To prove that $\lambda_1$ is surjective, it
suffices to show that $\lambda_2$ is bijective. Indeed, $\lambda_2$
is a composition of the following series of isomorphisms:
$$\Hom_{\K{R}}(R_{\Sigma}, \cpx{X})\simeq H^0(\Hom_R(R_{\Sigma},
\cpx{X}))= H^0(\Hom_{R_{\Sigma}}(R_{\Sigma}, \cpx{X}))\simeq
\Hom_{\K{R}}(R, \cpx{X}),
$$
where the equality follows from the fact that $\lambda$ is a ring
epimorphism. More precisely, for $\bar{\cpx{f}}:=\overline{(f^i)}\in
\Hom_{\K{R}}(R_{\Sigma},\cpx{X})$ with $(f^i)_{i\in {\mathbb Z}}$ a
chain map, the series of the above maps are defined by:
$$ \overline{(f^i)}\mapsto \overline{f^0} = \overline{f^0}\mapsto \overline{\lambda_*(\cpx{f})},$$
where $\lambda_*(\cpx{f})$ is a chain map from $R$ to $\cpx{X}$ with
$\lambda f^0$ in degree $0$ and zero in all other degrees. Thus
$\lambda_2$ is bijective, which implies that $\lambda_1$ is
surjective. Now, let $\lambda'$ be the map $\Hom_{\D{R}}(\lambda,
\cpx{T}):\Hom_{\D{R}}(R_{\Sigma},
\cpx{T})\to\Hom_{\D{R}}(R,\cpx{T})$.  Since $\cpx{T}\simeq \cpx{X}$
in $\D{R}$, we know that  $\lambda'$ also is surjective. Suppose
that  $\varphi:R\to \cpx{T}:={\bf{L}}(R)$ is the unit adjunction
morphism with respect to the adjoint pair $(\bf{L},\bf{j})$. Then
there exists $g:R_{\Sigma}\to\cpx{T}$ in $\D{R}$ such that
$\varphi=\lambda g$. Since $R_{\Sigma}$ belongs to $
\mathscr{D}(R)_{\Sigma^{\perp}}$, there exists $f:\cpx{T}\to
R_{\Sigma}$ in $\D{R}$ such that $\lambda=\varphi f$. This gives
rise to the following commutative diagram in $\D{R}$:
$$
\xymatrix{
R\ar@{=}[r]\ar[d]_-{\varphi}&R\ar@{=}[r]\ar[d]_-{\lambda}&R \ar[d]^-{\varphi}\\
\cpx{T}\ar@{-->}[r]^-{f}&R_{\Sigma}\ar@{-->}[r]^-{g}&\;\cpx{T}. }
$$
Consequently, $\varphi=\varphi fg$ and $\lambda=\lambda gf$. On the
one hand, since $\varphi$ is the unit adjunction morphism, we have
$fg=1{_{\cpx{T}}}$. On the other hand, it follows from \cite[Theorem
1.4]{HA} that $\lambda$ is an $R_{\Sigma}\Modcat$-reflection of $R$,
that is, the morphism of abelian groups
$\Hom_R(\lambda,Z):\Hom_R(R_{\Sigma},Z)\to\Hom_R(R,Z)$  is bijective
for any $Z\in R_{\Sigma}\Modcat$. This yields $gf=1_{R_{\,\Sigma}}$.
Thus $f$ is an isomorphism. In other words, $\cpx{T}\simeq
R_{\Sigma}$ in $\D{R}$. Therefore, (5) implies (4). $\square$

\bigskip
{\it Remark}. Note that every tilting module is associated to a
class of finitely presented modules of projective dimension at most
one (see \cite{HA, BZ2}) and thus to a universal localization since
each finitely presented module of projective dimension at most one
is determined by an injective homomorphism between finitely
generated projective modules. In Proposition \ref{th2}, we do not
require that each homomorphism in $\Sigma$ is injective. From this
point of view, the formulation of Proposition \ref{th2}(b) seems to
be more general than that in \cite[Thereom 4.8(3)]{HKL}.

\medskip
\begin{Koro} Let $R\subseteq S$ be an extension of rings, that is, $R$ is a subring of the ring $S$ with the same
identity, and let $B$ be the endomorphism ring of the $R$-module
$S\oplus S/R$. Then there is a recollement of triangulated
categories:

$$\xymatrix@C=1.2cm{\D{B}_{\Sigma^{\perp}}\ar[r]&\D{B}\ar[r]
\ar@/^1.2pc/[l]\ar@/_1.2pc/[l]
&\D{R}\ar@/^1.2pc/[l]\ar@/_1.2pc/[l]},\vspace{0.3cm}$$ where
$\Sigma: =\{\pi^*\}$, and the homomorphism  $\pi^*: \Hom_R(S\oplus
S/R,S)\to\Hom_R(S\oplus S/R,S/R)$ of $B$-modules is defined by
$f\mapsto f\pi $ for any $f\in\Hom_R(S\oplus S/R,S)$, which is
induced by the canonical map $\pi: S\ra S/R$. \label{3.5}
\end{Koro}

{\bf Proof.} It follows from Proposition \ref{th2}($a$) that we have
the following recollement:

$$\xymatrix@C=1.2cm{\D{B}_{\Sigma^{\perp}}\ar[r]&\D{B}\ar[r]
\ar@/^1.2pc/[l]\ar@/_1.2pc/[l] & {\rm
Tria}(\cpx{\Sigma})\ar@/^1.2pc/[l]\ar@/_1.2pc/[l]}.\vspace{0.3cm}$$

To show that ${\rm Tria}(\cpx{\Sigma})$ is equivalent to $\D{R}$ as
triangulated categories, it suffices to prove that the complex
$\cpx{\Sigma}\in \Kb{\pmodcat{B}}$ is exceptional with
$\End_{\D{B}}(\cpx{\Sigma})\simeq R$.

In fact, let $_RT:=S\oplus S/R$ and $\End_R(T) = B$. Then
$\add(_RT)$ and $\pmodcat{B}$ are equivalent, and therefore
$\Kb{\add(_RT)}$ and $\Kb{\pmodcat{B}}$ are equivalent as
triangulated categories via the functor $\Hom_R(T,-)$. Thus, to show
that the complex $\cpx{\Sigma}\in \Kb{\pmodcat{B}}$ is exceptional
with $\End_{\D{B}}(\cpx{\Sigma})\simeq R$, it is sufficient to show
that the complex
$$ \cpx{\Pi}:\quad 0\lra S\lraf{\pi} S/R\lra 0 $$
in $\Kb{\add(T)}$ is exceptional with
$\End_{\Kb{\add(T)}}(\cpx{\Pi})\simeq R$ since
$\Hom_R(T,\cpx{\Pi})=\cpx{\Sigma}$.

It is easy to see $\Hom_{\Kb{\add(T)}}(\cpx{\Pi},\cpx{\Pi}[-1]) =
0$. To see $\Hom_{\Kb{\add(T)}}(\cpx{\Pi},\cpx{\Pi}[1])=0$, we pick
up a homomorphism $f: S\ra S/R$ of $R$-modules, suppose $(1)f= s+
R\in S/R$ and define $g: {}_RS\ra {}_RS$ by $x\mapsto xs$ for $x\in
S$. Clearly, $g$ is a homomorphism of $R$-modules and $(f-g)|_R=0$.
Thus there exists a homomorphism $h: S/R\ra S/R$ such that $f-g=\pi
h$. This implies that $f$ is zero in $\Kb{\add(T)}$, that is,
$\Hom_{\Kb{\add(T)}}(\cpx{\Pi},\cpx{\Pi}[1])=0$. Hence we have shown
that $\cpx{\Pi}$ is exceptional.

Now, we define a ring homomorphism $\alpha$ from
$\End_{\Kb{\add(T)}}(\cpx{\Pi})$ to $R$ as follows: Given
$f=(f^0,f^1)\in \End_{\Kb{\add(T)}}(\cpx{\Pi})$, let $(f)\alpha$ be
the unique map determined in the following diagram of $R$-modules:
$$\begin{CD}
0 @>>> R@>{\lambda}>> S@>{\pi}>> S/R@>>> 0\\
@. @V{(f)\alpha}VV  @V{f^0}VV  @V{f^1}VV  \\
0 @>>> R@>{\lambda}>> S@>{\pi}>> S/R@>>> 0.
\end{CD} $$
Note that if $f$ is null-homotopic then $(f)\alpha$ is zero. This
means that $\alpha$ is well-defined. Clearly, $\alpha$ is a ring
homomorphism. We claim that $\alpha$ is an isomorphism of rings. It
is easy to check that $\alpha$ is injective. We shall show that
$\alpha$ is surjective. Let $r\in R$. We define $f^0: S\ra S$ to be
the right multiplication of $r$. Then there is a homomorphism $f^1:
S/R\ra S/R$ of $R$-modules such that $f^0\pi=\pi f^1$. This means
that $\alpha$ is surjective. Hence $\alpha$ is an isomorphism of
rings. So, $\cpx{\Sigma}$ is exceptional with
$\End_{\D{B}}(\cpx{\Sigma})\simeq R$. By Lemma \ref{compact-deriv},
we may identify Tria$(\cpx{\Sigma})$ with $\D{R}$. This proves
Corollary \ref{3.5}. $\square$

\medskip
As another corollary of Proposition \ref{th2}, we have the following
result.

\begin{Koro} If the weak dimension of $R$ is at most $1$, then
there is a recollement

$$\xymatrix@C=1.2cm{\mathscr{D}(R_{\Sigma})\ar[r]
&\D{R}\ar[r]\ar@/^1.2pc/[l]\ar @/_1.2pc/[l]
&{\rm{Tria}}{(\cpx{\Sigma})}
\ar@/^1.2pc/[l]\ar@/_1.2pc/[l]},\vspace{0.2cm}$$ where $\Sigma$ is a
set of homomorphisms between finitely generated projective
$R$-modules. \label{35}
\end{Koro}

{\it Proof.} Under the assumption, the universal localization map
$\lambda_{\Sigma}$ is trivially a homological ring epimorphism. So,
this corollary follows from Proposition \ref{th2}(b). $\square$

\medskip
As a consequence of Corollary \ref{35}, we have the following
result which is a generalization of \cite[Theorem 2.5, Corollary
3.3]{akl}.

\begin{Koro} Suppose that $R$ is a left semi-hereditary ring (that is, every finitely
generated submodule of a projective left $R$-module is projective).
If $\cpx{T}$ is a compact exceptional object in $\D{R}$, then there
is a ring $S$, a homological ring epimorphism $\lambda: R\ra S$ and
a recollement

$$\xymatrix@C=1.2cm{\D{S}\ar[r]
&\D{R}\ar[r]\ar@/^1.2pc/[l]\ar @/_1.2pc/[l] &
\D{\End_{\D{R}}(\cpx{T})} \ar@/^1.2pc/[l]\ar@/_1.2pc/[l]}.$$
\label{36}
\end{Koro}

{\it Proof.}  Since $\cpx{T}$ is a compact object in $\D{R}$, there
exists  a complex $\cpx{P}\in \Kb{\pmodcat{R}}$ such that
$\cpx{T}\simeq \cpx{P}$ in $\D{R}$. Suppose that $\cpx{P}$ is of the
following form
$$
\cdots\lra 0\lra P^s\lra\cdots\lra P^i\lraf{d^i} P^{i+1}\lra
\cdots\lra P^t\lra 0\lra \cdots, $$ where $P^i\in \add(_RR)$ for
$s\leq i\leq t$. Since $R$ is left semi-hereditary, we have
Im$(d^i)\in \add(_RR)$ for all $i$. This implies that $\cpx{P}$ is
isomorphic to a direct sum of finitely many two-term complexes in
$\Kb{\pmodcat{R}}$, say $\cpx{P}\simeq \oplus_{i=1}^n\cpx{P_i}$,
where $n\in {\mathbb N}$ and $\cpx{P_j}$ is of the form: $0\ra
P_j^{s_j-1}\lraf{d_j} P_j^{s_j}\ra 0$ with $s_j\in {\mathbb Z}$ for
each $1\le j\le n$. Now we choose $\Sigma:=\{d_j\mid 1\le j\le n\}$.
Then, by definition, we have $\cpx{\Sigma} = \{\cpx{P_j}[s_j]\mid
1\le j\le n\}$ (see notations at the beginning of Section
\ref{sect3}). Consequently, Tria$(\cpx{T})$ = Tria$(\cpx{P})$ =
Tria$(\cpx{\Sigma})$. Since $\cpx{T}$ is compact and
$\Hom_{\D{R}}(\cpx{T},\cpx{T}[i])=0$ for all $i\ne 0$, it follows
from Lemma \ref{compact-deriv} that Tria$(\cpx{T})$ is triangle
equivalent to the derived category $\D{\End_{\D{R}}(\cpx{T})}$ of
$\End_{\D{R}}(\cpx{T})$. Note that $R$ has weak dimension at most
$1$ because it is left semi-hereditary. Thus  the condition of
Corollary \ref{35} is fulfilled, and therefore Corollary \ref{36}
follows from Corollary \ref{35} if we define $S=R_{\Sigma}$.
$\square$

\medskip
In general, it is hard to compute $R_{\Sigma}$. However, if $\Sigma$
consists of only one element with an orthogonal assumption, one can
construct $R_{\Sigma}$ explicitly in terms of endomorphism rings of
modules. To this purpose, we first establish the following result
which generalizes \cite[Proposition 1.3]{CTT} where only stalk
complexes (or modules) were considered.

\begin{Prop}
Suppose that $\cpx{P}:=(P^i)_{i\in\mathbb{Z}}$ is a  complex in
$\C{\Pmodcat{R}}$ such that $P^n=0$ for all
$n\in\mathbb{Z}\backslash\{-1,0\}$. Define
${\cpx{P}}^{\perp}:=\{X\in R\Modcat~|~ \Hom_{\D{R}}(\cpx{P}, X[i])=0
\mbox{ for\, all \,} {i\in\mathbb{Z}}\}.$ If
\,$\Hom_{\D{R}}(\cpx{P}, \cpx{P}{^{(\delta)}}\,[1])=0$ for each
cardinal $\delta$, then the inclusion $j:{\cpx{P}}^{\perp} \to
R\Modcat$ admits a left adjoint $l: R\Modcat\to{\cpx{P}}^{\perp}$.

\label{prop32}
\end{Prop}

{\bf Proof.} The proof will be divided into three steps. We define
$\mathscr{X}$ to be the full subcategory of $R\Modcat$ consisting of
the objects $X$ such that $\Hom_{\D{R}}(\cpx{P}, X[1])=0$. Then, it
follows from $ \Hom_{\D{R}}(\cpx{P}, X[1])\simeq
\Hom_{\K{R}}(\cpx{P}, X[1])$ for $X\in R\Modcat$ that $\mathscr{X}$
is closed under quotients.

Let $M$ and $N$ be $R$-modules.

\textbf{Step $(1)$}. For the given $M$,  we shall  construct an
$R$-module, denoted by $l(M)$, which belongs to ${\cpx{P}}^{\perp}$
and is endowed with an $R$-homomorphism $\eta_M: M\to l(M)$.

Since $\cpx{P}\in\Cb{\Pmodcat{R}}$, we have $ \Hom_{\D{R}}(\cpx{P},
M[1])= \Hom_{\K{R}}(\cpx{P}, M[1])$. Let $\alpha$ be a generating
set of  $ \Hom_{\K{R}}(\cpx{P}, M[1])$ as an
$\End{_{\D{R}}}(\cpx{P})$-module. Thus each element of $\alpha$ is a
chain map from $\cpx{P}$ to $M[1]$. We define $\omega_M:
\cpx{P}{^{(\alpha)}}\to M[1]$ to be the coproduct of the elements of
$\alpha$. Then, it is clear that
$$
\omega_M':=\Hom_{\D{R}}(\cpx{P},\omega_M):\Hom_{\D{R}}(\cpx{P},
\cpx{P}{^{(\alpha)}})\lra \Hom_{\D{R}}(\cpx{P}, M[1])
$$
is a surjective homomorphism of $\End{_{\D{R}}}(\cpx{P})$-modules.
Let $\rm{cone}(\omega_M)$ be the mapping cone of $\omega_M$, and
$\cpx{\ol{M}}:=\rm{cone}(\omega_M)[-1]$.  Then we obtain the
following canonical triangle in $\D{R}$:
$$\begin{CD}
(\ast)\qquad
M@>{\varphi_M}>>\cpx{\ol{M}}@>{\psi_M}>>{\cpx{P}}^{(\alpha)}@>{\omega_M}>>M[1],
\end{CD}$$ where  $\varphi_M$ and $\psi_M$ can be constructed explicitly.
For more details, we refer the reader to any standard textbook of
homological algebra (for instance, \cite{We}). Note that the complex
$\cpx{\ol{M}}$ is of the form
$$ \cpx{\overline{M}}: \quad
 \,\cdots \lra 0\lra M^{-1}\xrightarrow{d_M} M^0\lra 0\lra
\cdots,
$$
where  $d_M$ is a homomorphism between $R$-modules $M^{-1}$ and
$M^0$,
 which are concentrated in the degrees $-1$ and $0$,
respectively. Let $\cpx{C}$ denote the complex: $0\ra \Img(d_M)\ra
M^0\ra 0$. Then we have an exact sequence of complexes: $0\ra
H^{-1}(\cpx{\ol{M}})[1]\ra \cpx{\ol{M}}\ra \cpx{C}\ra 0$, this gives
us an triangle
$$ H^{-1}(\cpx{\ol{M}})[1]\ra \cpx{\ol{M}}\ra \cpx{C}\ra H^{-1}(\cpx{\ol{M}})[2]. $$
Since $\cpx{C}$ is quasi-isomorphic to the stalk complex
$H^0(\cpx{\ol{M}})$, we get the following triangle in $\D{R}$:
$$
(\ast\ast)\quad H^{-1}(\cpx{\ol{M}})[1]\lra
\cpx{\ol{M}}\epa{\gamma_M} H^0(\cpx{\ol{M}})\lra
H^{-1}(\cpx{\ol{M}})[2],
$$
where the chain map $\gamma_M$ is induced by the homomorphism $d_M$
such that $H^{0}(\gamma_M)=1{_{H^0(\cpx{\ol{M}})}}$. Applying
$\Hom_{\D{R}}(P^\bullet, -)$ to $(\ast)$, we get a long exact
sequence of $\End{_{\D{R}}}(\cpx{P})$-modules:
$$
\Hom_{\D{R}}(\cpx{P}, \cpx{P}{^{(\alpha)}})\lraf{\omega_M'}
\Hom_{\D{R}}(\cpx{P}, M[1])\lra\Hom_{\D{R}}(\cpx{P}, \cpx{\ol{M}}[1]
)\lra \Hom_{\D{R}}(\cpx{P}, \cpx{P}{^{(\alpha)}}[1]).
$$
Since $\omega_M'$ is surjective and $\Hom_{\D{R}}(\cpx{P},
\cpx{P}{^{(\alpha)}}[1]) = 0$, we have $\Hom_{\D{R}}(\cpx{P},
\cpx{\ol{M}}[1]) = 0$.  Now it follows from Lemma \ref{lem31} that
$$\Hom_{\D{R}}(P^\bullet,H^0(\cpx{\ol{M}})[1])\simeq\Hom_{\D{R}}(\cpx{P},
\cpx{\ol{M}}[1])=0.$$ This implies
$H^0(\cpx{\ol{M}})\in\mathscr{X}$.

Now, we first fix a chain map $\omega_M$ for the given $M$, and then
define $l(M):=H^0(\cpx{\ol{M}})/t(M)$, where $t(M)$ denotes the
trace of $H^0(\cpx{P})$ in $H^0(\cpx{\ol{M}})$, that is, the sum of
the images  of all homomorphisms from $H^0(\cpx{P})$ to
$H^0(\cpx{\ol{M}})$. Thus $l(M)\in \mathscr X$ since $\mathscr X$ is
closed under quotients. By Lemma \ref{lem31}, we have
$\Hom_{\D{R}}(P^\bullet,H^0(P^\bullet)[1])$ $\simeq
\Hom_{\D{R}}(P^\bullet,P^\bullet[1]) =0$. This means
$H^0(P^\bullet)\in\mathscr{X}$. Since $H^0$ commutes with
coproducts, we infer from Lemma \ref{lem31} that  coproducts of
copies of $H^0(\cpx{P})$ lie in $\mathscr X$. This shows $t(M)\in
\mathscr X$ because it is an image of a coproduct of
$H^0(P^\bullet)$.

In the following, we shall prove $l(M)\in{\cpx{P}}^{\perp}$.
Clearly, we have $\Hom_{\D{R}}(P^\bullet,l(M)[i])=0$ for $i\not=0$,
and $\Hom_{\D{R}}(\cpx{P}, l(M)) \simeq \Hom_R(H^0(\cpx{P}), l(M))$.
So, in order to show $l(M)\in{\cpx{P}}^{\perp}$, it is sufficient to
prove $\Hom_R(H^0(\cpx{P}),l(M)) = 0$.

One the one hand, applying $\Hom_R(H^0(\cpx{P}),-)$ to the the
canonical exacts sequence
$$
0\lra t(M)\lra H^0(\cpx{\ol{M}})\lra  l(M)\lra 0
$$
$R$-modules, we can see that $\Hom_{R}(H^0(\cpx{P}), l(M))$ can be
embedded into $\Hom_{\D{R}}(H^0(\cpx{P}), t(M)[1])$ because
$\Hom_R(H^0(\cpx{P}), t(M))\lra\Hom_R(H^0(\cpx{P}),
H^0(\cpx{\ol{M}}))$ is always bijective by definition.  On the other
hand, applying $\Hom_{\D{R}}(-, t(M)[1])$ to the following canonical
triangle induced from $\cpx{P}$
$$
\quad H^{-1}(\cpx{P})[1]\lra \cpx{P}\lra H^0(\cpx{P})\lra
H^{-1}(\cpx{P})[2]
$$
in $\D{R}$ and using the fact that $\Hom_{\D{R}}(H^{-1}(\cpx{P})[1],
t(M)) = 0$ and that $t(M)\in \mathscr X$, we can deduce that
$\Hom_{\D{R}}(H^0(\cpx{P}), t(M)[1]) = 0$. Consequently,
$\Hom_{R}(H^0(\cpx{P}), l(M))=0$, as desired. Thus
$l(M)\in{\cpx{P}}^{\perp}$.

We define $\eta_M:=\varphi_M\gamma_M\pi_M:M\to l(M)$, which is
clearly a homomorphism of $R$-modules.

Similarly, for the module $N$, we fix, once and for all, a chain map
$\omega_N:\cpx{P}{^{(\beta)}}\to N[1]$, and then define $l(N)$ and
$\eta_N: N\to l(N)$, where $\beta$ is a cardinal. Clearly, we have
the following triangle in $\D{R}$:
$$
\xymatrix{N\ar[r]^-{\varphi_N}&\cpx{\ol{N}}\ar[r]^-{\psi_N}&\cpx{P}{^{(\beta)}}\ar[r]^-{\omega_N}&N[1]}
$$ with $\cpx{\ol{N}}:={\rm cone}(\omega_N)[-1]$.

\textbf{Step $(2)$}. For any  homomorphism $g: M\to N$ in
$R\Modcat$, we claim that there is a unique  homomorphism
$l(g):l(M)\to l(N)$ such that the following diagram commutes:
$$
\xymatrix{ M\ar[r]^-{\eta_M}\ar[d]_-{g}&l(M)\ar[d]^-{l(g)}\\
N\ar[r]^-{\eta_N}&l(N). }
$$
In fact, since $\Hom_{\D{R}}(P^\bullet, \cpx{\ol{N}}[1])=0$ by Step
(1), we know from homological algebra that
$\Hom_{\D{R}}(\cpx{P}{^{(\alpha)}}[-1]$, $\cpx{\ol{N}})$ $\simeq
\Pi_{\alpha}\Hom_{\D{R}}(\cpx{P}[-1],N)=0$. In particular, the
homomorphism $\omega_M[-1]g\varphi_N:\cpx{P}{^{(\alpha)}}[-1]\to
\cpx{\ol{N}}$ must be zero. Consequently, there is a homomorphism
$g':\cpx{\ol{M}}\to\cpx{\ol{N}}$ such that $g\varphi_N=\varphi_Mg'$.
So, we have the following commutative diagram
$$\xymatrix@C=1.2cm{\cpx{P}{^{(\alpha)}}[-1]\ar[r]^-{-\omega_M[-1]}&M\ar[r]^-{\varphi_M}\ar[d]^-{g}&
\cpx{\ol{M}}\ar[r]^-{\psi_M}\ar@{-->}[d]^-{g'}
&\cpx{P}{^{(\alpha)}}\\
P^{(\beta)}[-1]\ar[r]^-{-\omega_N[-1]}&N\ar[r]^-{\varphi_N}&\cpx{\ol{N}}\ar[r]^-{\psi_N}&
P^{(\beta)}.}$$
Notice that such $g'$ is not unique in general. Suppose that there
exists another morphism $h:\cpx{\ol{M}}\to\cpx{\ol{N}}$ such that
$g\varphi_N=\varphi_M h$. Then $g'-h=\psi_M s$  for some morphism
$s:\cpx{P}{^{(\alpha)}}\to\cpx{\ol{N}}$. Applying the 0-th
cohomological functor $H^0(-):\D{R}\to R\Modcat$ to this equality,
we get
$$H^0(g')-H^0(h)=H^0(\psi_M)H^0(s):\, H^0(\cpx{\ol{M}}) \lra
H^0(\cpx{\ol{N}}),$$ where $H^0(s):H^0(\cpx{P}{^{(\alpha)}})\to
H^0(\cpx{\ol{N}}).$ Since  the functor $H^0(-)$ commutes with
coproducts, we get  $H^0(\cpx{P}{^{(\alpha)}})\simeq
H^0(\cpx{P}){^{(\alpha)}}.$ This means that the image of
$H^0(g')-H^0(h)$ is contained in $t(N)$ by the definition of $t(N)$.
Thus
$$\ol{H^0(g')}=\ol{H^0(h)}:\, l(M) \lra
l(N),$$ which shows that $\ol{H^0(g')}$ depends on $g$ and not on
the choice of $g'$.  Thus, given $g: M\ra N$, we can define
$l(g):=\ol{H^0(g')}: l(M)\to l(N)$ which is a homomorphism  in
${\cpx{P}}^{\perp}$. To prove the equality $\eta_M l(g)=g\eta_N$, we
observe that ${H^0(g')}$ is the unique homomorphism from
$H^0(\cpx{\ol{M}})$ to $H^0(\cpx{\ol{N}})$ such that the following
diagram is commutative:
$$
\xymatrix{
 H^{-1}(\cpx{\ol{M}})[1]\ar[r]&\cpx{\ol{M}}\ar[r]^-{\gamma_M}\ar[d]_-{g'}
 &H^0(\cpx{\ol{M}})\ar[d]^-{H^0(g')}\ar[r]& H^{-1}(\cpx{\ol{M}})[2]\\
 H^{-1}(\cpx{\ol{N}})[1]\ar[r]&\cpx{\ol{N}}\ar[r]^-{\gamma_N}
 &H^0(\cpx{\ol{N}})\ar[r]& H^{-1}(\cpx{\ol{N}})[2]}
$$
because $\Hom_{\D{R}}(H^{-1}(\cpx{\ol{M}})[2],
H^0(\cpx{\ol{N}}))=0$. Then
$$
\eta_M
l(g)=\varphi_M\gamma_M\pi_Ml(g)=\varphi_M\gamma_M{H^0(g')}\pi_N=
\varphi_Mg'\gamma_N\pi_N=g\varphi_N\gamma_N\pi_N=g\eta_N.
$$
Hence, we have proved the existence of a homomorphism $l(g)$ from
$l(M)$ to $l(N)$ with the desired property.

Now we show the uniqueness of $l(g)$. Let $t_1, t_2:l(M)\to l(N)$ be
two homomorphisms such that $\eta_M t_{i}=g\eta_N$ for $i=1,2$. Set
$t:=t_1-t_2$. Then $\eta_Mt=\varphi_M(\gamma_M\pi_Mt)=0$. It follows
from the triangle $(\ast)$ that there exists a morphism
$u:\cpx{P}{^{(\alpha)}}\to l(N)$ in $\D{R}$ such that
$\gamma_M\pi_Mt=\psi_Mu$. But $\Hom_{\D{R}}(\cpx{P},l(N))=0$ since
$l(N)\in{\cpx{P}}^{\perp}$. This shows $u=0$ and $\gamma_M\pi_Mt=0$.
By the triangle $(\ast\ast)$, we know $\pi_Mt=0$ since
$\Hom_{\D{R}}(H^{-1}(\cpx{\ol{M}})[2],\, l(N))=0$. Note that $\pi_M$
is surjective. Hence $t=0$, that is, $t_1=t_2$. This finishes the
proof of the uniqueness.

Consequently, if $g$ is an isomorphism, then $l(g)$ is an
isomorphism. This shows also that, regardless of different choices
of $\omega_M$, the module $l(M)$ is unique up to isomorphism.

\textbf{Step $(3)$}. We define a functor $l:
R\Modcat\to{\cpx{P}}^{\perp}$ by sending $M$ in $R\Modcat$ to
$l(M)$, and homomorphism $g: M\to N$ to $l(g):l(M)\to l(N)$. By Step
$(2)$, $l$ is well-defined. Clearly, we have $l(Y)\simeq Y$ for any
$Y\in {\cpx{P}}^{\perp}$ by definition. Now, it follows from Step
$(2)$ that $\Hom_R(l(U), V)\simeq\Hom_R(U, j(V))$ for any $U\in
R$-Mod and $V\in {\cpx{P}}^{\perp}$. This isomorphism is natural in
$U$ and $V$. Thus $(l, j)$ is an adjoint pair of functors. $\square$

\medskip
As a consequence of  Proposition \ref{prop32}, we obtain  the
following promised result which provides an effective description of
$R_{\Sigma}$ as endomorphism rings.

\begin{Koro}
Let  $f: P^{-1}\to P^{0}$ be  a homomorphism between finitely
generated projective  $R$-modules and $\Sigma := \{f\}$. If
$\Hom_{\D{R}}(\cpx{P_f}, \cpx{P_f}[1])= 0$, then the inclusion
$j:\Sigma^{\perp} \to R\Modcat$ admits a left adjoint $l:
R\Modcat\to\Sigma^{\perp}$, which can be constructed explicitly. In
particular, $R_{\Sigma}$ is isomorphic to the endomorphism ring
$\End_{R}(l(R))$. \label{37}
\end{Koro}

{\bf Proof.} Since $\cpx{P_f}\in \Cb{\pmodcat R}$, we know that the
functor $\Hom_{\D{R}}(\cpx{P_f}, -)$ commutes with direct sums.
Consequently, if $\Hom_{\D{R}}(\cpx{P_f}, \cpx{P_f}[1])= 0$, then
$\Hom_{\D{R}}(\cpx{P_f}, \cpx{P_f}{^{(\delta)}}\,[1])=0$ for any
cardinal $\delta$. Thus the existence of the functor $l$ in
Corollary \ref{37} follows immediately from Proposition
\ref{prop32}.

In the following, we shall prove that $R_{\Sigma}$ is isomorphic to
the endomorphism ring $\End_{R}(l(R))$. Indeed, by Proposition
\ref{prop31}(1), we know that $\Sigma^{\perp}$ is closed under
extensions, kernels, cokernels, arbitrary direct sums and products.
Further, since the inclusion $j$ admits a left adjoint functor $l$,
the full subcategory $\Sigma^{\perp}$ of $R\Modcat$ satisfies all
assumptions of \cite[Proposition 3.8]{GL}. Define $S:=\End_R(l(R))$.
Then it follows  directly from \cite[Proposition 3.8]{GL}
 that $l(R)$ is a projective generator for $\Sigma^{\perp}$,
 and there exists a ring epimorphism $\rho$: $R\to S$ such that the
following diagram
$$
\xymatrix@C=1.2cm{\Sigma^{\perp}\ar@<0.4ex>[r]^-{j}\ar@<0.4ex>[d]^-{G}
&R\Modcat\ar@<0.6ex>[l]^-{l}\ar@<1.5ex>[ld]^-(.4){\rho^*}\\
S\Modcat\ar@<0.6ex>[u]^-{F} \ar@<-0.4ex>[ur]^-(.3){\rho_*}&}
$$
commutes, where $F:= {}_Rl(R)\otimes_S-$ and $G:=\Hom_R(l(R), -)$
are mutually inverse functors, and where $ \rho_*:= {}_RS\otimes_S-$
is the canonical embedding, and $\rho^*:= {}_SS\otimes_R-$ is a left
adjoint of $\rho_*$. Since the both ring epimorphisms $\lambda:R\to
R_{\Sigma}$ and $\rho$: $R\to S$ give rise to the same bireflective
subcategory $\Sigma^{\perp}$ of $R$-Mod, we conclude from Lemma
\ref{2.0} that $R_{\Sigma}$ is isomorphic to $S$. $\square$

\medskip
Finally, we remark that, in general, the two-term complex $\cpx{P}$
in Proposition \ref{prop32} cannot be replaced by a complex in
${\mathscr C}(\Pmodcat{R})$ with more than two terms. A
counterexample is the following:

Let $A$ be the algebra given by the following quiver with relations:
$$
\xymatrix{3\ar^-{\gamma}@/_1.0pc/[r]
&2\ar_-{\beta}@/_1.0pc/[l]\ar^-{\alpha}[r]&1,\qquad
\beta\gamma=\gamma\alpha=0.}
$$
We denote by $S_i$, $I_i$ and $P_i$ the simple, injective and
projective modules corresponding to the vertex $i$, respectively.
Let $\cpx{P}$ be the minimal projective resolution of $S_3$. Then
$\cpx{P}$ is a three-term complex. We can easily check that
${\cpx{P}}^{\perp}$ has only two indecomposable modules, they are
the indecomposable modules $I_1$ and $I_2$. Note that $A$ is
representation-finite and every indecomposable module is finitely
generated. Also, we have Ext$^i_A(S_3,S_3)=0$ for $i=1,2$. Actually,
this is true for all $i>0$. Since $\cpx{P}$ is compact, one has
$\Hom_{\D{A}}(\cpx{P},{\cpx{P}}^{(\alpha)}[1])=0$ for all cardinal
$\alpha$. If the inclusion functor $j$ from ${\cpx{P}}^{\perp}$ into
$A$-Mod would have a left adjoint, then $j$ would preserve injective
homomorphisms. One can verify that there is a non-zero homomorphism
from $I_1$ to $I_2$ which is a monomorphism in ${\cpx{P}}^{\perp}$,
but not a monomorphism in $A$-Mod. This is a contradiction and shows
that the inclusion functor from ${\cpx{P}}^{\perp}$ into $A$-Mod
cannot possess a left adjoint.

Note that the simple module  corresponding to the vertex $1$ is of
injective dimension $3$, and defines a $2$-APR-tilting module
$S_3\oplus P_2\oplus P_3$ (see \cite{hx2} for unexplained
definitions).

\section{Recollements of derived categories and infinitely generated tilting modules \label{sect4}}

In this section, we shall use our results in Section \ref{sect3} to
show the first statement of the main result, Theorem \ref{th01}.
More precisely, we first recall the definition of infinitely
generated tilting modules, and then discuss some of their
homological properties. Especially, we shall establish a crucial
result, Proposition \ref{prop1}, which will play a role in our proof
of the main result.

Let $A$ be a ring with identity.

\begin{Def}\rm \cite{ct} \label{def1}
An  $A$-module $T$ is called a tilting module (of projective
dimension at most one) if the following conditions are satisfied:

$(T1)$ the projective dimension of $T$ is at most $1$, that is,
there exists a projective resolution of $T$: $0\rightarrow
P_1\rightarrow P_0\rightarrow T\rightarrow 0,$ where $P_i$ is
projective for $i=0,1$.

$(T2)$ $\Ext^i_A(T,T^{({\alpha})})=0$ for each $i\geq 1$ and every
cardinal $\alpha$; where $T^{(\alpha)}$ stands for the direct sum of
$\alpha$ copies of $T$; and

$(T3)$ there exists  an exact sequence
$$0\lra {}_AA\lra
T_0\longrightarrow T_1\longrightarrow 0$$ of $A$-modules such that
$T_i\in\Add(T)$ for $ i=0,1$.

\medskip
If $P_1$ and $P_0$ in the condition $(T_1)$ are finitely generated,
then the tilting module $T$ is called a classical tilting module
(see \cite{BB} and \cite{HR}).

Two tilting $A$-modules $T$ and  $T'$ are said to be equivalent if
$\Add(T)=\Add(T')$, or equivalently, $\Gen(T)=\Gen(T')$, where
$\Gen(T)$ denotes the full subcategory of $A$-Mod generated by $T$.
Recall that an $A$-module $M$ is generated by $T$ if there is an
index set $I$ and a surjective homomorphism $f: T^{(I)}\ra M$.

An $A$-module $T$ is said to be good if it satisfies (1), (2) and

$(T3)'$ there is an exact sequence
$$0\lra {}_AA\lra
T_0\lra T_1\lra 0$$ in $A\Modcat$ such that $T_i\in\add(T)$ for $ i
= 0,1$.
\end{Def}

Note that each classical tilting module is good. Moreover, for any
given tilting module $_AT$ with (T1) and (T2), the module
$T':=T_0\oplus T_1$ is a good tilting module which is equivalent to
the given one.

\medskip
From now on, we assume in this section that $T$ is a {\bf{good}}
tilting $A$-module, namely, it satisfies $(T_1)$, $(T_2)$ and
$$(T_3)':\quad 0\lra A\longrightarrow T_0\lraf{\varphi} T_1\longrightarrow
0$$  with  $\varphi$ a homomorphism of $A$-modules between $T_i\in
\add(T)$. Let $B:=\End_A(T).$ We define
$$T^{\bot}:=\{X\in A\Modcat \mid \Ext^i _A(T,X)=0 \mbox{\;for\, all \,} i\geq
1\},\;\mathscr{E}:=\{Y\in B\Modcat \mid \Tor_i
^B(T,Y)=0\mbox{\;for\, all \,} i\geq0\}; $$
$$G:= {}_AT\otimesL_B-:
\D{B}\lra\D{A},\quad H:=\rHom_A(T,-): \D{A}\lra\D{B};$$
$$\mathcal{Y}:=\Ker(G),\qquad \mathcal{Z}:=\Img(H),$$
$$\cpx{Q}:= \quad \cdots \lra 0\lra\Hom_A(T, T_0)
\xrightarrow{\varphi{^*}}\Hom_A(T, T_1)\lra 0\lra
\cdots\,\in\Cb{\pmodcat B},$$ where
$\varphi{^*}:=\Hom_A(T,\varphi)$, and where the finitely generated
projective $B$-modules $\Hom_A(T,T_0)$ and $\Hom_A(T,T_1)$, as terms
of the complex $\cpx{Q}$, are concentrated on the degrees $0$ and
$1$, respectively. Clearly, $H(A)=\cpx{Q}$ in $\D{A}$.

In the next lemma we mention a few basic properties of tilting
modules. For proofs we refer to \cite[Proposition 1.4, Lemma
1.5]{Bz2} and \cite{Bz}.

\begin{Lem} Let $T$ be a tilting $A$-module. Then:

$(1)$ $T_B$ has a projective resolution
$$0\lra Q_1\lra Q_0\lra
T_B\lra 0$$ such that  $Q_i\in\add(B_B)$ for each $ 0\leq i\leq 1$.

$(2)$  $\End{_{B^{\opp}}}(T)\simeq A^{\opp}$ and $\Ext^i_B(T,T)=0$
for all $ i\geq 1$.

$(3)$  For each  $ Y\in\Add(_BB)$, we have $\Ext^i _A(T,
{}_AT\otimes_B Y)=0$ for each  $  i\geq 1$.

$(4)$ For each $X\in T^{\bot}$, we have $\Tor_i ^B(_AT_B,\,
\Hom_A(T,X))\simeq
\begin{cases}
X, & i=0,\\
0, & i > 0.
\end{cases}$

$(5)$ $T^{\bot}$ is closed under direct sums. \label{lem1}
\end{Lem}

The following result is shown in \cite[Theorem 5.1]{Bz}, which says
that the unbounded derived category of $B\Modcat$ is bigger than
that of $A\Modcat$ in general.

\begin{Lem}
The functor  $H$ is fully faithful, and the functor $G$ induces a
triangle equivalence between $\D{B}/\Ker(G)$ and $\D{A}$. Here we
denote by $\D{B}/\Ker(G)$ the Verdier quotient of $\D{B}$ by the
category $\Ker(G)$. \label{lem2}
\end{Lem}

The following lemma supplies a method to obtain modules in $\mathscr
E$, and is also useful for our later calculations.

\begin{Lem}
Suppose that  $I$ is  a cardinal and  $X_i\in T^{\bot}$ for each
$i\in I$. Consider the canonical exact sequence
$$0\lra\bigoplus_{i\in I}\Hom_A(T,
X_i)\lraf{\delta_I}\Hom_A(T, \bigoplus_{j\in I}X_j)\lra\Coker
(\delta_I)\lra 0$$ in $B\Modcat$, where $\delta_I$ is defined by
$(f_i)_{i\in I}\mapsto \sum_{i\in I}f_i\lambda_i$ with $f_i\in
\Hom_A(T, X_i)$ and $\lambda_i:X_i\to\bigoplus{_{j\in I}}X_j$ the
canonical inclusion for each $i\in I$. Then $\Coker(\delta_I)\in
\mathscr{E}$. Particularly, for each projective $B$-module $P$, the
unit adjunction morphism $\eta_P': P \to\Hom_A(T, T\otimes_BP)$ is
injective with $\Coker(\eta_P')\in \mathscr{E}$. \label{lem3}
\end{Lem}

{\it Proof.} Note that $\delta_I$ is well-defined. By the definition
of $\delta_I$, we can see easily that $\delta_I$ is injective. So,
there is a canonical exact sequence
$$(\ast)\quad
0\lra\bigoplus_{i\in I}\Hom_A(T, X_i)\lraf{\delta_I}\Hom_A(T,
\bigoplus_{j\in I}X_j)\lra\Coker (\delta_I)\lra 0 .
$$
Since $T^{\bot}$ is closed under direct sums by Lemma \ref{lem1}(5),
we have  $\bigoplus{_{j\in I}}X_j\in T^{\bot}$. It then follows from
Lemma \ref{lem1}(4) that
$$\Tor_m ^B(T,\Hom_A(T,\displaystyle\bigoplus_{j\in
I}X_j))\simeq
\begin{cases}
\displaystyle\bigoplus_{j\in
I}X_j, & m=0,\\
0, & m>0.
\end{cases}$$
Similarly, for any $i\in I$, we have
$$\Tor_n ^B(T,\Hom_A(T,X_i))\simeq
\begin{cases}
X_i, & n=0,\\
0, & n>0.
\end{cases}$$
Since the right module $T_B$ has projective dimension at most $1$,
we see that $\Tor_t ^B(T,\Coker(\delta_I))=0$ for any  $t>1$. By
applying the functor $_AT\otimes_B-:B\Modcat\to A\Modcat$ to the
sequence $(\ast)$, we can easily form the following exact
commutative diagram:
$$
\xymatrix@C=0.3cm {0\ar[r]&\Tor_1 ^B(T,\Coker
(\delta_I))\ar[r]&{T\otimes_B}{(\bigoplus_{i\in
I}\Hom_A(T,X_i}))\ar^-{\simeq}[d]\ar[r]&T\otimes_B\Hom_A(T,\bigoplus_{j\in
I}X_j)\ar[r]\ar^-{\simeq}[d]&T\otimes_B \Coker
(\delta_I)\ar[r]& 0\\
&&\bigoplus{_{i\in I}}X_i\ar@{=}[r]&\bigoplus{_{j\in I}}X_j.&}
$$
This implies $T\otimes_B\Coker (\delta_I)= 0 =\Tor_1
^B(T,\Coker(\delta_I))$. Hence $\Coker(\delta_I)\in\mathscr{E}$.

To prove the last statement of Lemma \ref{lem3}, we note that the
unit adjunction
$$
\eta':1_{B\Modcat}\lra \Hom_A(T,T\otimes_B-)
$$ is a natural transformation of functors from $B\Modcat$ to itself, and that
$\mathscr{E}$ is closed under direct summands. Thus, it is
sufficient to show that the statement holds for free $B$-modules.
Let $\alpha$ be any cardinal. Then we may form the following exact
commutative diagram:
$$
\xymatrix{ & B^{(\alpha)}\ar[r]^-{\eta_{B^{(\alpha)}}'}\ar@{=}[d]
 &\Hom_A(T,T\otimes_BB^{(\alpha)})\ar[d]^-{\simeq}\\
0\ar[r]&{\Hom_A(T,T)}^{(\alpha)}\ar[r]^-{\delta_{\alpha}}
&\Hom_A(T,T^{(\alpha)})\ar[r]&\Coker(\delta_{\alpha})\ar[r]&0. }
$$
Since $\delta_{\alpha}$ is injective, we conclude that
$\eta_{B^{(\alpha)}}'$ also is injective, and therefore
$\Coker(\eta_{B^{(\alpha)}}')\simeq
\Coker(\delta_{\alpha})\in\mathscr{E}$. This finishes the whole
proof. $\square$

\medskip
In the next lemma we give a description of the category $\mathscr
E$.

\begin{Lem} The following statements hold.

$(1)$
 $\mathscr{E}= \{X\in B\Modcat \mid \Hom{_{\D{B}}}(\cpx{Q},
X[i])= 0 \mbox{\;for\, all \,}i\in\mathbb{Z}\}.$ In particular,
$\mathscr{E}$ is closed under direct sums and products.

$(2)$ $\mathscr{E}$ of $B\Modcat$ is closed under isomorphic images,
extensions, kernels and cokernels. In particular, $\mathscr{E}$ is
an abelian subcategory of $B\Modcat$.

\label{lem4}
\end{Lem}

{\bf Proof.} $(1)$ Let $X$ be a $B$-module and $i$ be an integer.
Then
$$\Hom_{\D{B}}(Q^\bullet, X[i])\simeq\Hom_{\K{B}}(Q^\bullet, X[i])
\simeq H^i(\Hom_B(Q^\bullet, X))\simeq H^i(\Hom_B(Q^{\bullet},B)
\otimes_BX),$$ where the last isomorphism follows from the fact that
the restriction of the natural transformation $\Hom_B(-,
B)\otimes_BX\ra \Hom_B(-,X)$ to $\C{\pmodcat {B}}$ is a natural
isomorphism. By the definition of $Q^\bullet$, we know that
$\Hom_B(Q^{\bullet},B)$ is the complex:
$$
\,\cdots \lra 0\lra\Hom_A(T_1,
T)\xrightarrow{\varphi{_*}}\Hom_A(T_0, T)\lra 0\lra \cdots
$$
in $\Cb{\pmodcat {B^{\opp}}}$, where $\varphi{_*}:=\Hom_A(\varphi,
T)$, and where the finitely generated projective $B^{\opp}$-modules
$\Hom_A(T_1,T)$ and $\Hom_A(T_0,T)$ are of degrees $-1$ and $0$,
respectively. Note that the conditions $(T_2)$ and $(T_3)$ in
Definition \ref{def1} imply that the sequence
$$
0\lra\Hom_A(T_1, T)\xrightarrow{\varphi{_*}}\Hom_A(T_0, T)\lra T\lra
0
$$ is exact.
In other words, the complex $\Hom_B(Q^{\bullet}, {}_BB)$ is
quasi-isomorphic to $T_B$ (here we use the fact that the functor
$\Hom_A(-,T): \add(_AT)\ra \add(B_B)$ is an equivalence of
categories). It follows from the definition of $\Tor_i^B$ that
$$
H^i\big(\Hom_B(Q^{\bullet},B)\otimes_BX\big)\simeq \left\{
\begin{array}{ll}
0 & \mbox{ if } i>0,\\
\Tor^B_{-i}(T,X) & \mbox{ if } i\leq 0.
\end{array}
\right.
$$
This means that $\Hom_{\D{B}}(Q^\bullet, X[i])=0$ if and only if
$\Tor^B_{-i}(T,X)=0.$ Hence
$$\mathscr{E}= \{X\in B\Modcat \mid
\Hom{_{\D{B}}}(\cpx{Q}, X[i])= 0 \mbox{\;for\, all
\,}i\in\mathbb{Z}\}.$$ Consequently, $\mathscr{E}$ is closed under
 direct products. Further, since $\cpx{Q}$ is a bounded complex of finitely
generated projective  $B$-modules, we know that $\mathscr{E}$ is
closed under  direct sums, too.

$(2)$ This statement follows directly from Proposition
\ref{prop31}(1). $\square$
\medskip

The following proposition is crucial to the proof of Theorem
\ref{th01}(1).

\begin{Prop}
The triple $\big({\rm{Tria}}(\cpx{Q}),\Ker(G),\Img(H)\big)$ is a
{\rm TTF} tripe in $\D{B}$. Moreover,
$$\Ker(G)=\{\cpx{\overline{Y}}\in\D{B}\mid\cpx{\overline{Y}}\simeq \cpx{Y} {\rm \; in \;}\D{B}
{\rm \; with \; } Y^i\in\mathscr{E} {\rm \; for \ all \;}
i\in\mathbb{Z}\};
$$
$$\Img(H)=\{\cpx{\overline{Z}}\in\D{B}\mid\cpx{\overline{Z}}\simeq
\cpx{Z} \mbox{\; in \;}\D{B} \mbox{\; with\;}
Z^i\in\Hom_A(T,\Add(T))\mbox{\;for\, all \,} i\in\mathbb{Z}\},$$
where $\Hom_A(T,\Add(T))$ stands for the full subcategory of
$B\Modcat$ consisting of all the modules $\Hom_A(T,T')$ with $T'$ in
$\Add(T)$. \label{prop1}
\end{Prop}
{\it Proof.} Recall that we have denoted $\Ker(G)$ by $\mathcal Y$,
and Im$(H)$ by ${\mathcal Z}$. The whole proof of this proposition
will be divided into three steps.

\textbf{Step $(1)$}. We prove that  the  pair
$(\mathcal{Y},\mathcal{Z})$ is a torsion pair in $\D{B}$. In fact,
for any $\cpx{Y}\in\mathcal{Y}$ and $\cpx{W}\in\D{A}$, we have
$\Hom_{\D{B}}(\cpx{Y},H(\cpx{W}))\simeq\Hom_{\D{A}}(G(\cpx{Y}),\cpx{W})=
\Hom_{\D{A}}(0,\cpx{W})=0$ because the pair $(G,H)$  is an adjoint
pair of triangle functors by Lemma \ref{lem2}.  This shows
$\Hom{_{\D{B}}}(\mathcal{Y},\mathcal{Z})=0$. Let $\eta:
Id_{\D{B}}\to HG$ be the unit adjunction, and let $\varepsilon:GH\to
Id_{\D{A}}$ be the counit adjunction. By Lemma \ref{lem2}, we know
that $\varepsilon$ is invertible. For any $\cpx{M}$ in $\D{B}$, the
canonical morphism $\eta_{\cpx{M}}:\cpx{M}\to HG(\cpx{M})$ can be
extended to a triangle in  $\D{B}$:
$$\cpx{M}\epa{\eta_{\cpx{M}}} HG(\cpx{M})\lra \cpx{N}\lra \cpx{M}[1].$$
By applying the functor $G$ to the above triangle, we obtain a
triangle in $\D{A}$:
$$G(\cpx{M})\epa{G(\eta_{\cpx{M}})} GHG(\cpx{M})\lra
G(\cpx{N})\lra G(\cpx{M})[1].
$$  Since $\varepsilon$ is invertible, we see that $G(\eta_{\cpx{M}})$
is an isomorphism. This shows $G(\cpx{N})=0$, that is,
$\cpx{N}\in\mathcal{Y}$. Since  $\mathcal{Y}$ is a triangulated
subcategory of $\D{B}$, we have $\cpx{N}[-1]\in\mathcal{Y}$. Thus
the following triangle
$$
(\ast)\quad \cpx{N}[-1]\lra \cpx{M}\epa{\eta_{\cpx{M}}}
HG(\cpx{M})\lra \cpx{N}
$$
in $\D{B}$ with $HG(\cpx{M})\in\mathcal{Z}$ shows that the third
condition of Definition \ref{def02} is satisfied. Hence the pair
$(\mathcal{Y}, \mathcal{Z})$ is a torsion pair in $\D{B}$ by
Definition \ref{def02}. Since $\mathcal Y$ is a triangulated
category, the torsion pair $(\mathcal{Y, Z})$ is hereditary.

\textbf{Step $(2)$}. We calculate the categories $\mathcal{Y}$ and
$\mathcal{Z}$. Before starting our calculations, we mention the
following result in \cite[Theorem 10.5.9, Corollary 10.5.11]{We}:

For every complex $\cpx{X}$ in $\D{B}$, there exists a
quasi-isomorphism $\cpx{\ol{X}}\to\cpx{X}$ with $\cpx{\ol{X}}$ a
complex of $(_AT\otimes_B-)$-acyclic $B$-modules such that
$G(\cpx{X}) \simeq T\otimes_B\cpx{\ol{X}}$. Here, a $B$-module $N$
is said to be $(_AT\otimes_B-)$-acyclic if $\Tor^B_{i}(T,N)=0$ for
any $i>0$. Thus the action of the left derived functor $G$ on any
complex $\cpx{U}$  of $(_AT\otimes_B-)$-acyclic $B$-modules is the
same as that of the functor $_AT\otimes_B-$ which acts in component
wise on each term of $\cpx{U}$.

A similar statement holds for the right derived functor $H$.

\medskip
Now let us first interpret the triangle $(\ast)$ in terms of objects
in $\C{\Pmodcat B}$. For the complex $\cpx{M}$, we choose
$\cpx{P}\in\C{\Pmodcat B}$ such that $\cpx{P}$ is quasi-isomorphic
to the complex $\cpx{M}$. Then $G(\cpx{M}) \simeq
T\otimes_B\cpx{P}$. By Lemma \ref{lem1}(3), we have
$HG(\cpx{M})=\Hom_A(T, T\otimes_B\cpx{P})$ because the $A$-module
$T\otimes_BP$ is $\Hom_A(T, -)$-acyclic for any projective
$B$-module $P$. Note that the homomorphism $\eta_{\cpx{P}}$
coincides with $(\eta_{P^n}')_{n\in\mathbb{Z}}$, where $P^n$ is the
$n$-th term of the complex $\cpx{P}$ and
$\eta_{P^n}':P^n\to\Hom_A(T,T\otimes_BP^n)$ is the unit adjunction
 morphism for each $n\in\mathbb{Z}$. By Lemma \ref{lem3}, there is a short exact sequence of complexes
$$0\lra \cpx{P}\epa{\eta_{\cpx{P}}}\Hom_A(T, T\otimes_B \cpx{P})\lra \Coker(\eta_{\cpx{P}})\lra 0$$
such that
$(\Coker(\eta_{\cpx{P}}))^i=\Coker(\eta'_{P^i})\in\mathscr{E}$ for
each $i\in\mathbb{Z}$. Thus, we can  form the following commutative
diagram of triangles in $\D{B}$:
$$
\xymatrix{
\Coker(\eta_{\cpx{P}})[-1]\ar[r]\ar@{-->}^{\simeq}[d]&\cpx{P}\ar[r]^-{\eta_{\cpx{P}}}\ar[d]^-{\simeq}
&\Hom_A(T, T\otimes_B \cpx{P})\ar@{=}[d]\ar[r]& \Coker(\eta_{\cpx{P}})\ar@{-->}^{\simeq}[d]\\
 \cpx{N}[-1]\ar[r]&\cpx{M}\ar[r]^-{\eta_{\cpx{M}}}
 &\rHom_A(T, T\otimesL_B \cpx{M})\ar[r]&\;\; \; \cpx{N}.
 }$$
On the one hand, if $\cpx{M}\in\mathcal{Y}$, then
$T\otimesL_B\cpx{M}=0$ by definition, and so
$\cpx{M}\simeq\Coker(\eta_{\cpx{P}})[-1]$ in $\D{B}$. On the other
hand, if $\cpx{M}\simeq\cpx{Y}$ in $\D{B}$ for some complex
$\cpx{Y}$ with $Y^i\in\mathscr{E}$ for each $i\in\mathbb{Z}$, then
$T\otimesL_B\cpx{M}\simeq T\otimesL_B\cpx{Y}=T\otimes_B\cpx{Y}=0$ by
the above mentioned fact. This means $\cpx{M}\in\mathcal{Y}$. Hence
the first equality in Proposition \ref{prop1} holds.

To prove the second equality, we observe that, by Lemma
\ref{lem1}(4), $\Hom_A(T,T\otimes_B\Hom_A(T, T'))\simeq \Hom_A(T,
T')$ for any $T'\in\Add(T)$. Let $\cpx{Z}$ be a complex in $\D{B}$
such that $Z^i\in\Hom_A(T,\Add(T))$. Then
$HG(\cpx{Z})\simeq\Hom_A(T,T\otimes_B \cpx{Z})\simeq\cpx{Z}$  in
$\D{B}$ because every $B$-module in $\Hom_A(T, \Add(T))$ is
$(T\otimes_B-)$-acyclic by Lemma \ref{lem1}(3) and Lemma
\ref{lem1}(4). This implies $\cpx{Z}\in\mathcal{Z}$. Conversely, for
any $\cpx{W}\in\D{A}$, we can choose a complex
$\cpx{L}\in\C{\Pmodcat B}$ such that $\cpx{L}$ is quasi-isomorphic
to  $H(\cpx{W})$. By Lemma \ref{lem2}, we conclude that
$H(\cpx{W})\simeq HG(H(\cpx{W}))\simeq HG(\cpx{L})$ in $\D{B}$.
Since $HG(\cpx{L})= H(T\otimesL_B\cpx{L})=
H(T\otimes_B\cpx{L})\simeq \Hom_A(T,T\otimes_B\cpx{L})$, where the
last isomrphism follows from Lemma \ref{lem1}(3) and the above
mentioned fact about the functor $H$. Clearly, the complex
$\Hom_A(T,T\otimes_B\cpx{L})$ has each term in $\Hom_A(T,\Add(T))$.
Thus the second equality in Proposition \ref{prop1} holds.

\textbf{Step $(3)$}. We claim that there is a full subcategory
category $\mathcal X$ of $\D{B}$ such that
$(\mathcal{X},\mathcal{Y},\mathcal{Z})$ is a TTF triple in $\D{B}$.
Furthermore, we have $\mathcal{X} = {\rm{Tria}}(\cpx{Q})$.

Indeed, since $\mathscr{E}$ is closed under direct sums and products
by Lemma \ref{lem4}, we conclude that $\mathcal{Y}$ is closed under
all small coproducts and products. Then the existence of the TTF
triple $(\mathcal{X},\mathcal{Y},\mathcal{Z})$ in $\D{B}$ follows
straightforward from Lemma \ref{lem02}. Moreover,
$\mathcal{X}=\Ker(\Hom_{\D{B}}(-,\mathcal{Y}))$ and
$\mathcal{Y}=\Ker(\Hom_{\D{B}}(\mathcal{X},-))$. Now we shall prove
$\mathcal{X}={\rm{Tria}}(\cpx{Q})$. First, we show
$\cpx{Q}\in\mathcal{X}$. This is equivalent to verifying
$\Hom_{\D{B}}(\cpx{Q}, \mathcal{Y})$ = $0$. Let
$\mathcal{Y'}:=\Ker(\Hom_{\D{B}}({\rm{Tria}}(\cpx{Q}), -)).$ By
Lemma \ref{lem03}, we  see that $({\rm{Tria}}(\cpx{Q}),
\,\mathcal{Y'})$ is a torsion pair in $\D{B}$  with
$$
\mathcal{Y'}=\{\cpx{Y}\in\D{B}~|~\Hom_{\D{B}}(\cpx{Q}, \cpx{Y}[i])=0
\mbox{ for\, all \,} i\in\mathbb{Z}\}.
$$
Recall that  $\varphi^*:=\Hom_A(T,
\varphi):\Hom_A(T,T_0)\to\Hom_A(T, T_1)$ is a homomorphism between
finitely generated projective $B$-modules. We define
$\Sigma:=\{\varphi^*\}$. Then $\cpx{\Sigma}=\{\cpx{Q}[1]\}$ (see
notations in Section \ref{sect3}). By Lemma \ref{lem4} (1), we have
$\Sigma^{\perp}=\mathscr{E}$. Thus it follows from Proposition
\ref{prop31} that
 $$
 \mathcal{Y'}=\mathscr{D}(B)_{\mathscr{E}}:=\{\cpx{Y}\in\D{B}~|~H^i(\cpx{Y})\in\mathscr{E}
\mbox{ for\, all \,}i\in\mathbb{Z}\}.
$$  According to Lemma
\ref{lem4} (2), $\mathscr{E}$ is an abelian subcategory of
$B\Modcat$. This forces $\mathcal{Y}\subseteq\mathcal{Y'}$. In
particular,  we have $\Hom_{\D{B}}(\cpx{Q}, \mathcal{Y})=0$, which
yields $\cpx{Q}\in\mathcal{X}$. Therefore,
${\rm{Tria}}(\cpx{Q})\subseteq\mathcal{X}$ since $\mathcal X$ is a
full triangulated subcategory of $\D{B}$.

Let ${\bf{i}}:\mathcal{X}\to\D{B}$ and
${\bf{k}}:\mathcal{Z}\to\D{B}$ be the canonical inclusions.  Then
the functor  $\bf{i}$ has a right adjoint functor
${\bf{R}}:\D{B}\to\mathcal{X}$.
 Since $(\mathcal{X},\mathcal{Y}, \mathcal{Z})$ is a TTF
triple in $\D{B}$, the functor ${\bf{Rk}}:\mathcal{Z}\to\mathcal{X}$
is an equivalence (see the statements after Definition \ref{def03}
in Subsection \ref{sect2.2}). So the composition functor
${\bf{Rk}}H:\D{A}\to\mathcal{X}$ is an equivalence because
$H:\D{A}\to\mathcal{Z}$ is an equivalence. Since a functor
possessing a right adjoint functor preserves coproducts, we know
that the functor ${\bf{Rk}}H$ commutes with coproducts. Note that
$\D{A}$ admits all small coproducts and that the notion of
coproducts depends on the category where coproducts are taken. In
general, $\mathcal{Z}$ is not closed under coproducts.

Since $(\bf{i, R})$ is an adjoint pair, we know that the functor
$\bf{i}$ preserves coproducts. This means that a coproduct in
$\mathcal X$ is the same as that in $\D{B}$.

Since $H(A)\simeq\cpx{Q}\in\mathcal{X}$, we have ${\bf
Rk}H(A)\simeq{\bf R}(\cpx{Q})=\cpx{Q}$. Note that
$\D{A}={\rm{Tria}}(A)$ and that the triangle functor ${\bf
Rk}H:\D{A}\to\mathcal{X}$ is an equivalence under which ${\rm
Tria}(A)$ has the image ${\rm Tria}(\cpx{Q})$ since the functor
${\bf Rk}H$ commutes with  coproducts. It follows that
$\mathcal{X}={\rm{Tria}}(\cpx{Q})$ and $\mathcal{Y}=\mathcal{Y'}$.
Hence $\big({\rm{Tria}}(\cpx{Q}),\Ker(G),\Img(H)\big)$ is a \rm{TTF}
tripe in $\D{B}$. $\square$

\medskip
With the above preparations, now we prove Theorem \ref{th01}\,$(1)$.

\smallskip
{\bf Proof of Theorem \ref{th01}\,(1)}. By Proposition \ref{prop1},
we know that the triple $(\,{\rm{Tria}}(\cpx{Q}),\Ker(G),\Img(H)\,)$
is a \rm{TTF} tripe in $\D{B}$. Moreover, $\D{A}$ and
${\rm{Tria}}(\cpx{Q})$ are equivalent as triangulated categories.
According to the correspondence between recollements and TTF triples
in Lemma \ref{lem01}(2), we can form the following recollement
$$\xymatrix@C=1.2cm{\Ker(G)\ar^-{\bf{j}}[r]&\D{B}\ar[r]
\ar@/^1.2pc/[l]\ar_-{\bf{L}}@/_1.2pc/[l]
&\D{A}\ar@/^1.2pc/[l]\ar@/_1.2pc/[l]},$$
\medskip

\noindent where ${\bf{j}}$ is the canonical embedding and ${\bf{L}}$
is the left adjoint of ${\bf{j}}$. Recall that $\varphi^*:=\Hom_A(T,
\varphi)$ is the homomorphism between the finitely generated
projective $B$-modules $\Hom_A(T,T_0)$ and $\Hom_A(T, T_1)$. As in
Section \ref{sect3}, we define $\Sigma:=\{\varphi^*\}$. By Lemma
\ref{lem4}(1), we have $\Sigma^{\perp}=\mathscr{E}$. By Step (3) in
the proof of Proposition \ref{prop1}, we have
$\Ker(G)=\mathscr{D}(B){_{\Sigma^{\perp}}}$.  Let $\lambda:B\to
B_{\Sigma}$ be the universal localization of $B$ at $\Sigma$. Since
$\bf{L}$ is a functor from $\D{B}$ to $\Ker(G)$, we have
${\bf{L}}(B)\in\Ker(G)$, and therefore the condition (5) of
Proposition \ref{th2} is satisfied by ${\bf L}(B)$, according to
Proposition \ref{prop1}. Thus, by Proposition \ref{th2}, we know
that $\lambda_{*}: \D{B_{\Sigma}}\xrightarrow {\sim}
\mathscr{D}(B)_{\Sigma^{\perp}}$ is an equivalence of triangulated
categories, and  that the homomorphism $\lambda$ is a homological
ring epimorphism. Set $C:=B_{\Sigma}$. Then $\Ker(G)$ and $\D{C}$
are equivalent  as triangulated categories. Consequently, we can get
the following recollement from the above one:

$$\xymatrix@C=1.2cm{\D{C}\ar[r]&\D{B}\ar[r]
\ar@/^1.2pc/[l]\ar@/_1.2pc/[l]
&\D{A}\ar@/^1.2pc/[l]\ar@/_1.2pc/[l]}.$$\vspace{0.1cm}

\medskip
In the following, we shall explicitly describe the six triangle
functors arising in the above recollement.

Here, we follow the notations used in Definition \ref{def01}, and
take $\mathcal{D}=\D{B}, \mathcal{D'}=\D{A}$ and
$\mathcal{D''}=\D{C}$. Then it is not hard to see that
$i^*=C\otimesL_B-,\, i_*=\lambda_{*}$ and $i^!=\rHom_B(C,-)$. As for
the other three functors, we claim that $j_!={\bf{iRk}}H,\, j^!=G$
and $j_*=H$ up to natural ismorphisms.  Let
${\bf{U}}:\D{B}\to\mathcal{Z}$ be a left adjoint of the inclusion
${\bf{k}}:\mathcal{Z}\to\D{B}$. By Lemma \ref{lem01} and the proof
of Proposition \ref{prop1}, we get the following diagram
$$\xymatrix{\D{B}\ar^-{\bf{R}}[r]
&\mathcal{X}\ar^-{\vspace{0.2cm}\bf{kUi}}@/^1.0pc/[l]\ar_-{\bf{i}}@/_1.0pc/[l]&\D{A}\ar_-{{\bf{Rk}}H}[l]}$$
with the properties:

(i) $({\bf i, R})$ and $({\bf R, kUi})$ are adjoint pairs,

(ii) ${\bf Rk}H$ is an equivalence of triangulated categories.

\noindent This implies that $j_!={\bf{iRk}}H$ and
$j_*={\bf{(kUi)}}({\bf Rk}H)$. Note that
  the composition functor
${\bf{UiRk}}:\mathcal{Z}\to\mathcal{Z}$ of the functor ${\bf{Ui}}$
and ${\bf{Rk}}$ is natural isomorphic to the identity functor
$1_{\mathcal{Z}}$ by  the property (3) of a TTF triple (see
Subsection \ref{sect2.2}). Consequently, we can choose $j_*=H$.
Since $(G,H)$ is an adjoint pair of functors, we can choose $j^*=
G$. Thus the  proof of the first part of  Theorem \ref{th01} is
completed. $\square$

\medskip
{\it Remarks}. $(1)$ In the proof of Theorem \ref{th01}(1), we have
$\Sigma=\{\varphi^*\},\,\Sigma^{\perp}=\mathscr{E}$ and
$\Hom_{\D{B}}(\cpx{Q}, \cpx{Q}[1])= 0$. This means that the
homomorphism $\varphi^*$ satisfies the assumptions in Corollary
\ref{37}. Therefore we can explicitly construct a left adjoint
functor $l:B\Modcat\to\mathscr{E}$ of the inclusion $j:
\mathscr{E}\to B\Modcat$. In particular, we know that $C$ is
isomorphic to the endomorphism ring $\End_B(l(B))$ of $l(B)$.

$(2)$ The ring $C$ equals zero if and only if $T$ is a classical
tilting module. In fact, $C=0$ if and only if $\Ker(G)=0$ if and
only if $G$ is an equivalence if and only if $T$ is classical.

$(3)$ From the proof of Theorem \ref{th01}(1), we know that a good
tilting module $T$ has  the property:  the functor $G$ admits a
fully faithful left adjoint $j_!$. In the next section, we shall
show that this property guarantees that the tilting module $T$ is
good.

\section{Existence of recollements implies goodness of tilting modules\label{sect5}}

In this section, we shall  prove the second part of  Theorem
\ref{th01}, which is a converse of the first part in some sense. Our
proof depends on the property that the total left derived functor
$G$ admits a fully faithful left adjoint $j_!$.

\medskip
{\bf Proof of Theorem \ref{th01}\,(2)}.

Let $T$ be a tilting $A$-module and $B$ the endomorphism algebra of
$T$. Recall that $G$ and $H$ stand for the triangle functors
$T\otimesL_B-: \D{B}\ra \D{A}$ and $\rHom_A(T,-): \D{A}\ra\D{B}$,
respectively. Suppose that $G$ admits a fully faithful left adjoint
$j_!: \D{A}\to\D{B}$. We want to show that $T$ is a good tilting
module.

To prove that $T$ is good, it suffices to find a short exact
sequence of $A$-modules,
$$0\lra A\lra
T_0\lra T_1\lra 0,$$ such that $T_i\in\add(T)$ for $ i=0,1$.

First, we observe some consequences of the assumption that $j_!$ is
fully faithful. Set $\cpx{W}:=j_!(A)$. Since the total left derived
functor $G$ commutes with coproducts, we can easily show that the
functor $j_!$ preserves compact objects. In particular, the complex
$\cpx{W}$ is compact in $\D{B}$, which implies $\cpx{W}\simeq
\cpx{Q}$ in $\D{B}$ for some $\cpx{Q}\in\Cb{\pmodcat B}$. Since the
Hom-functor $\Hom_A(T,-)$ induces an equivalence between $\add(T)$
and $B\pmodcat$, we can assume that $\cpx{Q}=\Hom_A(T,\cpx{X})$,
where $\cpx{X}\in\Cb{\add(T)}$ is of the following form
$$
0\lra X^s\lra\cdots\lra X^i\lraf{d^i} X^{i+1}\lra \cdots\lra X^t\lra
0$$  for $s\le 0\le t$. Since the functor $j_!$ is fully faithful,
we conclude from  \cite[Chapter IV, Section 3, Theorem 1, p.90]{m}
that the unit adjunction morphism $\widetilde{\eta}:Id_{\D{A}}\to
Gj_!$ is invertible. Thus $A\simeq G(\cpx{W})\simeq G(\cpx{Q})$ in
$\D{A}$. Note that $T\otimes_B\Hom_A(T,\cpx{X})\simeq\cpx{X}$ in
$\Cb{\pmodcat A}$ since $X^i\in\add(T)$ for each $s\leq i \leq t$.
Consequently, we have $A\simeq\cpx{X}$ in $\D{A}$. It follows that
$H^0(\cpx{X})\simeq A$ and $H^i(\cpx{X})=0$ for any $ i\neq 0$.

Second, if $t=0$, then the homomorphism $X^0\ra H^0(\cpx{X})$
splits, this implies $A\in\add(T)$. Hence $T$ is a good tilting
module. Now we assume $t\neq0$. Then we can decompose $\cpx{X}$ into
two long exact sequence of $A$-modules:
$$
0\lra X^s\lraf{d^s}\cdots\lra X^{-1}\lraf{d^{-1}}
X^0\lraf{\pi}M\lra0,
$$
$$ 0\lra A\lra M\lraf{\mu}X^1\lraf{d^1}\cdots\lraf{d^t}
X^t\lra 0;
$$
where $d^0=\pi\mu$ and $M$ is the image of $d^0$.  We claim
$\Img(\mu)\in\add(T)$. In fact, we have a long exact sequence
$$
0\lra \Img(\mu)\lraf{\nu} X^1\lraf{d^1}\cdots\lraf{d^t} X^t\lra 0,
$$
where $\nu$ is the canonical inclusion. For each $1\leq i\leq t$,
since $X^i\in\add(T)$,  we have $\Img(d^i)\in\Gen(T)$. As we know,
$T^{\bot}=\Gen(T)$ for a tilting module $T$. Consequently, we see
that $\Ext_A^1(T,\Img(d^i))=0$ for any $1\leq i\leq t$. Note that
$\Img(d^t)=X^t\in\add(T)$. Thus we can easily show
$\Img(\mu)\in\add(T)$ by induction on $t$.

Finally, we shall prove $M\in\add(T)$. If $s=0$, then
$M=X^0\in\add(T)$. Suppose $s<0$. Since $\Img(\mu)\in\add(T)$ and
the sequence  $0\to A\to M\to \Img(\mu)\to 0$ is exact,  we know
that $\Ext_A^1(M,T)=0$ and $M$ has projective dimension at most $1$.
In addition, $\Img(d^{-1})$ is a quotient module of $X^{-1}$. It
follows that  $\Ext_A^1(M,\Img(d^{-1}))=0$, which implies that the
homomorphism $\pi$ splits. Thus $M\in\add(X^0)\subseteq\add(T)$.

Now we define $T_0=M$ and $T_1=\Img(\mu)$. Then the sequence $0\to
A\to T_0\lraf{\mu} T_1\to 0$ satisfies $T_i\in\add(T)$ for $i=0,1$.
Thus $T$ is a good tilting module, and the proof is completed.
$\square$

\medskip
{\it Remark.} Suppose that $G$ admits a fully faithful left adjoint
$j_!: \D{A}\to\D{B}$. Then there exists a {\rm{TTF}} triple
$\big(j_!(\D{A}),\Ker(G), H(\D{A})\big)$ in $\D{B}$ (see
\cite[Chapter I, Proposition 2.11]{BI} for details), where
$j_!(\D{A})$ and $H(\D{A})$ denote the images of $j_!$ and $H$,
respectively. By Lemma \ref{lem01}, we know that the derived
category $\D{B}$ is a recollement of the derived category $\D{A}$
and $\Ker(G)$. Since $T$ is good by Theorem \ref{th01}(2), it
follows from Theorem \ref{th01}(1) that $\Ker(G)$ is triangle
equivalent to the derived category $\D{C}$ of a ring $C$. Thus we
get a recollement of derived module categories as in Theorem
\ref{th01}(1).

\section{Applications to tilting modules arising from ring epimorphisms\label{6}}

In this section we apply our main result Theorem \ref{th01} to
tilting modules arising from ring epimorphisms. In this case we
shall describe the universal localization rings appearing in the
main result by coproducts defined by Cohn in \cite{cohn}. In fact,
our discussion in this section will be implemented in the general
setup of injective homological ring epimorphisms between arbitrary
rings, which is weaker than conditions of being tilting modules.

We start with recalling of some definitions.

Let $R_0$ be a ring with identity. An $R_0$-ring is a ring $R$
together with a ring homomorphism $\lambda_R: R_0\ra R$. An
$R_0$-homomorphism from an $R_0$-ring $R$ to another $R_0$-ring $S$
is a ring homomorphism $f: R\ra S$ such that $\lambda_S=\lambda_Rf$.
If $R_0$ is commutative and the image of $\lambda_R: R_0\ra R$ is
contained in the center $Z(R)$ of $R$, then we say that $R$ is an
$R_0$-algebra.

Recall that the coproduct of a family $\{R_i\mid i\in I\}$ of
$R_0$-rings with $I$ an index set is an $R_0$-ring $R$ together with
a family $\{\rho_i: R_i\ra R\mid i\in I \}$ of $R_0$-homomorphisms
such that, for  any $R_0$-ring $S$ with a family of
$R_0$-homomorphisms $\{\tau_i: R_i\ra S\mid i\in I\}$, there is a
unique $R_0$-homomorphism $\delta: R\ra S$ such that
$\tau_i=\rho_i\delta$ for all $i\in I$.

It is well-known that the coproduct of a family $\{R_i\mid i\in I\}$
of $R_0$-rings exists. In this case, we denote their coproduct by
$\sqcup_{R_0}R_{i}$. For example, the coproduct of the polynomial
rings $k[x]$ and $k[y]$ over $k$ is the free ring $k<x,y>$ in two
variables over $k$.  Note that $R_0\sqcup_{R_0}S=S=S\sqcup_{R_0}R_0$
for every $R_0$-ring $S$.

Let $R_i$ be an $R_0$-ring for $i=1,2$. We denote by
$B$ the matrix ring $\left(\begin{array}{cc} R_1 & R_1\otimes_{R_0}R_2\\
0 & R_2\end{array}\right).$ Let $e_1 = \left(\begin{array}{ll} 1 &0\\
0 & 0\end{array}\right)$, $e_2=\left(\begin{array}{ll} 0 &0\\
0 & 1\end{array}\right)\in B$, and $\varphi: Be_1\ra Be_2$ be
the map that sends $\left(\begin{array}{l} r_1 \\
0 \end{array}\right)$ to $\left(\begin{array}{c} r_1\otimes 1\\
0 \end{array}\right)$ for $r_1\in R_1$. Let $\rho_i: R_i\ra
R_1\sqcup_{R_0}R_2$ be the canonical $R_0$-homomorphism for $i=1,2.$

\medskip
The following lemma reveals a relationship between coproducts and
localizations.

\begin{Lem} {\rm\cite[ Theorem 4.10, p. 59]{Sch}} The universal localization $B_{\varphi}$
of $B$ at  $\varphi$ is equal to $M_2(R_1\sqcup_{R_0}R_2)$, the
$2\times 2$ matrix ring over the coproduct $R_1\sqcup_{R_0}R_2$ of
$R_1$ and $R_2$ over $R_0$. Furthermore, the corresponding ring
homomorphism
$\lambda_{\varphi}: B\to B_{\varphi}$ is given by $\left(\begin{array}{cc} r_1 & x_1\otimes x_2\\
0 & r_2\end{array}\right)\mapsto\left(\begin{array}{cc} (r_1)\rho_1
&(x_1)\rho_1(x_2)\rho_2\\
0 & (r_2)\rho_2 \end{array}\right)$ for  $r_i, x_i\in R_i$ with
$i=1,2$. \label{epi0}
\end{Lem}

The next result says, in some sense, that coproducts of rings
preserve universal localizations.

\begin{Lem}
Let $R_0$ be a ring, $\Sigma$ a set of homomorphisms between
finitely generated projective $R_0$-modules, and $\lambda_{\Sigma}:
R_0\to R_1:=(R_0)_{\Sigma}$ the universal localization of $R_0$ at
$\Sigma$. Then, for any $R_0$-ring $R_2$, the coproduct
$R_1\sqcup_{R_0}R_2$ is isomorphic to the universal localization
$(R_2)_{\Delta}$ of $R_2$ at the set $\Delta$, where
$\Delta:=\{R_2\otimes_{R_0}f\mid f\in\Sigma\}$. \label{epi01}
\end{Lem}

{\it Proof.} Let $R:=(R_2)_{\Delta}$, and let
$\lambda_{\Delta}:R_2\to R$ be the universal localization of $R_2$
at $\Delta$. Suppose that $\lambda_{R_2}:R_0\to R_2$ is the ring
homomorphism defining the $R_0$-ring $R_2$. Then $R$ is an $R_0$-
ring via the composition $\lambda_{R_2}\lambda_{\Delta}$ of
$\lambda_{R_2}$ with $\lambda_{\Delta}$. Moreover, we shall prove
that there is a unique $R_0$-ring homomorphism $\nu: R_1\ra R$, that
is, a ring homomorphism $\nu$ with
$\lambda_{R_2}\lambda_{\Delta}=\lambda_{\Sigma}\nu$. In fact, for
any $f:P_1\to P_0$ in $\Sigma$, the map
$R\otimes_{R_0}f:R\otimes_{R_0}P_1\to R \otimes_{R_0}P_0$ of
$R$-modules is an isomorphism because $R\otimes_{R_0}f\simeq
R\otimes_{R_2}(R_2\otimes_{R_0}f)$, where the later is an
isomorphism. Thus, by the property of the universal localization
$\lambda_{\Sigma}$, there is a unique ring homomorphism $\nu: R_1\ra
R$ such that $\lambda_{R_2}\lambda_{\Delta}=\lambda_{\Sigma}\nu$, as
desired.

Now, we show that $R$ together with the two ring homomorphisms
$\lambda_{\Delta}$ and $\nu$ satisfies the definition of coproducts,
and therefore $R_1\sqcup_{R_0}R_2$ is isomorphic to $R$.

Indeed, suppose that $S$ is an arbitrary $R_0$-ring with two
$R_0$-homomorphisms $\tau_i: R_i\ra S$ for $i=1,2$. Then
$\lambda_{\Sigma}\tau_1=\lambda_{R_2}\tau_2$. Further, since we have
$$S\otimes_{R_2}(R_2\otimes_{R_0}h)\simeq S\otimes_{R_0}h\simeq
S\otimes_{R_1}(R_1\otimes_{R_0}h),$$ and since $R_1\otimes_{R_0}h$
is an isomorphism  for any $h\in\Sigma$, we infer that
$S\otimes_{R_2}(R_2\otimes_{R_0}h)$ is an isomorphism for any $h\in
\Sigma$. It follows from the property of universal localizations
that there is a unique ring homomorphism $\delta:R\to S$ such that
$\tau_2=\lambda_{\Delta}\delta$. Clearly,
$\lambda_{\Sigma}\tau_1=\lambda_{\Sigma}\,\nu\delta$, and
$\tau_1=\nu\delta$ since $\lambda_{\Sigma}$ is a ring epimorphism.
Note that $\delta$ is also an $R_0$-ring homomorphism. Thus, there
is a unique $R_0$-homomorphism $\delta: R\ra S$ such that
$\tau_1=\mu\delta$ and $\tau_2=\lambda_{\Delta}\delta$. This shows
that the coproduct $R_1\sqcup_{R_0}R_2$ of $R_1$ and $R_2$ over
$R_0$ is isomorphic to $R$. $\square$

\medskip
Sometimes, coproducts can be interpreted as tensor products of
rings.

\begin{Lem}
Let $R_0$ be a commutative ring and let $R_i$ be an $R_0$-algebra
for $i=1,2$.  If  one of the  homomorphisms $\lambda_{R_1}:R_0\to
R_1$ and $\lambda_{R_2}:R_0\to R_2$  is a ring epimorphism, then the
coproduct $R_1\sqcup_{R_0}R_2$  is isomorphic to the tensor product
$R_1\otimes_{R_0}R_2$.

\label{epi00}
\end{Lem}

{\it Proof.} It is known that the tensor product
$R_1\otimes_{R_0}R_2$ of two rings $R_1$ and $R_2$ over $R_0$ has
the following universal property: If $f_i: R_i\ra R$ is a
homomorphism of $R_0$-rings for $i=1,2,$ such that
$(r_2)f_2(r_1)f_1=(r_1)f_1(r_2)f_2$ for all $r_i\in R_i$ with
$i=1,2$, then there is a unique ring homomorphism $f: R_1\otimes_{
R_0}R_2\ra R$ of $R_0$-rings that satisfies $(x_1\otimes
x_2)f=(x_1)f_1(x_2)f_2$ for $x_i\in R_i$ with $i=1,2$. In
particular, if $\lambda_1: R_1\ra R_1\otimes_{R_0}R_2$ is the map
given by $r_1\mapsto r_1\otimes 1$ for $r_1\in R_1$, and if
$\lambda_2: R_2\ra R_1\otimes_{R_0}R_2$ is the one given by
$r_2\mapsto 1\otimes r_2$ for $r_2\in R_2$, then $f_i=\lambda_if$
for $i=1,2.$

To prove Lemma \ref{epi00}, it suffices to show that, for any
$R_0$-homomorphisms $f_i: R_i\ra R$ for $i=1,2$,  the condition
$(r_2)f_2(r_1)f_1=(r_1)f_1(r_2)f_2$ holds true  for all $r_i\in R_i$
with $i=1,2$.

Assume that $\lambda_{R_1}:R_0\to R_1$ is a ring epimorphism. For
any element $y\in R_2$, we define two ring homomorphisms
$\theta_1:R_1\to \rm M_2(R) $ and  $\theta_{2}:R_1\to \rm M_2(R)$ as
follows:
$$(x)\theta_1=\left(\begin{array}{cc} (x)f_1 &0\\
0 & (x)f_1 \end{array}\right)$$ and
$$(x)\theta_{2}=\left(\begin{array}{cc} 1 &0\\
(y)f_2& 1
\end{array}\right)\left(\begin{array}{cc} (x)f_1 &0\\
0 & (x)f_1 \end{array}\right)\left(\begin{array}{cc} 1 &0\\
-(y)f_2& 1
\end{array}\right)=\left(\begin{array}{cc} (x)f_1 &0\\
(y)f_2(x)f_1-(x)f_1(y)f_2 & (x)f_1
\end{array}\right)$$ for  $x\in R_1$. Now, we verify
$\lambda_{R_1}\theta_1=\lambda_{R_1}\theta_2$. This is equivalent to
showing that, if $x=(r)\lambda_{R_1}$ with $r\in R_0$, then
$(y)f_2(x)f_1=(x)f_1(y)f_2 $. In fact, we always have
$$
(y)f_2(x)f_1=(y)f_2\big((r)\lambda_{R_1}\big)f_1=(y)f_2\big((r)\lambda_{R_2}\big)f_2
=\big(y(r)\lambda_{R_2}\big)f_2;
$$
$$
(x)f_1(y)f_2=\big((r)\lambda_{R_1}\big)f_1(y)f_2=\big((r)\lambda_{R_2}\big)f_2
(y)f_2=\big((r)\lambda_{R_2}y\big)f_2.
$$
Since $R_2$ is an $R_0$-algebra, it follows from
$\Img(\lambda_{R_2})\subseteq Z(R_2)$ that
$y(r)\lambda_{R_2}=(r)\lambda_{R_2}y$, and so
$(y)f_2(x)f_1=(x)f_1(y)f_2$ whenever $x=(r)\lambda_{R_1}$ with $r\in
R_0$. This shows $\lambda_{R_1}\theta_1=\lambda_{R_1}\theta_2$ and
$\theta_1=\theta_2$ since $\lambda_{R_1}: R_0\to R_1$ is a ring
epimorphism. Thus $(y)f_2(x)f_1=(x)f_1(y)f_2$ for any $x\in R_1$.
Note that $y$ is an arbitrary element of $R_2$. Hence
$(y)f_2(x)f_1=(x)f_1(y)f_2$ for any $x\in R_1$ and $y\in R_2$.

As a result, the tensor product $R_1\otimes_{R_0}R_2$ together with
the two ring homomorphisms $\lambda_i$ satisfies the definition of
coproducts, and we therefore have proved that the coproduct
$R_1\sqcup_{R_0}R_2$  is isomorphic to the tensor product
$R_1\otimes_{R_0}R_2$. Similarly, we can prove Lemma \ref{epi00}
under the assumption that $\lambda_{R_2}:R_0\to R_2$  is a ring
epimorphism. $\square$

From now on, let $\lambda: R\ra S$ denote an injective ring
homomorphism from $R$ to $S$. We define $B$ to be the endomorphism
ring of the $R$-module $S\oplus S/R$, and $S'$ the endomorphism ring
of the $R$-module $S/R$. Let $\pi$ stands for the canonical
surjective map $S\ra S/R$ of $R$-modules. Then we have an exact
sequence of $R$-modules:
$$(*)\quad 0\lra R\lra S\lraf{\pi} S/R\lra 0.$$

In the next lemmas, we collect some facts on ring epimorphisms.
\begin{Lem} Let $\lambda: R\ra S$ be an injective ring epimorphism with $\Tor_1^R(S,S)=0$.
Then,

$(1)$ an $R$-module $X$ belongs to $S\Modcat$ if and only if
$\Ext^i_R(S/R,X)=0$ for $i=0,1$.

$(2)$ Let $T:= S\oplus S/R$. Then
$$\End_R(T)\simeq\left(\begin{array}{lc} S & \Hom_R(S,S/R)\\
0 & \End_R(S/R)\end{array}\right). $$ Moreover, if $e_1$ and $e_2$
are the idempotent elements in $\End_R(T)$ corresponding to the
summands $S$ and $S/R$, respectively, then the homomorphism
$\pi^*:\End_R(T)e_1\ra \End_R(T)e_2$ induced from the canonical
surjection $\pi: S\ra S/R$ is given by $\left(\begin{array}{l} s \\
0 \end{array}\right)\mapsto \left(\begin{array}{c} (x\mapsto (xs)\pi\\
0 \end{array}\right)$ for $s,x\in S$. \label{epi1}
\end{Lem}

{\it Proof.} (1) follows from \cite{GL}. For (2), it follows from
(1) that $\Hom_R(S/R,S)=0$. By applying $\Hom_R(-,S)$ to the exact
sequence ($*$), we get $\Hom_R(S,S)\simeq \Hom_R(R,S)\simeq S$.
$\square$

\begin{Lem} Suppose that $\lambda: R\ra S$ is an injective ring epimorphism
with $\Tor_1^R(S,S)=0$.

$(1)$ The right multiplication map $\mu:R\to S'$ defined by
$r\mapsto(x\mapsto xr)$ for $r\in R$ and $x\in S/R$, is a ring
homomorphism. Consequently, $S'$ can be regarded as an $R$-ring via
the map $\mu$. Further, $\mu$ is an isomorphism if and only if
$\Ext^i_R(S,R)=0$ for $i=0,1$.

$(2)$ There is an isomorphism $\theta: S\otimes_RS'\simeq
\Hom_R(S,S/R)$ of $S$-$S'$-bimodules such that $1\otimes 1$ is
mapped to the canonical surjection $\pi: S\ra S/R$.

$(3)$  There is an exact sequence of $R$-$S'$-modules:
$$0\ra S'\lraf{\lambda'} S\otimes_RS'\lraf{\pi\otimes S'}(S/R)\otimes_RS'\ra 0,$$
where the map $\lambda'$ is defined by $f\mapsto
 1\otimes f$ for any $f\in S'$. Moreover, the evaluation map
$\psi:(S/R)\otimes_RS'\to S/R$ defined by $y\otimes g\mapsto (y)g$
for $y\in S/R$ and $g\in S'$, is an isomorphism of
$R$-$S'$-bimodules.

$(4)$ If $\lambda: R\ra S$ is  homological, then $\Tor^R_i(S,S')=0$
for any $i>0$.

$(5)$ If $R$ is commutative, then so is $S'$.
\label{epi2}
\end{Lem}

{\it Proof.} $(1)$  It is easy to checked that the right
multiplication map $\mu$ is a ring homomorphism since  $S/R$ is an
$R$-$R$-bimodule. Clearly, $\mu$ is injective if and only if
$\Hom_R(S,R)=0$. For the condition of $\mu$ being surjective, we use
the following exact sequence: $$ 0\lra \Hom_R(S,R)\lra
\Hom_R(S,S)\lra \Hom_R(S,S/R)\lra \Ext^1_R(S,R)\lra \Ext^1_R(S,S),
$$ where $\Ext^1_R(S,S)=0$ by Lemma \ref{epi1}(1). Thus (1) follows.

$(2)$ Note that a ring homomorphism is an epimorphism if and only if
the multiplication map $S\otimes_RS\ra S$ is an isomorphism as
$S$-$S$-bimodules. Since $\lambda$ is injective, it follows from the
exact sequence ($*$) that we have a long exact sequence of
$S$-$R$-bimodules:
$$ 0\lra \Tor_1^R(S,S)\lra \Tor_1^R(S,S/R)\lra S\otimes_R
R\lraf{1\otimes \lambda} S\otimes_RS\lra S\otimes_R(S/R)\lra 0.$$
 Since $\Tor_1^R(S,S)=0$ and
$1\otimes_R\lambda$ is an isomorphism of $S$-$R$-modules, we have
$S\otimes_R(S/R)=0=\Tor_1^R(S,S/R)$.

Now, by applying $\Hom_R(-,S/R)$ to ($*$), we can get another exact
sequence of $R$-$\End_R(S/R)$-bimodules:
$$(**)\quad 0\lra \Hom_R(S/R,S/R)\lra \Hom_R(S,S/R)\lra \Hom_R(R,S/R). $$
One can check that the last homomorphism in the above sequence is
surjective because each element $x+R$ in $S/R$ gives rise to at
least one homomorphism from the $R$-module $S$ to the $R$-module
$S/R$ by $s\mapsto sx+R$ for $s\in S$. This yields the following
exact sequence of $S$-$\End_R(S/R)$-bimodules:
$$0\lra S\otimes_R\Hom_R(S/R,S/R)\lra S\otimes_R\Hom_R(S,S/R)\lra S\otimes_R(S/R)\lra 0, $$
which shows that $S\otimes_R\Hom_R(S/R,S/R)\lraf{\sim}
S\otimes_R\Hom_R(S,S/R)$. Clearly, under this isomorphism the
element $1\otimes_R 1$ in $S\otimes \Hom_R(S/R,S/R)$ is sent to
$1\otimes \pi$. Since the multiplication map: $S\otimes_RS\lra S$ is
an isomorphism of $S$-$S$-bimodules, we see that the multiplication
map: $S\otimes_RX\ra X$ is an isomorphism for every $S$-module $X$.
Clearly, $\Hom_R(_RS_S, S/R)$ is an $S$-module. So, it follows that
$S\otimes_R\Hom_R(S,S/R)\ra \Hom_R(S,S/R)$ is an isomorphism under
which $1\otimes \pi$ is sent to $\pi$. As a result,  the map
$\theta:S\otimes_RS'\ra \Hom_R(S,S/R)$ defined by $s\otimes f\mapsto
\big(t\mapsto (ts)(\pi f)\big)$ for $s,t\in S$ and $f\in S'$, is an
isomorphism of $S$-$S'$-bimodules. Clearly, under this isomorphism,
the element $1\otimes 1$ in $S\otimes_RS'$ is sent to $\pi$.

$(3)$ Applying $-\otimes_RS'$ to the sequence $(*)$ and identifying
$R\otimes_RS'$ with $S'$, we then obtain the following right exact
sequence of $R$-$S'$-bimodules:
$$ (\spadesuit)\quad S'\lraf{\lambda'} S\otimes_RS'\lraf{\pi\otimes S'}(S/R)\otimes_RS'\ra 0,$$
where the map $\lambda'$ is defined by $f\mapsto
 1\otimes f$ for any $f\in S'$. Combining this sequence with $(**)$,
one can check  that the following diagram of $R$-$S'$-bimodules is
exact and commutative:
$$\xymatrix{
&S'\ar@{=}[d]\ar[r]^-{\lambda'}&S\otimes_RS'\ar[d]_-{\theta}\ar[r]^-{\pi\otimes
S'}&(S/R)\otimes_RS'\ar[d]_-{\psi}\ar[r]&0\\
0\ar[r]&\Hom_R(S/R,S/R)\ar[r]^-{\pi_*}&\Hom_R(S,S/R)\ar[r]&S/R\ar[r]&0
}$$ is commutative, where $\psi$ is the evaluation map, and where
$\Hom_R(R,S/R)$ is identified with $S/R$ as $R$-$S'$-bimodules.
Since $\theta$ is an isomorphism, we infer that $\lambda'$ is
injective, and that $\psi$ is an isomorphism of $R$-$S'$-bimodules.

$(4)$ Suppose that $\lambda$ is an injective homological ring
epimorphism. Let $i$ be any positive integer. Then
$\Tor^R_i(S,S)=0$.
 Recall that we have proved that
$S\otimes_R(S/R)=0=\Tor_1^R(S,S/R)$ in $(2)$. Thus, by applying the
tensor functor $S\otimes_R-$ to the canonical sequence $(*)$, we
conclude that $\Tor^R_i(S,S/R)=0$. By $(3)$, we know that
$(S/R)\otimes_RS'\simeq S/R$ as left $R$-modules. Thus
$\Tor^R_i(S,(S/R)\otimes_RS')=0$. Since $S\otimes_RS'$ is a left
$S$-module, it follows from Lemma \ref{2.1}(2) that
$\Tor^R_i(S,S\otimes_RS')=0$. Now,  applying the tensor functor
$S\otimes_R-$ to the exact sequence $(\spadesuit)$, we obtain
$\Tor^R_i(S,S')=0$.

(5) Since $R$ is commutative, the tensor product $S\otimes_RS'$ of
$S$ and $S'$ over $R$ is a ring, which is well-defined. By Lemma
\ref{epi2}(3), there exists an exact sequence of $R$-$S'$-modules:
$$0\ra S'\lraf{\lambda'} S\otimes_RS'\lraf{\pi\otimes S'}(S/R)\otimes_RS'\ra 0.$$
Since $\lambda'$ is a ring homomorphism, the ring $S\otimes_RS'$ can
be considered an $S'$-$S'$-bimodule via $\lambda'$, and therefore,
$(S/R)\otimes_RS'$ can also be regarded as an $S'$-$S'$-bimodule. In
addition, by Lemma \ref{epi2}(3), the evaluation map $\psi:
(S/R)\otimes_RS'\to S/R$ defined by $y\otimes g\mapsto (y)g$ for any
$y\in S/R$ and $g\in S'$, is an isomorphism of $R$-$S'$-bimodules.
Since the image of $(y)g\otimes 1$ under $\psi$ is also equal to
$(y)g$, we have $(y)g\otimes 1=y\otimes g$ in $(S/R)\otimes_RS'$.
Consequently, for any $f,g\in S'$ and $y\in S/R$, we get
 $y\otimes fg=f(y\otimes g)=f((y)g\otimes
1)=(y)g\otimes f$ in $(S/R)\otimes_RS'$, where the first  and third
equalities follow from the left $S'$-module structure of
$(S/R)\otimes_RS'$. This yields that $(y)fg=\big(y\otimes
fg\big)\psi=\big((y)g\otimes f\big)\psi=(y)gf$ in $S/R$. Thus
$fg=gf$. Since $f$ and $g$ are arbitrary elements in $S'$, we see
that $S'$ is a commutative ring. $\square$

\medskip
As a consequence of Theorem \ref{th01}, we have the following
corollary.

\begin{Koro} $(1)$ Let $\lambda: R\ra S$ be an injective
ring epimorphism such that $\Tor_1^R(S,S)=0$. If $_RS$ has
projective dimension at most one, then there is a recollement of
derived module categories:

$$\xymatrix@C=1.2cm{\D{S\sqcup_RS'}\ar[r]&\D{B}\ar[r]
\ar@/^1.2pc/[l]\ar@/_1.2pc/[l]
&\D{R}\ar@/^1.2pc/[l]\ar@/_1.2pc/[l]},\vspace{0.3cm}$$ where
$S\sqcup_RS'$ is the coproduct of $S$ and $S'$ over $R$.

$(2)$ Let R be a ring,  $\Sigma$ a left Ore set of regular elements
of R, and $S:=\Sigma^{-1}R$ the localization of $R$ at $\Sigma$. If
$_RS$ has projective dimension at most one, then the recollement in
$(1)$ exists. \label{ringepi}
\end{Koro}

{\it Proof.} (1) Now, let $R_0=R$, $R_1=S$, $R_2=\End_R(S/R)$,
$T:=S\oplus S/R$, and $B=\End_R(S\oplus S/R)$. It is well known
that, under the above assumptions, $T:=S\oplus S/R$ is a good
tilting $R$-module (see \cite{HJ}). By Lemma \ref{epi2}, the map
$\varphi$ in Lemma \ref{epi0} is precisely the map $\pi^*$ in Lemma
\ref{epi1}(2) under the identification of $\theta$ in Lemma
\ref{epi2}. Thus the localization of $B$ at $\pi^*$ is isomorphic to
the $2\times 2$ matrix ring over the coproduct of $S$ and $S'$ over
$R$ by Lemma \ref{epi0}. Thus Corollary \ref{ringepi}(1) follows
from Theorem \ref{th01}(1).

(2) follows from (1). $\square$

The tilting module $S\oplus \; S/R$ in Corollary \ref{ringepi} has
an equivalent form (see \cite[Theorem 2.10]{HJ}), by which we can
restate Corollary \ref{ringepi} in the following form. Here we
present explicitly the $R$-ring homomorphisms which will be used for
later calculations in Section \ref{sect8}.

\begin{Koro}  Let $R$ be a ring, and let $_RT$ be a tilting $R$-module with an exact sequence
$$ 0\lra R\lra T_0\lra T_1\lra 0$$ of $R$-modules such that $T_i\in \Add(T)$ and $\Hom_R(T_1,T_0)=0$.
Set $S:=\End_R(T_0)$, $S':=\End_R(T_1)$ and $B:=\End_R(T_0\oplus
T_1)$. Then there is the following recollement:

$$\xymatrix@C=1.2cm{\D{S\sqcup_RS'}\ar[r]&\D{B}\ar[r]
\ar@/^1.2pc/[l]\ar@/_1.2pc/[l]
&\D{R}\ar@/^1.2pc/[l]\ar@/_1.2pc/[l]},\vspace{0.3cm}$$ where
$S\sqcup_RS'$ is the coproduct of $S$ and $S'$ over $R$.
\label{6.15}
\end{Koro}

{\it Proof.} First of all, we show that $S$ and $S'$ can be regarded
as $R$-rings, namely, we construct two ring homomorphisms
$\lambda:R\to S$ and $\mu:R\to S'$ (see Lemma \ref{epi2}(1)). For
any $r\in R$, we denote by $\rho_r:R\to R$ the right multiplication
map by the element $r$. It follows from $\Hom_R(T_1,T_0)=0=
\Ext_R^1(T_1,T_0)$ that there exists a unique homomorphism $f:T_0\to
T_0$ and therefore a unique homomorphism $g:T_1\to T_1$ such that
the following exact diagram
$$\xymatrix{
 0\ar[r]& R\ar[r]\ar[d]_-{\rho_r}&
T_0\ar[r]\ar@{-->}_-{f}[d]&T_1\ar[r]\ar@{-->}_-{g}[d]&0\\
0\ar[r]& R\ar[r]&T_0\ar[r]&T_1\ar[r]&0 }
$$
commutes. Now, we define $\lambda:R\to S$  and $\mu:R\to S'$ by
sending $r$ to $f$ and sending $r$ to $g$, respectively. One can
check directly that $\lambda$ is injective, and that $\lambda$ and
$\mu$ are ring homomorphisms. Furthermore, by the  proof of
\cite[Theorem 2.10]{HJ}, there are isomorphisms $\varphi:T_0\to S$
and $\psi:T_1\to S/R$ of $R$-modules, such that the following exact
diagram of $R$-modules
$$\xymatrix{
 0\ar[r]& R\ar[r]\ar@{=}[d]&
T_0\ar[r]\ar_-{\varphi}[d]^-{\wr}&T_1\ar[r]\ar_-{\psi}[d]^-{\wr}&0\\
0\ar[r]& R\ar[r]^-{\lambda}&S\ar[r]^-{\pi}&S/R\ar[r]&0 }
$$
is commutative. Now, Corollary \ref{6.15} follows from Corollary
\ref{ringepi}. $\square$.

\medskip
{\it Remarks.} 
(1) In general, Corollary \ref{ringepi}(1) supplies us a class of
recollements which cannot be obtained from the structure of
triangular matrix rings.

If $\Ext^i_R(S,R)=0$ for $0\le i\le 1$ in Corollary
\ref{ringepi}(1), then $S'\simeq R$ by Lemma \ref{epi2}(1), and
therefore $S\sqcup_RS'\simeq S\sqcup_R R=S$. Even in this case, the
recollement in Corollary \ref{ringepi}(1) is not equivalent to the
one induced from the triangular matrix ring (see Lemma
\ref{epi1}(2)) since they are induced by non-equivalent homological
ring epimorphisms $B\ra M_2(S)$ and $B\ra S$, respectively. See
Section \ref{sect7} below for more examples.

\medskip
(2) In Corollary \ref{ringepi}(1), the condition that $``_RS$ has
projective dimension at most one" ensures the category
$\D{B}_{\Sigma^{\perp}}$ in Corollary \ref{3.5} can be replaced by
the derived module category of a ring. However, this condition is
not necessary for getting such a recollement. In fact, we have the
following result.

Let $\lambda: R\ra S$ be an injective ring epimorphism such that
$\Tor_1^R(S,S)=0$. Suppose that  the right multiplication map $\mu:
R\ra S'$ defined in Lemma \ref{epi2}(1) is an isomorphism of rings.
Then  the universal localization $\lambda_{\pi^*}:B\to B_{\pi^*}$ of
$B$ at $\pi^*$ is  homological. In particular,
$\D{B}_{\{\pi^*\}^{\perp}}$ in Corollary \ref{3.5} can be replaced
by $\D{S}$.

{\it Proof.} Combining Lemmata \ref{epi0} and \ref{epi1}(2) with
\ref{epi2}(2), we know that  $B\simeq\left(\begin{array}{lc} S & S\\
0 & R\end{array}\right)$ and $B_{\pi^*}\simeq M_2(S)$ as rings.
Under these isomorphisms,  the universal localization
$\lambda_{\pi^*}:B\to B_{\pi^*}$ is equivalent to the canonical
injection $\chi:R_1:=\left(\begin{array}{lc} S & S\\
0 & R\end{array}\right)\to M_2(S)$ induced by the injection
$\lambda:R\to S$. Clearly,  $M_2(S)$ is projective as a right
$R_1$-module. Thus $\chi$ is homological, and consequently,
$\lambda_{\pi^*}$ is homological, and $\D{S}$ is triangle equivalent
to $\D{B}_{\{\pi^*\}^{\perp}}$ by Proposition \ref{th2}. $\square$

\medskip
Let us give a concrete example satisfying all conditions in
the remark (2).

Let $R=\left(\begin{array}{ccc} k & 0          & 0\\
                                  k[x]/(x^2) & k          & 0 \\
                                  k[x]/(x^2) & k[x]/(x^2) & k
                         \end{array}\right),$\vspace{0.2cm}
where $k$ is a field and $k[x]$ is the polynomial algebra over $k$
in one variable $x$, and let $S$ be the $3$ by $3$ matrix ring
$M_3(k[x]/(x^2))$. Then the inclusion $\lambda$ of $R$ into $S$ is a
universal localization of $R$, so it is a ring epimorphism. Further,
we have $\Tor_1^R(S,S)=0\neq\Tor_2^R(S,S)$ (see \cite{nrs}). Thus it
is not homological. So, $_RS$ cannot have projective dimension less
than or equal to one. But one can check that $\mu$ defined in Lemma
\ref{epi2} is an isomorphism.

This example also shows that Proposition \ref{new} below may be
false if the injective ring epimorphism $\lambda : R\ra S$ is not
homological.

\medskip
(3) Under the conditions of Remark (2), one can get another pattern
of recollements, namely, we have the following result.

\medskip
Let $\lambda: R\ra S$ be an injective ring epimorphism such that
$\Tor_1^R(S,S)=0$. Suppose that  the right multiplication map $\mu:
R\ra \End_R(S/R)$ defined in Lemma \ref{epi2}(1) is an isomorphism
of rings. Then there is a recollement of derived module categories
of $S$, $B$ and $R$:

$$\xymatrix@C=1.2cm{\D{R}\ar[r]&\D{B}\ar[r]
\ar@/^1.2pc/[l]\ar@/_1.2pc/[l]
&\D{S}\ar@/^1.2pc/[l]\ar@/_1.2pc/[l]}.\vspace{0.3cm}$$

\medskip
{\it Proof}.  Note that the sequence $0\lra
R\lraf{\lambda}S\lraf{\pi}S/R\lra 0$ is an $\add(_RS)$-split
sequence in $R$-Mod.  By \cite[Theorem 3.5]{hx2}, we conclude that
$B$ is derived equivalent to the endomorphism ring $\End_R(R\oplus
S)$ which is isomorphic to the
triangular matrix ring $\left(\begin{array}{lc} R & S\\
0 & S\end{array}\right)$. Consequently, we can get the above
recollement of derived module categories of $S$, $B$ and $R$ by the
structure of triangular matrix rings. $\square$

\medskip
A variation of Corollary \ref{ringepi}(1) is the following
proposition in which we relax the condition of being tilting
modules, and require ring epimorphisms to be homological.

\bigskip
Let $\lambda: R\ra S$ be an injective ring homomorphism between
rings $R$ and $S$. We consider $S':=\End_R(S/R)$ as an $R$-ring via
$\mu$ defined in Lemma \ref{epi2}. Furthermore, let $\rho: S\ra
S\sqcup_{R}S'$ and $\rho':S'\ra S\sqcup_{R}S'$  be the canonical
$R$-homomorphisms in the definition of coproducts of $R$-rings.

\begin{Prop}  If $\lambda: R\to S$ is an injective homological ring epimorphism, then the following assertions  are equivalent:

 $(1)$ The universal
localization $\lambda_{\pi^*}:B\to B_{\pi^*}$ of $B$ at $\pi^*$ is
homological.

$(2)$ The ring homomorphism $\rho':S'\ra S\sqcup_{R}S'$ is
homological.

In particular, if one of the above assertions holds, then there is a
recollement of derived module categories:
$$\xymatrix@C=1.2cm{\D{S\sqcup_RS'}\ar[r]&\D{B}\ar[r]
\ar@/^1.2pc/[l]\ar@/_1.2pc/[l]
&\D{R}\ar@/^1.2pc/[l]\ar@/_1.2pc/[l]}.$$
\label{new}
\end{Prop}

{\it Proof.} Recall that $S'$ is an $R$-ring via the right
multiplication map $\mu:R\to S'$, defined by $r\mapsto(x\mapsto xr)$
for $r\in R$ and $x\in S/R$ (see Lemma \ref{epi2}(1)). Then it
follows from the definition of coproducts of rings that
$\lambda\rho=\mu\rho':R\to S\sqcup_RS'$.

\textbf{Step $(1)$}.  We claim that, for any $S\sqcup_RS'$-module
$W$, if we regard $W$ as a left $S'$-module via the ring
homomorphism $\rho'$, then $(S\otimes_RS')\otimes_{S'}W\simeq W$ as
$S$-modules, and $\Tor^{S'}_i(S\otimes_RS', W)=0$ for any $i>0$.

To prove this, we fix a projective resolution  $\cpx{Q}$ of $S_R$:
$$\cdots
\lra Q^n\lra Q^{n-1}\lra\cdots \lra Q^1\lra Q^0 \lra S_R\lra 0$$
with $Q^i$ projective right $R$-modules. By Lemma \ref{epi2}(4), we
have $\Tor^R_j(S,S')=0$ for any $j>0$. It follows that  the complex
$\cpx{Q}\otimes_RS'$ is a projective resolution of the right
$S'$-module $S\otimes_RS'$. Note that we have the following
isomorphisms of complexes of abelian groups:
$$\big(\cpx{Q}\otimes_RS'\big)\otimes_{S'}W\simeq
\cpx{Q}\otimes_R\big(S'\otimes_{S'}W\big)\simeq\cpx{Q}\otimes_RW.$$
This implies that $\Tor^{S'}_i(S\otimes_RS', W)\simeq\Tor^{R}_i(S,
W)$ for any $i>0$. Clearly, $W$ admits an $S$-module structure via
the map $\rho$. Moreover, it follows from $\lambda\rho=\mu\rho'$
that the $R$-module structure of $W$ endowed via the ring
homomorphism $\mu\rho'$ is the same as that endowed via the ring
homomorphism $\lambda\rho$. Then, by Lemma \ref{2.1}, we conclude
that $S\otimes_RW\simeq W$ as $S$-modules, and that $\Tor^R_i(S,
W)=0$ for $i>0$. Therefore, $\Tor^{S'}_i(S\otimes_RS',W)=0$ for
$i>0$. Note that $(S\otimes_RS')\otimes_{S'}W\simeq S\otimes_RW$ as
$S$-modules. As a result, we have $(S\otimes_RS')\otimes_{S'}W\simeq
W.$ This finish the proof of Step (1).

 \textbf{Step $(2)$}.  We shall prove that $B_{\pi^*}$ is Morita
equivalent to the ring $S\sqcup_RS'$.

By Lemmata \ref{epi1}(2) and \ref{epi2}(2), there are isomorphisms
of rings $$B:=\End_R(T)
\simeq\left(\begin{array}{lc} S & \Hom_R(S,S/R)\\
0 & \End_R(S/R)\end{array}\right)\simeq C:=\left(\begin{array}{lc} S & S\otimes_RS'\\
0 & S'\end{array}\right),$$ where the second isomorphism sends
$\left(\begin{array}{cc} 0 & \pi\\
0 & 0\end{array}\right)$
to $\left(\begin{array}{cc} 0 & 1\otimes 1\\
0 & 0\end{array}\right)$. Let $e_1:= \left(\begin{array}{ll} 1 &0\\
0 & 0\end{array}\right)$, $e:=\left(\begin{array}{ll} 0 &0\\
0 & 1\end{array}\right)\in C$, and
$\varphi: Ce_1\to Ce$ to be the map that sends $\left(\begin{array}{l} s \\
0 \end{array}\right)$ to $\left(\begin{array}{c} s\otimes 1\\
0 \end{array}\right)$ for $s\in S$. Then $\pi^*$ corresponds to
$\varphi$ under the isomorphism $B\simeq C$. Let
$\lambda_{\varphi}:C\to C_{\varphi}$ be the universal localization
of $C$ at $\varphi$. Then, we have $B_{\pi^*}\simeq C_{\varphi}$.
Note that $\lambda_{\pi^*}$ is homological if and only if
$\lambda_{\varphi}$ is homological.

By Lemma \ref{epi0}, we know that $C_{\varphi}=M_2(S\sqcup_RS')$,
the $2\times 2$ matrix ring over $S\sqcup_RS'$, and that the
corresponding ring epimorphism $\lambda_{\varphi}: C\to
M_2(S\sqcup_RS')$ is given by
$\left(\begin{array}{cc} s & t\otimes f\\
0 & g\end{array}\right)\mapsto\left(\begin{array}{cc} (s)\rho
&(t)\rho(f)\rho'\\
0 & (g)\rho' \end{array}\right)$ for  $s,t\in S$ and $f,g\in S'$.
Hence, $B_{\pi^*}\simeq \Lambda$, and therefore $B_{\pi^*}$ is
Morita equivalent to $S\sqcup_RS'$.

\textbf{Step $(3)$}. We shall prove that the ring homomorphism
$\lambda_{\varphi}: C\to M_2(S\sqcup_RS')$ is homological if and
only if so is the ring homomorphism $\rho':S'\to S\sqcup_RS'$.

Before starting our proof, we  mention a general result:  if
$F:\mathcal {C}\to \mathcal{E}$ is an exact functor between abelian
categories $\mathcal C$ and $\mathcal E$, then $F$ can be extended
to a canonical triangle functor
$\overline{F}:\D{\mathcal{C}}\to\D{\mathcal{E}}$, which sends the
complex  $\cpx{X}:=\big(X^i, d_X^i \big)_{i\in\mathbb{Z}}$ over
$\mathcal{C}$ to the complex $\ol{F}(\cpx{X}):= \big(F(X^i),
F(d_X^i)\big)_{i\in\mathbb{Z}}$ over $\mathcal{E}$. This is due to
the fact that $F$ preserves quasi-isomorphisms. Since $\ol{F}$ is
completely determined by $F$, we may write $F$ for $\ol{F}$.

Set $\Gamma:=S\sqcup_RS'$, $\Lambda :=M_2(\Gamma)$ and
$e':=(e)\lambda_{\varphi}\in\Lambda$. Then
$e'=(e')^2,\,\End_{\Lambda}(\Lambda e')\simeq \Gamma$ and
$\End_{C}(Ce)\simeq S'$. Observe that $\Lambda e'$ is a projective
generator for $\Lambda$-Mod. Then, by Morita theory, we know that
the tensor functor $ e'\Lambda\otimes_{\Lambda}-:\Lambda\Modcat\lra
\Gamma\Modcat$ is an equivalence of module categories, which can be
extended to a canonical triangle equivalence from $\D{\Lambda}$ to
$\D{\Gamma}$.

Note that  $eC\otimes{_C}\Lambda\simeq e\Lambda$ as
 $S'$-$\Lambda$-bimodules, where the left $S'$-module structure of $e'\Lambda$
is induced by the ring homomorphism $\rho':S'\to \Gamma$. It follows
that  the following diagram of functors between module categories
$$\xymatrix{
\Lambda\Modcat\ar^-{e'\Lambda\otimes_{\Lambda}-}[rr]\ar_-{(\lambda_{\varphi})_*}[d]
&&\Gamma\Modcat\ar^{(\rho')_{*}}[d] \\
C\Modcat\ar^{eC\otimes_{C}-}[rr] &&S'\Modcat}
$$
is commutative, where  $(\lambda_{\varphi})_*$ and $(\rho')_{*}$
stand for the restriction functors induced by the ring homomorphisms
$\lambda_{\varphi}:C\to\Lambda$ and $\rho':S'\to\Gamma$,
respectively. Since all of the functors appearing in the diagram are
exact, we can  form the following commutative diagram of functors
between derived module categories:
$$\xymatrix{
\D{\Lambda}\ar^-{e'\Lambda\otimes_{\Lambda}-}[rr]\ar_-{(\lambda_{\varphi})_*}[d]
&&\D{\Gamma}\ar^{(\rho')_{*}}[d] \\
\D{C}\ar^{eC\otimes_{C}-}[rr] &&\D{S'},}
$$
where the upper tensor functor $e'\Lambda\otimes_{\Lambda}-$ is a
triangle equivalence.

From the triangular structure of $C$ it follows 
that there is a recollement of derived module categories:
$$\xymatrix{\D{S}\ar[r]&\D{C}\ar^-{eC\otimes_C-}[r]\ar@/^1.4pc/[l]\ar@/_1.4pc/[l]
&\D{S'}\ar@/^1.4pc/[l]\ar_-{Ce\otimesL_{S'}-}@/_1.4pc/[l]}\vspace{0.3cm}.$$
In particular, the pair
$\big({\rm{Tria}}(Ce),{\rm{Tria}}(Ce_1)\big)$ is a torsion pair in
$\D{C}$, and the  functor $eC\otimes_C-$ induces a triangle
equivalence between ${\rm{Tria}}(Ce)$ and $\D{S'}$.

We claim that the image of the restriction functor
$(\lambda_{\varphi})_*$ belongs to ${\rm{Tria}}(Ce)$. This implies
that, for any complexes
$\cpx{X},\cpx{Y}\in\Img\big((\lambda_{\varphi})_*\big)$, we have
$$\Hom_{\D{C}}(\cpx{X},\cpx{Y})\lraf{\sim}\Hom_{\D{S'}}(e\cpx{X},
e\cpx{Y}),$$ where the  functor $eC\otimes_C-$ is identified with
the left multiplication functor by the element $e$. Clearly, the
functor $(\lambda_{\varphi})_*$  can commute with small coproducts
since it admits a right adjoint. In addition,
$\D{\Lambda}=\rm{Tria}(\Lambda e')$. Therefore, to prove this claim,
it suffices to prove  $\Lambda e'\in {\rm Tria}(Ce)$. This is
equivalent to showing that $Ce\otimesL_{S'}e\Lambda e'\simeq \Lambda
e'$ in $\D{C}$.

Set $M:=S\otimes_RS'$. Note that $C$ is a triangular matrix ring.
Recall that any $C$-module can be expressed in the form of the
triple $(X,Y,h)$ with $X\in S'\Modcat,\,Y\in S\Modcat$ and
$h:M\otimes_{S'}X\to Y$ an $S$-morphism. The morphisms between two
modules $(X,Y,h)$ and $(X',Y',h')$ are pairs of morphisms $(\alpha,
\beta)$, where $\alpha:X\to X'$ and $\beta:Y\to Y'$ are
homomorphisms in $S'$-Mod and $S$-Mod, respectively, such that
$h\beta=(M\otimes_{S'}\alpha )h'$.

 With these interpretations, we rewrite $\Lambda
e'=(\Gamma, \,\Gamma, u)\in C\Modcat$, where
$u:M\otimes_{S'}\Gamma\to \Gamma$ is defined by $(s\otimes
f)\otimes\gamma\mapsto (s)\rho(f)\rho'\gamma\;$ for $s\in S, f\in
S'$ and $\gamma\in\Gamma$. Then $e\Lambda e'\simeq \Gamma$ as left
$S'$-modules. Clearly, $Ce\simeq M\oplus S'$ as right $S'$-modules.
Recall that we have proved in Step $(1)$ that $u$ is an isomorphism
of $S$-modules  and $\Tor^{S'}_i(M, \Gamma)=0$ for any $i>0$. It
follows that $\Tor^{S'}_i(Ce, e\Lambda e')=0$. Then we get the
following isomorphisms in $\D{C}$:
$$Ce\otimesL_{S'}e\Lambda e'\simeq Ce\otimes_{S'}e\Lambda
e'\simeq(\Gamma,\,M\otimes_{S'}\Gamma,\,1)\simeq\Lambda e'.$$ Thus
$\Lambda e'\in {\rm Tria}(Ce)$,  and we have finished the claim that
the image of the restriction functor $(\lambda_{\varphi})_*$ belongs
to ${\rm{Tria}}(Ce)$.

 With the above preparations, we now can prove
that the ring homomorphism $\lambda_{\varphi}: C\to \Lambda$ is
homological if and only if so is the ring homomorphism
$\rho':S'\to\Gamma$.

In fact, this can be concluded from the following commutative
diagram of
 functors between triangulated categories:
$$\xymatrix{
&\ar[ld]_-{(\lambda_{\varphi})_*}\D{\Lambda}\ar^-{e'\Lambda\otimes_{\Lambda}-}_-{\simeq}[rr]
\ar^-{(\lambda_{\varphi})_*}[d]
&&\D{\Gamma}\ar^{(\rho')_{*}}[d] \\
\D{C}&\ar@{_{(}->}[l]\;{\rm{Tria}}(Ce)\ar^-{eC\otimes_{C}-}_-{\simeq}[rr]
&&\D{S'},}
$$
which implies that $(\lambda_{\varphi})_*$ is fully faithful if and
only if so is $(\rho')_{*}$. It is known that $\lambda_{\varphi}:
C\to \Lambda$ (respectively, $\rho':S'\to\Gamma$) is homological if
and only if $(\lambda_{\varphi})_*$ (respectively, $(\rho')_{*}$) is
fully faithful.

Thus, we have proved that $\lambda_{\pi^*}:B\to B_{\pi^*}$ is
homological if and only if  $\rho':S'\ra S\sqcup_{R}S'$ is
homological. This finishes the proof of the first part of
Proposition \ref{new}. Clearly, the second part of Proposition
\ref{new} follows directly from Proposition \ref{th2}. $\square$

\medskip
Under all conditions in Corollary \ref{ringepi}(1), we see that both
$\lambda$ and $\lambda_{\pi^*}$ are homological, and therefore Lemma
\ref{new} implies Corollary \ref{ringepi}(1). However, for an
injective homological ring epimorphism $\lambda: R\ra S$, the
projective dimension of $_RS$ may not be at most one in general (see
the example at the end of Corollary \ref{new2} below). So, from this
point of view, Lemma \ref{new} may be regarded as a generalization
of Corollary \ref{ringepi}(1).

\medskip
Combining Lemma \ref{epi01} with Proposition \ref{new}, we get the
following criterion for $\lambda_{\pi^*}$ to be homological.

\begin{Koro}
Let $\Sigma$ be a set of homomorphisms between finitely generated
projective $R$-modules. Suppose that the universal localization
$\lambda_{\Sigma}:R\to R_{\Sigma}$ is an injective homological ring
homomorphism. Set $S:=R_{\Sigma}$ and $\lambda:=\lambda_{\Sigma}$,
and $\Gamma:=\{S'\otimes_{R}f\mid f\in\Sigma\}$ . Then the universal
localization $\lambda_{\pi^*}:B\to B_{\pi^*}$ of $B$ at $\pi^*$ is
homological if and only if the universal localization
$\lambda'_{\Gamma}:S'\to S'_{\Gamma}$ of $S'$ at $\Gamma$ is
homological. In particular, if one of the above equivalent
conditions holds , then there is a recollement of derived module
categories:

$$\xymatrix@C=1.2cm{\D{S'_{\Gamma}}\ar[r]&\D{B}\ar[r]
\ar@/^1.2pc/[l]\ar@/_1.2pc/[l]
&\D{R}\ar@/^1.2pc/[l]\ar@/_1.2pc/[l]}.$$
\label{new1}
\end{Koro}

As a consequence of Corollary \ref{new1}, we obtain the following
result which can be used to show when the universal localization
$\lambda_{\pi^*}:B\to B_{\pi^*}$ of $B$ at $\pi^*$ in Proposition
\ref{new} is homological in some special cases.

\begin{Koro}
Let $C\subseteq D$ be an extension of rings. Set $R:=\left(\begin{array}{cc} D & D\\
0 & C\end{array}\right)$ and $S:=M_2(D)$. Let $\lambda:R\to S$ be
the canonical injective ring homomorphism. Then the universal
localization $\lambda_{\pi^*}:B\to B_{\pi^*}$ of $B$ at $\pi^*$ is
homological if and only if the universal localization
$\lambda_{\omega^*}:E\to E_{\omega^*}$ of $E$ at ${\omega^*}$ is
homological, where $E:=\End_C(D\oplus D/C)$, and the homomorphism
$\omega^*: \Hom_C(D\oplus D/C,D)\to\Hom_C(D\oplus D/C, D/C)$  of
projective $E$-modules  is induced by the canonical epimorphism
$\omega: D\ra D/C$. \label{new2}
\end{Koro}
{\it Proof.} Recall that the right multiplication map $\mu: R\ra S'$
is a ring homomorphism (see Lemma\ref{epi2}(1)).
Set $e_1:= \left(\begin{array}{ll} 1 &0\\
0 & 0\end{array}\right)$, $e_2:=\left(\begin{array}{ll} 0 &0\\
0 & 1\end{array}\right)$ and $e_{12}:=\left(\begin{array}{ll} 0 &1\\
0 & 0\end{array}\right)\in R$. Furthermore, let $\varphi: Re_1\ra
Re_2$ and $\varphi':S'(e_1)\mu\to S'(e_2)\mu$ be the right
multiplication maps by $e_{12}$ and $(e_{12})\mu$, respectively.

It follows from Lemma \ref{epi0} and $D\sqcup_{C}C=D$ that the map
$\lambda:R\to S$ is  the universal localization of $R$ at $\varphi$.
In particular, $\lambda$ is a ring epimorphism. Since $S\simeq
e_1R\oplus e_1R$ as right $R$-modules, the embedding $\lambda$ is
always homological. Note that $S'\otimes_R\varphi$ can be identified
with $\varphi'$. By Lemma \ref{new1}, the universal localization
$\lambda_{\pi^*}:B\to B_{\pi^*}$ of $B$ at $\pi^*$  is homological
if and only if the universal localization $\lambda'_{\varphi'}:S'\to
S'_{\varphi'}$ of $S'$ at $\varphi'$ is homological.

Clearly, $R/Re_1R\simeq C$ as rings. So, every $C$-module can be
regarded as an $R$-module. In particular,  $D\oplus D/C$ can be seen
as an $R$-module. Further, one can check that the map
$\alpha:D\oplus D/C\to S/R $ defined by
$$(d,t+C)\mapsto \left(\begin{array}{lc} 0 & 0\\
d & t\end{array}\right)+R$$ for $d,t\in D$, is an isomorphism of
$R$-modules. Thus $S'\simeq E$, $\varphi'$ corresponds to $\omega^*$
under this isomorphism, and therefore $S'_{\varphi'}\simeq
E_{\omega^*}$. It follows that $\lambda'_{\varphi'}:S'\to
S'_{\varphi'}$ is homological if and only if so is
$\lambda_{\omega^*}:E\to E_{\omega^*}$. This finishes the proof.
$\square$

\medskip
{\it Remark.} In general, the special form of the universal
localization $\lambda_{\omega^*}:E\to E_{\omega^*}$ of $E$ at
${\omega^*}$(or equivalently, the universal localization
$\lambda_{\pi^*}:B\to B_{\pi^*}$ of $B$ at $\pi^*$) in Corollary
\ref{new2} may not be homological, through the $\lambda$ is always
homological.

For example, let $C=\big\{\left(\begin{array}{lc} a & 0\\
b & a\end{array}\right)\mid a,b\in k\big\}$ and $D=\left(\begin{array}{lc} k & 0\\
k & k\end{array}\right)$ with $k$ a field.  Then one can check that
the canonical map $\omega: D\ra D/C$ is a split epimorphism in
$C$-Mod, and therefore $_CD\simeq C\oplus D/C$. Let $e$ be the
idempotent of $E$ corresponding the direct summand $C$ of the
$C$-module $D\oplus D/C$. Then $E_{\omega^*}\simeq E/EeE\simeq
M_2(k)$. Furthermore, the universal localization
$\lambda_{\omega^*}:E\to E_{\omega^*}$ of $E$ at ${\omega^*}$ is
equivalent to  the canonical projection $\tau:E\to E/EeE$. Since
$\Ext^2_E(E/EeE,E/EeE)\neq 0$, we see that $\tau$ is not
homological. This implies that $\lambda_{\omega^*}$ is not
homological, too. Thus $\lambda_{\pi^*}:B\to B_{\pi^*}$ is not
homological by Corollary \ref{new2}, that is, the restriction
functor $(\lambda_{\pi^*})_*:\D{B_{\pi^*}}\to \D{B}$ is not fully
faithful.

\section {Commutative rings and recollements of derived module
categories \label{sect7}}

In this section, we shall first discuss recollements of derived
module categories arising from injective homological ring
epimorphisms $\lambda: R\ra S$ between arbitrary commutative rings
without the assumption that the modules $S\oplus S/R$ are tilting
modules, and then turn to the special case of one-Gorenstein rings.
We shall see that, for commutative rings, the universal
localizations appearing in the main result Theorem \ref{th01} will
be further strengthened as tensor products. As a consequence, we can
produce examples to show that two different stratifications of a
derived module category by derived module categories of rings may
have different derived composition factors, which answers negatively
a question in \cite{akl} and shows that the Jordan-H\"older theorem
fails for derived module categories with simple derived module
categories as composition factors.

Note that if $R$ is a commutative ring and $\lambda: R\ra S$ is a
ring epimorphism, then $S$ must be commutative. So, in the
following, we can assume that both rings $R$ and $S$ are commutative
rings.

\subsection{General case: Arbitrary commutative rings}

The main purpose of this subsection is to prove the following
existence result for recollements arising form injective ring
epimorphisms between commutative rings. Here we remove the condition
of tilting modules.

\begin{Prop} Suppose that $\lambda: R\ra S$ is an injective
homological ring epimorphism between commutative rings $R$ and $S$.
Let $B$ be the endomorphism ring of the $R$-module $S\oplus\;
S\!/\!R$. Then there is a recollement of derived module categories:

$$\xymatrix@C=1.2cm{\D{S\otimes_RS'}\ar[r]&\D{B}\ar[r]
\ar@/^1.2pc/[l]\ar@/_1.2pc/[l]
&\D{R}\ar@/^1.2pc/[l]\ar@/_1.2pc/[l]},\vspace{0.3cm}$$ where
$S':=\End_R(S/R)$ is commutative, and $S\otimes_RS'$ is the tensor
product of $S$ and $S'$ over $R$.

In particular, the left global dimension of $B$ is finite if and
only the left global dimensions of $R$ and $S\otimes_RS'$ both are
finite. \label{com1}
\end{Prop}

{\it Proof.} Since $S$ is a commutative ring, we see that $S'$ is an
$R$-algebra via the right multiplication map $\mu:R\to S'$, defined
by $r\mapsto(x\mapsto xr)$ for $r\in R$ and $x\in S/R$ (see Lemma
\ref{epi2}(1)). Then the tensor product $S\otimes_RS'$ of $S$ and
$S'$ over $R$ makes sense. Moreover, the map $\lambda': S'\to
S\otimes_RS'$, defined by $s'\mapsto
 1\otimes s'$ for any $s'\in S'$,  and the map $\mu': S\ra
S\otimes_RS'$,  defined by $s\mapsto s\otimes 1$ for $s\in S$, are
ring homomorphisms. So, $S\otimes_RS'$ is an $S'$-$S'$-bimodule via
$\lambda'$. In addition, it follows from Lemma \ref{epi2}(5) that
$S'$ is a commutative ring.

Since $\lambda$ is a ring epimorphism, we know from Lemma
\ref{epi00} that the tensor product $S\otimes_RS'$ together with the
two ring homomorphisms $\lambda'$ and  $\mu'$ satisfies the
definition of coproducts. In other words, the coproduct
$S\sqcup_{R}S'$  of $R$-rings $S$ and $S'$ over $R$ is isomorphic to
the tensor product $S\otimes_{R}S'$.

By Proposition \ref{new}, to get the recollement of derived module
categories in Proposition \ref{com1}, it is sufficient to
demonstrate that $\lambda': S'\to S\otimes_RS'$ is homological.

In fact, we have  the following composition of a series of
isomorphisms of $S$-$S'$bimodules:
$$ (S\otimes_RS')\otimes_{S'}(S\otimes_RS')\simeq
S\otimes_R\big(S'\otimes_{S'}(S\otimes_RS')\big)\simeq
S\otimes_R(S\otimes_RS')\simeq (S\otimes_RS)\otimes_RS'\simeq
S\otimes_RS'.$$ This composition map is actually the multiplication
map from $(S\otimes_RS')\otimes_{S'}(S\otimes_RS')$ to
$S\otimes_RS'$. Thus it is an
$(S\otimes_RS')$-$(S\otimes_RS')$-bimodule isomorphism. Hence
$\lambda'$ is a ring epimorphism.

It remains to show $\Tor^{S'}_i(S\otimes_RS',S\otimes_RS')=0$ for
all $i>0$. However, this follows immediately from Step $(1)$ in the
proof of Proposition \ref{new}. Thus, we finish the proof of
Proposition \ref{com1}.

The last statement follows from a result in \cite{Happel}. $\square$

\subsection{Special case: One-Gorenstein rings\label{subsect6.2}}

Throughout this subsection, $R$ will stand for a commutative ring.
We denote by Spec$(R)$ (respectively, mSpec$(R)$ ) the set of all
prime (respectively, maximal) ideals of $R$. For each non-negative
integer $i$, we denote by $P_i$ the set of all prime ideals of $R$
with height $i$.

Let $M$ be an $R$-module. We denote by $E(M)$ the injective envelope
of $M$, and by $\pd(M),\id(M)$ and $\fld(M)$ the projective,
injective and flat dimensions of $_RM$, respectively.

For a multiplication subset $\Sigma$ of $R$, we denote by
$\Sigma^{-1}R$ the localization of $R$ at $\Sigma$, and by
$f_{\Sigma}:R\to \Sigma^{-1}R$ the canonical homological ring
epimorphism. In general, the homomorphism $f_{\Sigma}$ is not
injective. But, if $\Phi$ is the multiplicative set of all non-zero
divisors of $R$, then the localization map $f_{\Phi}: R\ra
{\Phi}^{-1}R$ is always injective. In this case, the ring
${\Phi}^{-1}R$ is called the total quotient ring of $R$, denoted by
$Q$. In fact, $Q$ is the largest localization of $R$ for which the
canonical map is injective, that is, if the map $f_{\Sigma}:R\to
\Sigma^{-1}R$ is injective, then $\Sigma\subseteq \Phi$, and there
is a unique injective ring homomorphism $h:\Sigma^{-1}R\to Q$ such
that $f_{\Phi}=f_\Sigma\,h$. In addition, if $R$ is noetherian, then
$P_0$ is finite and  $\Phi =R\setminus \cup_{\mathfrak p\in
P_0}\mathfrak p$.

As usual, for a prime ideal $\mathfrak p$ of $R$, we always write
$R_{\mathfrak p}$ for $(R\setminus\mathfrak p)^{-1}R$, and
$f_{\mathfrak p}$ for $f_{R\setminus {\mathfrak p}}$, and say that
$R_{\mathfrak p}$ is the localization of $R$ at ${\mathfrak p}$.

Note that the localization ${\mathbb Z}_{\mathfrak p}$ of $\mathbb
Z$ at the maximal ideal ${\mathfrak p}=p\mathbb Z$ is ${\mathbb
Q}_{(p)}$ for every prime $p\in {\mathbb N}$, where ${\mathbb
Q}_{(p)}$ is the set of $p$-integers. Recall that $q=n/m \in
{\mathbb Q}$ with $m,n\in \mathbb Z$ is called a $p$-integer if $p$
does not divide $m$.

Let $\mathfrak p_1, \cdots, \mathfrak p_n$ be prime ideals in $R$,
and set $\Sigma :=R\setminus\cup_{i=1}^n\mathfrak p_i$. Then
$\Sigma=\cap_{i=1}^n(R\setminus \mathfrak p_i)$ is a multiplicative
subset of $R$, and the prime ideals of the localization
$\Sigma^{-1}R$ are in one-to one correspondence with the prime
ideals $\mathfrak p$ of $R$ with $\mathfrak p\cap \Sigma
=\emptyset$, that is, with the prime ideals of $R$ contained in
$\cup_{i=1}^nP_i$. By prime avoidance theorem, any such prime ideal
is contained in one of the $\mathfrak p_i$. Hence,
$\{\Sigma^{-1}\mathfrak p_j\mid 1\le j\le n\}$ contains the maximal
ideals of $\Sigma^{-1}R$. If all $\mathfrak p_i$ are pairwise
incomparable, that is $\mathfrak p_i\not\subseteq \mathfrak p_j$ for
$i\ne j$, then this is exactly the set of all maximal ideals of
$\Sigma^{-1}R$.

Now, let us mention the following homological results about
commutative noetherian rings, which are needed for our discussions
in this section. For more details, we refer the reader, for
instance, to \cite[Theorem 3.3.8, Theorem 3.4.1, Lemma 6.7.7]{EJ}
and \cite[Corollary 11.2]{ne}.

\begin{Lem}
Suppose that $R$ is a noetherian ring. Let $\mathfrak p, \mathfrak
q\in Spec(R)$. We define $J_{\mathfrak p}:=\varprojlim_iR/{\mathfrak
p^i}$. Then the following hold:

$(1)$ $\Hom_R(E(R/\mathfrak p), E(R/\mathfrak q))\neq 0$ if and only
if $\mathfrak p \subseteq \mathfrak q$. In particular,
$E(R/\mathfrak p)\simeq E(R/\mathfrak q)$ if and only if $\mathfrak
p=\mathfrak q$.

$(2)$ If $\Sigma$ is a multiplication subset of $R$, then, as
$R$-modules,

$$\Sigma^{-1}E(R/\mathfrak p)\simeq
\left\{\begin{array}{ll} E(R/\mathfrak p) & \mbox{if } \Sigma\cap\mathfrak p=\emptyset,\\
0 & \mbox{if } \Sigma\cap\mathfrak p\neq\emptyset.\end{array}
\right. $$

$(3)$ If $\mathfrak p$ is a maximal ideal of $R$, then
$$\End_R(E(R/\mathfrak p))\simeq\varprojlim_iR_{\mathfrak
p}/{\mathfrak p^iR_{\mathfrak p}}\simeq J_{\mathfrak p}.$$

$(4)$ Let $P$ be a set of maximal ideals of $R$. If $\mathfrak q$ is
a maximal ideal of $R$, which does not belong to $P$, then
$$
E(R/\mathfrak q)\otimes_R\prod_{\mathfrak p\in P}J_{\mathfrak
p}=0.$$

$(5)$ Define $S_n:=\{x\in E(R/\mathfrak p)\,|\, \mathfrak p^n x=0\}$
for each  $n>0$. Then $E(R/\mathfrak p)=\cup_nS_n$. \label{lem6.9}
\end{Lem}

Of our interest is the class of $n$-Gorenstein rings. Recall that,
for a non-negative integer $n$, the ring $R$ is called
$n$-Gorenstein if $R$ is noetherian and $\id(R)\leq n$. The
following homological properties of $n$-Gorenstein rings are
well-known, for their proofs, we refer to \cite[Theorem 1, Theorem
6.2]{Bs}, \cite[Theorem 9.1.10, Theorem 9.1.11]{EJ} and
\cite[Introduction]{TP}.

\begin{Lem} Let $n$ be a non-negative integer. Assume that $R$ is an $n$-Gorenstein ring.

$(1)$ The regular module $R$ has a minimal injective resolution of
the form:
$$
0\lra R\lra\bigoplus_{\mathfrak p\in P_0}E(R/\mathfrak
p)\lra\cdots\lra \bigoplus_{\mathfrak p\in P_n}E(R/\mathfrak p)\lra
0 .
$$ Moreover, the total quotient ring $Q$ of $R$ is isomorphic to
$E(R) $ as an $R$-module.

 $(2)$ Let $M$ be an $R$-module. Then the following are equivalent:
$$\begin{array}{rlcrlcrl}
 (i)  & \pd(M)<\infty ; &\qquad & (ii) & \;\id(M)<\infty ;& \qquad & (iii) & \;\fld(M)<\infty ; \\
 (iv) & \;\pd(M)\leq n; &\qquad & (v)  & \;\id(M)\leq n;  & \qquad & (vi)  & \; \fld(M)\leq n.
\end{array}$$

$(3)$ The $R$-module
$$T_{(P_n)}:=\bigoplus_{i\leq
n}\bigoplus_{\mathfrak p\in P_i}E(R/\mathfrak p)$$ is an (infinitely
generated) $n$-tilting module, that is, it is of projective
dimension at most $n$, and satisfies (T2) and (T3) (replaced by a
longer exact sequence) in Definition \ref{def1}.

$(4)$ If $\Sigma$ is a multiplication subset of $R$, then
$\Sigma^{-1}R$ is an $n$-Gorenstein  ring. \label{lem6.10}
\end{Lem}

From now on, we assume that \textbf{$R$ is a $1$-Gorenstein ring}.
Then $P_1$ consists of all maximal ideals of $R$ which are not
minimal prime ideals. By Lemma \ref{lem6.10}, one gets a tilting
$R$-module $T_{(P_1)}$ of projective dimension at most one. This
construction of tilting modules from $1$-Gorenstein rings can be
generalized to obtain the so-called Bass tilting modules, as
mentioned in \cite {HHT}. Now, let us recall the construction.

Let
$$0\lra R\lraf {f_\Phi}Q\lraf{\pi}\bigoplus_{\mathfrak p\in P_1}E(R/\mathfrak p)\lra 0,$$
be a minimal injective resolution of $R$, where $\pi$ is the
canonical surjective map which is regarded as a homomorphism of
$R$-modules. Let $\Delta$ be a subset of $P_1$. Then we define
$$
R_{(\Delta)}:=\pi^{-1}(\bigoplus_{\mathfrak p\in
\Delta}E(R/\mathfrak p))\quad  \mbox{and}\quad
T_{(\Delta)}:=R_{(\Delta)}\oplus\bigoplus_{\mathfrak p\in
\Delta}E(R/\mathfrak p).
$$
Clearly, we get two associated exact sequences of $R$-modules
$$
(a)\quad \quad 0\lra R\lra
R_{(\Delta)}\lraf{\pi}\bigoplus_{\mathfrak p\in \Delta}E(R/\mathfrak
p)\lra 0\,;
$$
$$
(b)\quad \quad 0\lra R_{(\Delta)}\lra Q\lra\bigoplus_{\mathfrak p\in
{P_1\setminus\Delta}}E(R/\mathfrak p)\lra 0.
$$
Note that $R_{(\Delta)}$  is just an $R$-submodule of $Q$. It is
shown in \cite [Section 4]{HHT} that the $R$-module $T_{(\Delta)}$
is a tilting module, which is called a Bass tilting module over $R$.
Further, the authors of \cite{TP} prove that  every tilting module
over $R$ is equivalent to a Bass tilting module. Note that the
sequence $(a)$ implies that the $R$-tilting module $T_{(\Delta)}$ is
good.

The next lemma describes some properties relevant to Bass tilting
modules. Note that the conclusions (1) and (2) of Lemma
\ref{lem6.11} below are mentioned in \cite{tw2} for Dedekind
domains.

\begin{Lem}  let $\Delta$ be a subset of $P_1$. Assume that each prime
ideal belonging to the set $P_1\setminus\Delta$ contains all zero
divisors of $R$. Then we have the following:

$(1)$ For each $\mathfrak p\in {P_1\setminus\Delta}$, the canonical
ring homomorphism $f_{\mathfrak p}:R\to R_{\mathfrak p}$ is
injective.

$(2)$ $R_{(\Delta)}=\bigcap_{{\mathfrak p}\in {P_1\setminus\Delta}}
R_{\mathfrak p}$, which is a flat $R$-module, where $R_{\mathfrak
p}$ is regarded as a subring of the total ring $Q$ of $R$. Hence
$R_{(\Delta)}$ can be regarded as a subring of $Q$ containing $R$.
In particular, the total quotient ring of $R_{(\Delta)}$ also equals
$Q$. (Note that we set $\bigcap_{{\mathfrak p}\in \emptyset
}R_{\mathfrak p}=Q$.)

$(3)$ The canonical inclusions  $\lambda_{\Delta}:R\to R_{(\Delta)}$
and $\mu_{\Delta}:R_{(\Delta)}\to Q$ are  homological ring
epimorphisms.

$(4)$ The $R_{(\Delta)}$-module
$$T'_{(\Delta)}:=Q\oplus\bigoplus_{\mathfrak p\in
P_1\setminus\Delta}E(R/\mathfrak p)$$ is a good tilting
$R_{(\Delta)}$-module.  If  $R_{(\Delta)}$ is noetherian, then
$R_{(\Delta)}$ is $1$-Gorenstein.

$(5)$
$$B_\Delta:=\End_R(T_{(\Delta)})\simeq
\left(\begin{array}{lc}
R_{(\Delta)}& R_{(\Delta)}\otimes_RJ_{\Delta}\\
0 & J_{\Delta}
\end{array}\right),
$$
where $J_{\Delta}:=\End_R(R_{\Delta)}/R)\simeq \prod_{\mathfrak p\in
\Delta}J_{\mathfrak p}$.

$(6)$
$$B'_\Delta:=\End_{R_{(\Delta)}}(T'_{(\Delta)})\simeq
\left(\begin{array}{lc}
Q& Q\otimes_RJ'_{\Delta}\\
0 & J'_{\Delta}
\end{array}\right),
$$
where
$J'_{\Delta}:=\End_{R_{(\Delta)}}(Q/R_{(\Delta)})\simeq\prod_{\mathfrak
p\in {P_1\setminus\Delta}}J_{\mathfrak p}.$

$(7)$ The ring homomorphism $\mu_\Delta$ induces a ring isomorphism
$$R_{(\Delta)}\otimes_RJ_{\Delta}\simeq Q\otimes_RJ_{\Delta}.$$

$(8)$ For any subset $P$ of $P_1$, the canonical map
$$\Theta_P:\, Q\otimes_R\prod_{\mathfrak p\in P}J_{\mathfrak
p}\lra\prod_{\mathfrak p\in P}Q\otimes_RJ_{\mathfrak p},$$ defined
by $q\otimes\big(x_{\mathfrak p}\big)_{\mathfrak p\in P}\mapsto
\big(q\otimes x_{\mathfrak p}\big)_{\mathfrak p\in P}$ for $q\in Q$
and $x_\mathfrak p\in J_{\mathfrak p}$, is an injective ring
homomorphism. \label{lem6.11}
\end{Lem}

{\it Proof.} $(1)$ Note that, for each $r\in \Ker(f_{\mathfrak p})$,
there exists an element $x\in R\setminus\mathfrak p $ such that
$rx=0$. Since $\mathfrak p$ contains all zero divisors of $R$, we
know that $x$ is non-zero divisor of $R$. This implies $r=0$, and so
the map $f_{\mathfrak p}$ is injective.

$(2)$ Let ${\mathfrak q}\in P_1\setminus \Delta$. Since the
localization map $f_{\mathfrak q}: R\ra R_{\mathfrak q}$ is
injective by (1), there is a unique injective homomorphism
$\mu_{\mathfrak q}: R_{\mathfrak q}\ra Q$ such that
$f_{\Phi}=f_{\mathfrak q}\mu$ by the universal property of the total
quotient ring of $R$. So, we can think of $R_{\mathfrak q}$ as a
subring of $Q$. Under this identification, we can speak of the
intersection of $R_{\mathfrak p}$ defined in (2).

First, we  show that  if $\Delta=P_1\setminus\{\mathfrak p\}$ for
some $\mathfrak p\in P_1$, then $R_{(\Delta)}=R_{\mathfrak p}$. By
Lemma \ref{lem6.10}(1), we have the following  exact commutative
diagram:
\begin{eqnarray*}
\xymatrix{
&       & 0\ar[d]             &0\ar[d]         &\\
0\ar[r]&R\ar[r]^-{f_{\mathfrak p}}\ar@{=}[d]&R_\mathfrak
p\ar[d]^-{\mu_{\mathfrak p}}
\ar[r]^-{\mu_{\mathfrak p}\pi}\ar[d]&Y\ar[r]\ar[d]&0\\
0\ar[r]& R\ar[r]^-{f_\Phi} &Q\ar[d]\ar[r]^-{\pi}&
\bigoplus_{\mathfrak q\in P_1}E(R/\mathfrak q)\ar[r]\ar[d]\ar[d]^-{(g_{\mathfrak q})_{{\mathfrak q}\in P_1}}&0.\\
&        &                         E(R/\mathfrak p)\ar@{=}[r]\ar[d]
&
E(R/\mathfrak p)\ar[d]& \\
&        &                     0 &  0& \\
}
\end{eqnarray*}
where $Y$ denotes the image of the restriction of $\pi$ to
$R_\mathfrak p$, and where $g_{\mathfrak q}: E(R/\mathfrak q)\ra
E(R/\mathfrak p)$ is a homomorphism of $R$-modules. Clearly, the
localization $Y_\mathfrak p$ of $Y$ at $\mathfrak p$ is zero.  Let
$\mathfrak a\in\Delta$. By Lemma \ref{lem6.9}(1), we know that
$\Hom_R(E(R/\mathfrak a), E(R/\mathfrak p))=0$ since both $\mathfrak
a$ and $\mathfrak p$ are maximal ideals of $R$. Consequently,
$g_{\mathfrak a}=0$ and $Y= \Ker(g_{\mathfrak
p})\oplus\bigoplus_{\mathfrak q\in \Delta}E(R/\mathfrak q)$, where
$g=g_{\mathfrak p}: E(R/\mathfrak p)\to E(R/\mathfrak p)$ is a
surjective homomorphism of $R$-modules. We claim $\Ker(g)=0$. In
fact, by Lemma \ref{lem6.9}(2), we know  that $E(R/\mathfrak
p))\simeq (E(R/\mathfrak p))_{\mathfrak p}$ as $R$-modules. This
implies that $\Ker(g)\simeq \Ker(g)_{\mathfrak p}$ as $R$-modules.
Then it follows from $Y_\mathfrak p=0$ that $\Ker(g)_\mathfrak p=0$.
Thus $R_{\mathfrak p}=R_{(\Delta)}$ under our identification of
$R_{\mathfrak p}$ in $Q$.

Second, in the general case, we observe that
$$\Delta=\bigcap_{\mathfrak b\in
P_1\setminus\Delta}(P_1\setminus\{\mathfrak b\})\quad
\mbox{and}\quad \bigoplus_{\mathfrak p\in \Delta}E(R/\mathfrak
p)=\bigcap_{\mathfrak b\in P_1\setminus\Delta}(\bigoplus_{\mathfrak
q\in P_1\setminus\{\mathfrak b\}}E(R/\mathfrak q)\,).$$ Thus
$$ R_{(\Delta)}= \pi^{-1}\big(\bigoplus_{{\mathfrak p}\in \Delta}E(R/{\mathfrak p})\big)=\pi^{-1}\bigg(\bigcap_{\mathfrak b\in P_1\setminus\Delta}\;
\big(\bigoplus_{q\in P_1\setminus\{\mathfrak b\}}E(R/{\mathfrak q}
)\big)\bigg)=\bigcap_{\mathfrak p\in P_1\setminus\Delta}R_{\mathfrak
p}.$$

Third, to prove that $R_{(\Delta)}$ is a flat $R$-module, we use the
exact sequence $(b)$. Since $R$ is noetherian, the arbitrary direct
sum of injective $R$-modules  is injective. Thus
$\bigoplus_{\mathfrak p\in {P_1\setminus\Delta}}E(R/\mathfrak p)$ is
injective. Note that $_RQ$ is flat. By Lemma \ref{lem6.10}(2), we
deduce that  $R_{(\Delta)}$ is a flat $R$-module.

Finally, note that the total quotient ring of $R_{(\Delta)}$ also
equals $Q$.

$(3)$ Recall that $Q$ is the localization of $R$ at the
multiplication set $\Phi$ consisting of all non-zero divisors of
$R$. Clearly, the map $f_\Phi :R\to Q$ is a ring epimorphism. Note
that $Q$ and $R_{(\Delta)}$ are flat $R$-modules. It follows from
Lemma \ref{prop6.7} that both $\lambda_{\Delta}$ and $\mu_{\Delta}$
are homological ring epimorphisms.

$(4)$ For simplicity, we set $W:=T'_{(\Delta)}$. Note that the
$R_{(\Delta)}$-module $W$ is injective as an $R$-modules. Since
$\lambda_{\Delta}$ is a homological ring epimorphism, it follows
from Lemma \ref{2.1} that $W$ is an injective $R_{(\Delta)}$-module.
In particular, the regular module $R_{(\Delta)}$ has injective
dimension at most $1$. If $R_{(\Delta)}$ is noetherian, then
$R_{(\Delta)}$ is $1$-Gorenstein by definition.  In this case, it
directly follows from \cite[Section 4]{HHT} that the module $W$ is a
tilting $R_{(\Delta)}$-module. However, in general, we do not know
whether $R_{(\Delta)}$ is noetherian or not. Because of this reason,
we have to prove, in the following, that $W$ is a tilting
$R_{(\Delta)}$-module. Indeed, since $W$ is an injective $R$-module
and $R$ is $1$-Gorenstein, we know from Lemma \ref{lem6.10}(2) that
$\pd(_RW)\leq 1$. Note that $\lambda_{\Delta}$ is a homological ring
epimorphism. It follows from Lemma \ref{2.1} that
$\pd(_{R_{(\Delta)}}W)\leq\pd(_RW)\leq 1$. Clearly, the module $W$
satisfies the condition $(T_3)$ in Definition \ref{def1}. To see
that $W$ satisfies the condition $(T_2)$ in Definition \ref{def1},
we observe that
$$\Ext^i_{R_{(\Delta)}}(W,W^{({\alpha})})\simeq\Ext^i_R(W,W^{({\alpha})})=0$$
for every  $i\geq 1$ and every cardinal $\alpha$, where the first
isomorphism   follows from Lemma \ref{2.1}(3), and  the second
equality follows from the fact that every direct sum of injective
$R$-modules is injective since $R$ is noetherian.  Thus $W$ is a
tilting $R_{(\Delta)}$-module. Clearly, the exact sequence $(b)$
implies that  $W$ is a good tilting.

$(5)$ By Lemma \ref{lem6.9}(1) and Lemma \ref{lem6.9}(3), we have
$J_{\Delta}\simeq \End_R(\bigoplus_{\mathfrak p\in
\Delta}E(R/\mathfrak p))\simeq \prod_{\mathfrak p\in
\Delta}J_{\mathfrak p}$. According to Lemma \ref{epi2}, we have
$\Hom_R\big(R_{(\Delta)}, \bigoplus_{\mathfrak p\in
\Delta}E(R/\mathfrak p\big) $ $\simeq
R_{(\Delta)}\otimes_RJ_{\Delta}$ as
$R_{(\Delta)}$-$J_{\Delta}$-bimodules. Now, (5) follows from Lemma
\ref{epi1}(2) immediately.

$(6)$ We first observe that
$$
\Hom_R(X, Y)\simeq\Hom_{R_{(\Delta)}}(X, Y)\quad \mbox{and}\quad
X\otimes_RY\simeq X\otimes_{R_{(\Delta)}}Y $$ for any
$R_{(\Delta)}$-modules $X$ and $Y$ since  $\lambda_{\Delta}:R\to
R_{(\Delta)}$ is a ring epimorphism, and then use Lemma
\ref{epi1}(2), Lemma \ref{epi2} and Lemma \ref{lem6.9}. We omit the
details here.

$(7)$  Note that if $R$ is a commutative noetherian ring and if $I$
is an ideal of $R$, then (i) the $I$-adic completion of $R$ is a
flat $R$-module, and (ii) the product of flat $R$-modules is flat
(see \cite[Corollary 2.5.15, Theorem 3.2.24]{EJ}). Hnece
$J_{\Delta}$ is a flat $R$-module. In order to prove that
$\mu_\Delta\otimes_RJ_\Delta:R_{(\Delta)}\otimes_RJ_{\Delta}\to
Q\otimes_RJ_{\Delta}$ is an isomorphism, it is sufficient to show
$(\bigoplus_{\mathfrak p\in {P_1\setminus\Delta}}E(R/\mathfrak
p))\otimes_RJ_{\Delta}=0$. This is equivalent to $E(R/\mathfrak
p)\otimes_RJ_{\Delta}=0$ for any ${\mathfrak p\in
{P_1\setminus\Delta}}$. However, the latter is a direct consequence
of Lemma \ref{lem6.9}(4).

$(8)$ Clearly, the map $\Theta_P$ is a ring homomorphism. Applying
the tensor functors $-\otimes_R\prod_{\mathfrak p\in P}J_{\mathfrak
p}$ and $-\otimes_RJ_{\mathfrak p}$ to the minimal injective
coresolution of $R$, respectively, we can get the following
 exact commutative diagram of $R$-modules:
$$\xymatrix{0\ar[r]&R\otimes_R\prod_{\mathfrak p\in
P}J_{\mathfrak p}\ar^-{\wr}[d]\ar[r]&Q\otimes_R\prod_{\mathfrak p\in
P}J_{\mathfrak p}\ar[d]^-{\Theta_P}\ar[r]&\big(\bigoplus_{\mathfrak
q\in P_1}E(R/\mathfrak q)\big)\otimes_R\prod_{\mathfrak p\in
P}J_{\mathfrak p}\ar[d]^-{\Theta'_P}\ar[r]&0\\
0\ar[r]& \prod_{\mathfrak p\in P}R\otimes_RJ_{\mathfrak p}\ar[r]
&\prod_{\mathfrak p\in P}Q\otimes_RJ_{\mathfrak p} \ar[r]
&\prod_{\mathfrak p\in P}\big(\bigoplus_{\mathfrak q\in
P_1}E(R/\mathfrak q)\big)\otimes_RJ_{\mathfrak p}\ar[r]&0, } $$
where we define the homomorphism $\Theta'_P$ of $R$-modules in the
same way as we did for $\Theta_P$. We claim that $\Theta'_P$ is
injective. In fact, since tensor functor commutes with direct sums,
it follows from Lemma \ref{lem6.9}(4) that we can embed $\Theta'_P$
into the following commutative diagram:
$$\xymatrix{
\big(\bigoplus_{\mathfrak q\in P_1}E(R/\mathfrak
q)\big)\otimes_R\prod_{\mathfrak p\in P}J_{\mathfrak
p}\ar[r]^-{\sim}\ar[d]^-{\Theta'_P}&\bigoplus_{\mathfrak p\in
P}E(R/\mathfrak p)\otimes_RJ_{\mathfrak p}\ar[d]^-{\lambda}\\
\prod_{\mathfrak p\in P}\big(\bigoplus_{\mathfrak q\in
P_1}E(R/\mathfrak q)\big)\otimes_RJ_{\mathfrak p}\ar[r]^-{\sim}
&\prod_{\mathfrak p\in P}E(R/\mathfrak p)\otimes_RJ_{\mathfrak p}\,,
}
$$
where the map $\lambda$ is the canonical inclusion. This shows that
$\Theta'_P$ is injective, which implies that $\Theta_P$ also is
injective. $\square$

\medskip
Note that if $R$ is  a local ring or a domain, then the assumption
in Lemma \ref{lem6.11} holds. It is well-known that Dedekind domains
are $1$-Gorenstein rings. Recall that a commutative ring $R$ is
called a Dedekind domain if $R$ is a domain in which every ideal
($\neq R$) is the product of a finite numbers of prime ideals. This
is equivalent to saying that $R_{\mathfrak p}$ is a discrete
valuation ring for each prime ideal $\mathfrak p$ of $R$. A typical
example of Dedekind domain is the ring $\mathbb{Z}$ of rational
integers.

The assumption of Lemma \ref{lem6.11} cannot be droped for (1) to
hold true. For example, if $R$ is a local $1$-Gorenstein ring, then
the direct sum $S:= R\oplus R$ of two copies of $R$ is again
$1$-Gorenstein. If we take $\mathfrak m$ to be the unique maximal
ideal of $R$, then the localization of $S$ at the maximal ideal
${\mathfrak p}:=({\mathfrak m},R)$ is isomorphic to $R_{\mathfrak
m}$. This shows that the localization map $S\ra S_{\mathfrak p}$ is
not injective.

 By Lemma \ref{lem6.11}(2), we know that $R_{(\Delta)}$ is always
an intersection of localizations. But, in general, it may not be a
localization of $R$ at any multiplication set. For a counterexample,
we refer the reader to \cite{NS1}. A natural question arises: when
is $R_{(\Delta)}$ itself a localization of $R$ at some
multiplication set ? The following result provides some partial
answers to this question.

\begin{Lem}
 Let $\Delta$ be a subset of $P_1$. Assume that each prime
ideal belonging to $P_1\setminus\Delta$ contains all zero divisors
of $R$. Define $\Sigma:=R\setminus\bigcup_{\mathfrak q\in
{P_1\setminus\Delta}}\mathfrak q$ and $\Delta_1:=\{\mathfrak a\in
P_1\mid \mathfrak a\subseteq\bigcup_{\mathfrak q\in
{P_1\setminus\Delta}}\mathfrak q\}$.  Then we have the following:

$(1)$ $\Sigma^{-1}R=\pi^{-1}\big(\bigoplus_{\mathfrak
p\in{P_1\setminus \Delta_1}}E(R/\mathfrak p)\big)\subseteq
R_{(\Delta)}\subseteq Q$.

$(2)$ $R_{(\Delta)}=\Sigma_1^{-1}R$ for some multiplication subset
$\Sigma_1$ of $R$ if and only if  $R_{(\Delta)}=\Sigma^{-1}R $ if
and only if $\Delta_1=P_1\setminus\Delta$.

$(3)$ If $P_1\setminus\Delta$ is a finite set, or if each  ideal in
$\Delta$ is a principal ideal of $R$, then $
R_{(\Delta)}=\Sigma^{-1}R$. \label{lem6.11'}
\end{Lem}
{\it Proof.} $(1)$ Clearly, we have $\Sigma\subseteq \Phi$ and
$\Sigma$ is a multiplicative set. Thus the canonical  map
$f_{\Sigma}:R\to \Sigma^{-1}R$ is injective, and there is a unique
injective ring homomorphism $h:\Sigma^{-1}R\to Q$ such that
$f_{\Phi}=f_\Sigma\,h$. In this sense, we may regard $\Sigma^{-1}R$
as a subring of the total quotient ring $Q$ containing $R$.
Moreover, the total quotient ring of $\Sigma^{-1}R$ equals $Q$.
Since $R$ is a $1$-Gorenstein ring, it follows from Lemma
\ref{lem6.10}(4) that $\Sigma^{-1}R$ also is a $1$-Gorenstein ring.
In addition, it follows from standard commutative algebra that the
map $\varphi:\Delta_1\to \mbox{Spec}(\Sigma^{-1}R)$ sending
$\mathfrak q$ to $\Sigma^{-1}\mathfrak q$ for $\mathfrak
q\in\Delta_1$ is a bijection. This shows that we can have the
following exact sequence of $R$-modules:
$$0\lra R_{\Sigma}\lraf{h} Q\lra \bigoplus_{\mathfrak q\in
\Delta_1}E(R/\mathfrak q)\lra 0.$$ By Lemma \ref{lem6.10}, we can
further form the following exact commutative diagram of $R$-modules:
\begin{eqnarray*}
\xymatrix{
&       & 0\ar[d]             &0\ar[d]         &\\
0\ar[r]&R\ar[r]^-{f_{\Sigma}}\ar@{=}[d]&\Sigma^{-1}R \ar[d]^-{h}
\ar[r]^-{h\pi}\ar[d]&Y'\ar[r]\ar[d]&0\\
0\ar[r]& R\ar[r]^-{f_\Phi} &Q\ar[d]\ar[r]^-{\pi}&
\bigoplus_{\mathfrak q\in P_1}E(R/\mathfrak q)\ar[r]\ar[d]\ar[d]^-{(g'_{\mathfrak q})_{{\mathfrak q}\in P_1}}&0,\\
&        &                         \bigoplus_{\mathfrak q\in
\Delta_1}E(R/\mathfrak q)\ar@{=}[r]\ar[d] & \bigoplus_{\mathfrak
q\in
\Delta_1}E(R/\mathfrak q)\ar[d]& \\
&        &                     0 &  0& \\
}
\end{eqnarray*}
where $Y'$ denotes the image of $h\pi$, and where $g'_{\mathfrak q}:
E(R/\mathfrak q)\ra \bigoplus_{\mathfrak q\in \Delta_1}E(R/\mathfrak
q)$ is a homomorphism of $R$-modules. By the same argument as in the
proof of  Lemma \ref{lem6.11}(2), we can prove
$Y'=\bigoplus_{\mathfrak p\in {P_1\setminus\Delta_1}}E(R/\mathfrak
p)$. Thus $\Sigma^{-1}R=\pi^{-1}(\bigoplus_{\mathfrak
p\in{P_1\setminus \Delta_1}}E(R/\mathfrak p))$.  By definition, we
have $P_1\setminus\Delta\subseteq\Delta_1$, and so
$\Sigma^{-1}R\subseteq R_{(\Delta)}$.

$(2)$ It follows from the definition of $ R_{(\Delta)}$ and the
statement $(1)$ that $R_{(\Delta)}=\Sigma^{-1}R $ if and only if
$\Delta_1=P_1\setminus\Delta$. To prove the statement $(2)$, it
suffices to show that if $R_{(\Delta)}=\Sigma_1^{-1}R$ for some
multiplication subset $\Sigma_1$ of $R$, then
$R_{(\Delta)}=\Sigma^{-1}R$. Now, assume
$R_{(\Delta)}=\Sigma_1^{-1}R$. By Lemma \ref{lem6.10}(4),
$R_{(\Delta)}$ is a $1$-Gorenstein ring. Note that $Q/R_{(\Delta)}
\simeq\bigoplus_{\mathfrak p\in {P_1\setminus\Delta}}E(R/\mathfrak
p)$ as $R$-modules. Then it follows from Lemma \ref{lem6.10}(1) and
Lemma \ref{lem6.9}(2) that, for any $\mathfrak p\in
{P_1\setminus\Delta}$, we have $\mathfrak p\cap\Sigma_1=\emptyset$,
and so $\Sigma_1\subseteq R\setminus\mathfrak p$. Since
$\Sigma:=R\setminus\bigcup_{\mathfrak q\in
{P_1\setminus\Delta}}\mathfrak q=\bigcap_{\mathfrak q\in
{P_1\setminus\Delta}}R\setminus\mathfrak q$, we have
$\Sigma_1\subseteq \Sigma$, and so
$R_{(\Delta)}=\Sigma_1^{-1}R\subseteq \Sigma^{-1}R$. Thanks to the
statement $(1)$, we get  $\Sigma^{-1}R\subseteq
R_{(\Delta)}\subseteq Q$. Thus $R_{(\Delta)}=\Sigma^{-1}R$.

$(3)$ It suffices to show $\Delta_1=P_1\setminus\Delta$. Clearly, we
have $P_1\setminus\Delta\subseteq\Delta_1$ by definition. Now we
show $\Delta_1\subseteq P_1\setminus \Delta$. In fact, if $\mathfrak
a\in\Delta_1$, then $\mathfrak a\subseteq\bigcup_{\mathfrak q\in
{P_1\setminus\Delta}}\mathfrak q$. Thus,

if, in addition, $P_1\setminus\Delta$ is finite, then $\mathfrak
a\subseteq\mathfrak q_1$ for some $\mathfrak
q_1\in{P_1\setminus\Delta}$ by prime avoidance theorem. Since
$\mathfrak a$ is a maximal ideal of $R$, it follows that $\mathfrak
a=\mathfrak q_1$. Hence $\Delta_1=P_1\setminus\Delta$.

If we assume that each ideal  in $\Delta$ is principal, then
$\mathfrak a$ must be in $P_1\setminus\Delta$. If it is not the
case, then $\mathfrak a\in \Delta$, and so there exists  an $r\in R$
such that $\mathfrak a=Rr$. Since $\mathfrak
a\subseteq\bigcup_{\mathfrak q\in {P_1\setminus\Delta}}\mathfrak q$,
we know that  $r\in\mathfrak q$ for some $\mathfrak q\in
{P_1\setminus\Delta}$, and so $\mathfrak a\subseteq \mathfrak q$. By
the maximality of $\mathfrak a$, we have $\mathfrak a=\mathfrak q$.
This is impossible because the intersection of $\Delta$ and
$P_1\setminus\Delta$ is empty. Hence $\Delta_1 =
P_1\setminus\Delta$. By $(2)$, we have $R_{(\Delta)}=\Sigma^{-1}R $
for either case. $\square$

\medskip
Combining  Corollary \ref{ringepi}(1) and Lemma \ref{epi00} with
Lemma \ref{lem6.11}(7), we have the following result on recollements
of derived module categories of endomorphism rings.

\begin{Prop}
Let $R$ be a $1$-Gorenstein ring, and let $\Delta$ be a subset of
$P_1$. Assume that each prime ideal in $P_1\setminus\Delta$ contains
all zero divisors of $R$. Then we get the following recollements of
derived module categories:

$$
\xymatrix@C=1.2cm{\D{Q\otimes_RJ_{\Delta}}\ar[r]&\D{B_\Delta}\ar[r]
\ar@/^1.2pc/[l]\ar@/_1.2pc/[l]
&\D{R}\ar@/^1.2pc/[l]\ar@/_1.2pc/[l]},\vspace{0.3cm}
$$
$$
\xymatrix@C=1.2cm{\D{Q\otimes_RJ'_{\Delta}}\ar[r]&\D{B'_\Delta}\ar[r]
\ar@/^1.2pc/[l]\ar@/_1.2pc/[l]
&\D{R_{(\Delta)}}\ar@/^1.2pc/[l]\ar@/_1.2pc/[l]}.\vspace{0.3cm}
$$
 \label{prop6.13}
\end{Prop}
{\it Proof.} Here we provide another proof. We consider the
injective homological ring epimorphism $\lambda_{\Delta}: R\ra
R_{(\Delta)}$ defined in Lemma \ref{lem6.11}(3). Then, we have
$R_{(\Delta)}\otimes_RJ_{\Delta}\simeq Q\otimes_RJ_{\Delta}$ by
Lemma \ref{lem6.11}(7). Now, the first recollement in Proposition
\ref{prop6.13} follows immediately from Proposition \ref{com1}.

The proof of the existence of the second recollement in Proposition
\ref{prop6.13} can be implemented similarly as we did for the first
one.
 $\square$

\medskip
In the rest of this subsection, we consider the ring $\mathbb{Z}$,
it is a Dedekind domain and, of course, a $1$-Gorenstein ring.
Clearly, it fulfills the assumption of Proposition \ref{prop6.13}.
In this case, we can have a more explicit formulation for
Proposition \ref{prop6.13}. Our discussion below uses some basic
results on $p$-adic numbers in algebraic number theory.

Fix a prime number $p\ge 2$. A $p$-adic integer is a formal infinite
series $\sum^{\infty}_{i=0}a_ip^i$, where $0\leq a_i<p$ for all $
i\geq 0$.  A $p$-adic number is a formal infinite series of the form
$\sum^{\infty}_{j=-m}a_jp^j$, where $m\in\mathbb{Z}$ and $0\leq
a_j<p$ for all $j\geq -m$. The sets of all $p$-adic integers and
$p$-adic numbers are denoted by ${\mathbb Z}_p$ and ${\mathbb Q}_p$,
respectively. Note that  ${\mathbb Z}_p$ is a discrete valuation
ring of global dimension $1$ with the unique maximal ideal
$p{\mathbb Z}_p$, and that ${\mathbb Q}_p$ is a field.

If $f\in\mathbb{Q}$ is a rational number, then we can write
$$ f=\frac{g}{h}p^{-m} \quad \mbox{where}\quad g, h\in\mathbb{Z}, \,(gh, p)=1.$$
Since the rational number $\frac{g}{h}$ always belongs to $\mathbb
Z_p$, that is, there are $0\leq a_i<p$  for all $ i\geq 0$ such that
$\frac{g}{h}=\sum^{\infty}_{i=0}a_ip^i.$ Consequently, we have
$$
f=\sum^{\infty}_{i=0}a_ip^{-m + i}\in{\mathbb Q}_p.
$$
In this way, we can regard $\mathbb{Q}$ as a subfield of ${\mathbb
Q}_p$.  This implies that, for $f\in\mathbb{Q}$, there are at most
finitely many prime numbers $q$ such that
$f\in{\mathbb{Q}_q\setminus\mathbb{Z}_q}$, or equivalently,
$f\in{\mathbb{Z}_q}$ for almost all prime number $q$. It is
well-known that ${\mathbb Q}\otimes_{\mathbb Z}{\mathbb Z}_p\simeq
{\mathbb Q}_p$ by multiplication map since ${\mathbb Q}_p
=\{p^my\mid m\in {\mathbb Z}, y\in {\mathbb Z}\}$. Clearly, $\mathbb
Q\subset {\mathbb Q}_p$ and $\mathbb Z\subset {\mathbb
Q}_{(p)}\subset {\mathbb Z}_p\subset{\mathbb Q}_p$ for every prime
$p\in {\mathbb N}:=\{0,1,2, \cdots \}$. It is known that
$\End_{\mathbb Z}({\mathbb Q}/{\mathbb Z})\simeq\prod_{p}{\mathbb
Z}_p$ as rings, where $p$ goes through all prime numbers.

An alternative definition of ${\mathbb Z}_p$ is that ${\mathbb Z}_p$
is the $p$-adic completion $\varprojlim_i{\mathbb Z}/p^i{\mathbb Z}$
of $\mathbb Z$. Another algebraic definition of ${\mathbb Z}_p$ is
that ${\mathbb Z}_p$ is isomorphic to the quotient of the formal
power series ring ${\mathbb Z}[[X]]$ by the ideal generated by
$X-p$. Note that ${\mathbb Q}_p$ is the field of fractions of
${\mathbb Z}_p$. For more details about $p$-adic numbers, one may
refer to \cite[Chapter IV, Section 2]{ne}. We denote by
$\widehat{\mathbb Z}$ the product $\prod_{p}{\mathbb Z}_p$ of
 all ${\mathbb Z}_p$ with $p$ positive prime numbers. This is a commutative ring.

Now, let $\Lambda$ be the set of all prime numbers in $\mathbb{N}$,
and let $I$ be a subset of $\Lambda$. Set
 $I':=\Lambda\setminus I,\, \Delta:=\{\mathfrak p\mid p\in I\}$
  and $\mathbb{Z}_{(I)}:=\mathbb{Z}_{(\Delta)}.$

\begin{Lem} The following statements hold true for the ring $\mathbb{Z}$ of integers.

$(1)$ Let $\Sigma:=\mathbb{Z}\setminus \cup_{q\in I'}\,\mathfrak q.$
Then $\mathbb{Z}_{(I)}=\Sigma^{-1}\mathbb{Z}$, which is the smallest
subring of $\mathbb{Q}$ containing $\frac{1}{p}$ for all $p\in I$.

$(2)$The injective ring homomorphism
$$\Theta_I:\, \mathbb{Q}\otimes_\mathbb{Z}\prod_ {p\in I}\mathbb{Z}_p
\lra\prod_ {p\in I}\mathbb{Q}_p$$ defined by $q\otimes(x_ p)_{p\in
I}\mapsto (qx_p)_{ p\in I}$ for $q\in \mathbb{Q}$ and $x_
p\in\mathbb{Z}_p$ satisfies that
$$\Img(\Theta_I)=\mathbb{A}_I:=\big\{(y_p)_ {p\in I} \in\prod_ {p\in
I}\mathbb{Q}_p\mid y_p\in\mathbb{Z}_p\; \mbox{for almost all } p\in
I\big\}.$$ In particular, if $I$ is a finite set, then
$\Img(\Theta_I)=\mathbb{A}_I=\prod_ {p\in I}\mathbb{Q}_p$. Note that
$\mathbb{A}_I$ is a kind of ad\'ele in global class field theory
(see \cite[Chapter VI]{ne}).

$(3)$ There are the following ring isomorphisms:
$$\mathbb{Q}\otimes_\mathbb{Z}\prod_ {p\in
I}\mathbb{Z}_p\simeq\mathbb{A}_I, \qquad
\mathbb{Q}\otimes_\mathbb{Z}\prod_ {p\in
I'}\mathbb{Z}_p\simeq\mathbb{A}_{I'}.
$$\label{6.10}
\end{Lem}
\bigskip

{\it Proof.} $(1)$ Let $q\in I'$. By Lemma \ref{lem6.11}(2), we have
$\mathbb{Z}_{(I)}=\bigcap_{q\in I'} \mathbb{Z}_{\mathfrak q}$, where
$\mathbb{Z}_{\mathfrak q}$ is the localization of $\mathbb{Z}$ at
$\mathfrak q$ with ${\mathfrak q}= q\mathbb{Z}$. It follows from
${\mathbb Z}_{\mathfrak q}={\mathbb Q}_{(q)}$ that
$$\mathbb{Z}_{(I)}=\bigcap_{q\in I'}{\mathbb Q}_{(q)}={\mathbb
Z}[p^{-1}\mid p\in I]=\Sigma^{-1}{\mathbb Z}.$$

$(2)$ For each prime number $p$, the canonical ring homomorphism
$\mu:{\mathbb Q}\otimes_{\mathbb Z}{\mathbb Z}_p\to{\mathbb Q}_p$,
defined by $f\otimes x_p\mapsto fx_p$ for any $f\in\mathbb{Q},\,
x_p\in{\mathbb Z}_p$, is an isomorphism. Moreover, for such
$f\in\mathbb{Q}$, there are at most finitely many prime numbers $q$
such that $f\in{\mathbb{Q}_q\setminus\mathbb{Z}_q}$. In other words,
$f\in{\mathbb{Z}_q}$ for almost all prime number $q$. This implies
$\Img(\Theta_I)=\mathbb{A}_I$.

$(3)$ This follows from (2). $\square$

\medskip
With help of Lemma \ref{6.10}, we can state Proposition
\ref{prop6.13} for $R = \mathbb{Z}$ more explicitly.

\begin{Koro}
We have the following recollements of derived module categories:
$$
\xymatrix@C=1.2cm{\D{\mathbb{A}_I}\ar[r]&\D{B_I}\ar[r]
\ar@/^1.2pc/[l]\ar@/_1.2pc/[l]
&\D{\mathbb{Z}}\ar@/^1.2pc/[l]\ar@/_1.2pc/[l]},\vspace{0.3cm}
$$
$$
\xymatrix@C=1.2cm{\D{\mathbb{A}_{I'}}\ar[r]&\D{B'_I}\ar[r]
\ar@/^1.2pc/[l]\ar@/_1.2pc/[l]
&\D{\mathbb{Z}_{(I)}}\ar@/^1.2pc/[l]\ar@/_1.2pc/[l]},\vspace{0.3cm}
$$
where
$B_I:=\End_\mathbb{Z}(\mathbb{Z}_{(I)}\oplus\mathbb{Z}_{(I)}/\mathbb{Z})$
and $B'_I:=\End_{\mathbb{Z}_{(I)}}(\mathbb{Q}\oplus
\mathbb{Q}/\mathbb{Z}_{(I)}).$ \label{6.11}
\end{Koro}

\subsection{Examples\label{stratification}}

\medskip In the following we shall exploit Corollary
\ref{6.11} to give a couple of examples of derived module categories
that have two different stratifications by derived module categories
of rings with different composition factors. This is related to the
following problem proposed in \cite{akl}:

{\bf Problem:} Given a ring $R$, do all stratifications of $\D{R}$
by derived module categories of rings have the same finite number of
factors, and are these factors the same for all stratifications, up
to ordering and up to derived equivalence?

A negative partial solution to this problem can be seen from
Examples \ref{6.6} and \ref{6.8} below.

Let us first recall the definition of a stratification of $\D{R}$
for $R$ a ring in \cite{akl}.

Let $R$ be a ring. If there are rings $R_1$ and $R_2$ such that a
recollement
$$(*)\qquad \xymatrix@C=1.2cm{\D{R_1}\ar[r]&\D{R}\ar[r]
\ar@/^1.2pc/[l]\ar@/_1.2pc/[l]
&\D{R_2}\ar@/^1.2pc/[l]\ar@/_1.2pc/[l]}\vspace{0.3cm}$$ exists, then
$R_i$ or $\D{R_i}$ are called factors of $R$ or $\D{R}$. In this
case, we also say that ($*$) is a recollement of $R$. The ring R is
called derived simple if $\D{R}$ does not admit any non-trivial
recollement whose factors are derived categories of rings. It is
pointed out in \cite{HKL} that every Dedekind ring ( thus every
discrete valuation ring) is derived simple.

A stratification of $\D{R}$ is defined to be a sequence of iterated
recollements of the following form: a recollement of $R$, if it is
not derived simple,
$$\xymatrix@C=1.2cm{\D{R_0}\ar[r]&\D{R}\ar[r]
\ar@/^1.2pc/[l]\ar@/_1.2pc/[l]
&\D{R_1}\ar@/^1.2pc/[l]\ar@/_1.2pc/[l]},\vspace{0.3cm}$$  a
recollement of $R_0$, if it is not derived simple,
$$\xymatrix@C=1.2cm{\D{R_{00}}\ar[r]&\D{R_0}\ar[r]
\ar@/^1.2pc/[l]\ar@/_1.2pc/[l]
&\D{R_{01}}\ar@/^1.2pc/[l]\ar@/_1.2pc/[l]},\vspace{0.3cm}$$ and a
reollement of $R_2$, if it is not derived simple,
$$\xymatrix@C=1.2cm{\D{R_{10}}\ar[r]&\D{R_1}\ar[r]
\ar@/^1.2pc/[l]\ar@/_1.2pc/[l]
&\D{R_{11}}\ar@/^1.2pc/[l]\ar@/_1.2pc/[l]}\vspace{0.3cm}$$ and
recollements of $R_{ij}$ with $0\le i,j \le 1$, if they are not
derived simple, and so on, until one arrives at derived simple rings
at all positions, or continue to infinitum. All the derived simple
rings appearing in this procedure are called composition factors of
the stratification. The cardinality of the set of all composition
factors (counting the multiplicity) is called the length of the
stratification. If this procedure stops after finitely many steps,
we say that this stratification is finite or of finite length.

The first example shows two stratifications of a derived module
category with infinitely many different derived simple module
categories as composition factors

\medskip
\begin{Bsp}{\rm Let $\mathbb Z\hookrightarrow \mathbb Q$ be the
inclusion. Then $T={\mathbb Q}\oplus {\mathbb Q}/{\mathbb Z}$ is a
tilting $\mathbb Z$-module, and $$B:= \End_R(T)=\left(\begin{array}{cc} {\mathbb Q} & \Hom_{\mathbb Z}({\mathbb Q},{\mathbb Q}/{\mathbb Z})\\
0 & \widehat{\mathbb Z}\end{array}\right).
$$}
\label{6.6}
\end{Bsp}
\noindent Note that $\Hom_{\mathbb Z}({\mathbb Q},{\mathbb
Q}/{\mathbb Z})\simeq {\mathbb R}$ as abelian groups, where $\mathbb
R$ is the field of real numbers.

We take $\Delta :={\rm mSpec}({\mathbb Z})$. By Proposition
\ref{prop6.13} and Lemma \ref{6.10}(3), we have a recollement:

$$\xymatrix@C=1.2cm{\D{{\mathbb Q}\otimes_{\mathbb Z}\widehat{\mathbb Z}}\ar[r]&\D{B}\ar[r]
\ar@/^1.2pc/[l]\ar@/_1.2pc/[l] &\D{{\mathbb
Z}}\ar@/^1.2pc/[l]\ar@/_1.2pc/[l]}.\vspace{0.3cm}$$

Let $e_2=(1, 0, \cdots, )\in \widehat{\mathbb Z}$. Then
$\widehat{\mathbb Z} = \widehat{\mathbb Z}e_2\oplus \widehat{\mathbb
Z}(1-e_2)$. This is a decomposition of ideals of $\widehat{\mathbb
Z}$. Thus we have a decomposition of ideals of the ring ${\mathbb
Q}\otimes_{\mathbb Z}\widehat{\mathbb Z}$:
$${\mathbb Q}\otimes_{\mathbb Z}\widehat{\mathbb Z} = {\mathbb Q}\otimes_{\mathbb Z}{\mathbb Z}_p\oplus {\mathbb Q}
\otimes_{\mathbb Z}\prod_{p\ge 3}{\mathbb Z}_p= {\mathbb Q}_p\oplus
{\mathbb Q} \otimes_{\mathbb Z}\prod_{p\ge 3}{\mathbb Z}_p.$$ This
procedure can be repeated infinitely many times. Then it follows
that $\D{{\mathbb Q}\otimes_{\mathbb Z}\widehat{\mathbb Z}}$ has a
derived composition series with infinitely many simple factors
$\D{{\mathbb Q}_p}$. This shows that $\D{B}$ has a stratification
with derived composition factors equivalent to either $\D{\mathbb
Z}$ or $\D{{\mathbb Q}_p}$, both are derived simple, that is, they
are not a middle term in any proper recollement of derived module
categories of rings.

Transparently, it follows from the triangular form of $B$ that
$\D{B}$ has a stratification with infinitely many composition
factors equivalent to either $\D{\mathbb Q}$ or $\D{{\mathbb Q}_p}$.
Clearly, $\D{\mathbb Z}$ and $\D{\mathbb Q}$ are not equivalent as
triangulated categories since the global dimension of $\mathbb Z$ is
one and the global dimension of $\mathbb Q$ is zero. Thus $\D{B}$
has two stratifications which have different composition factors.
This gives negatively an answer to the second question of the above
mentioned problem.

\medskip
In Example \ref{6.6} the two stratifications of the category $\D{B}$
by derived module categories have infinite many composition factors.
In the next example we shall see that even one requires finiteness
of stratifications of a derived module category, their composition
factors still may be different. This is contrary to the well-known
Jordan-H\"older theorem which says that any two (finite) composition
series of a group have the same list of composition factors (up to
the ordering and up to isomorphism).

\begin{Bsp} {\rm  (1) Let $I$ be a non-empty
finite subset of mSpec$({\mathbb Z})$. We consider the exact
sequence $$0\ra {\mathbb Z}\ra {\mathbb Z}_{(I)} \ra \bigoplus_{p\in
I} E({\mathbb Z}/p)\ra 0$$ of abelian groups. Then $T:= {\mathbb
Z}_{(I)}\oplus \bigoplus_{p\in I}E({\mathbb Z}/p{\mathbb Z})$ is a
tilting module. On the one hand, by Lemmata \ref{epi1}(2) and
\ref{lem6.11}(5), we have

$$ \End_{\mathbb Z}(T)\simeq \left(\begin{array}{lc} {\mathbb Z}_{(I)} & \Hom_{\mathbb Z}({\mathbb Z}_I,{\mathbb Z}_{(I)}/{\mathbb Z})\\
0 & \bigoplus_{p\in I}{\mathbb Z}_p\end{array}\right).
$$ }
\label{6.8}
\end{Bsp}
On the other hand, since $I$ is a finite set, by Corollary
\ref{6.11}, $\End_{\mathbb Z}(T)$ admits a recollement

$$\xymatrix@C=1.2cm{\D{\bigoplus_{p\in I}{\mathbb Q}_p}\ar[r]&\D{\End_{\mathbb Z}(T)}\ar[r]
\ar@/^1.2pc/[l]\ar@/_1.2pc/[l] &\D{{\mathbb
Z}}\ar@/^1.2pc/[l]\ar@/_1.2pc/[l]}.\vspace{0.3cm}$$ Thus,
$\D{\End_{\mathbb Z}(T)}$ admits two stratifications, one has the
composition factors ${\mathbb Z}_{(I)}$ and ${\mathbb Z}_p$ with
$p\in I$, and the other has the composition factors ${\mathbb Z}$
and ${\mathbb Q}_p$ with $p\in I$. Since ${\mathbb Z}_{(I)}$ is a
localization of $\mathbb Z$ by Lemma \ref{6.10}, it is of global
dimension one. Note that derived equivalences preserve the centers
of rings. This shows that all rings $\mathbb{Z, Z}_{(I)},
\mathbb{Z}_p$ and $\mathbb{Q}_p$ are pairwise not
derived-equivalent. Hence the two stratifications have completely
different composition factors.

(2) Let ${\mathfrak p}= p{\mathbb Z}\subset {\mathbb Z}$ with $p$ a
prime number in $\mathbb N$. We consider the exact sequence of
${\mathbb Z}_{\mathfrak p}$-modules:
$$ 0\ra {\mathbb Z}_{\mathfrak p}\ra {\mathbb Q}\ra E({\mathbb Z}_{\mathfrak p}/p{\mathbb Z}_{\mathfrak p})\ra 0.$$
Define $T:= {\mathbb Q}\oplus E({\mathbb Z}_{\mathfrak p}/p{\mathbb
Z}_{\mathfrak p}) $.  Thus, by Lemmata \ref{lem6.11} and \ref{6.11},
we have

$$ \End_{{\mathbb Z}_{\mathfrak p}}(T)\simeq \End_{\mathbb Z}(T)\simeq \left(\begin{array}{lc} {\mathbb Q} & {\mathbb Q}_p\\
0 & {\mathbb Z}_p\end{array}\right),$$ and a recollement:

$$\xymatrix@C=1.2cm{\D{{\mathbb Q}_p}\ar[r]&\D{\End_{{\mathbb Z}_{\mathfrak p}}(T)}\ar[r]
\ar@/^1.2pc/[l]\ar@/_1.2pc/[l] &\D{{\mathbb Z}_{\mathfrak
p}}\ar@/^1.2pc/[l]\ar@/_1.2pc/[l]}.\vspace{0.3cm}$$ Note that the
ring $\End_{{\mathbb Z}_{\mathfrak p}}(T)$ is left hereditary, but
not left noetherian.

On the one hand,  $\D{\End_{{\mathbb Z}_{\mathfrak p}}(T)}$ has
clearly a stratification of length $2$ with the composition factors
$\mathbb Q$ and ${\mathbb Z}_p$. On the other hand, it admits
another stratification of length $2$ with the composition factors
${\mathbb Q}_p$ and ${\mathbb Z}_{\mathfrak p}$. Note that ${\mathbb
Z}_{\mathfrak p}={\mathbb Q}_{(p)}$. Since ${\mathbb Z}_p$ and
${\mathbb Q}_p$ are uncountable sets and since derived equivalences
preserve the centers of rings, we deduce that neither ${\mathbb Q}$
and ${\mathbb Q}_{(p)}$, nor ${\mathbb Z}_p$ and ${\mathbb Q}_{(p)}$
are derived equivalent. Clearly, the global dimensions of ${\mathbb
Z}_p$ and ${\mathbb Q}_{(p)}$ are one. Thus we have proved that the
derived category of the ring $\End_{{\mathbb Z}_{\mathfrak p}}(T)$
has two stratifications of length two without any common composition
factors.

Thus, this example shows also that the main result in \cite[Theorem
6.1]{akl} for hereditary artin algebras  cannot be extended to left
hereditary rings.

Note that in each example given in this section the sets of
composition factors of the two stratifications of the derived module
category have the same cardinalities. In the next section we shall
see that this phenomenon is not always true.

\section{Further examples and open questions\label{sect8}}

The main purpose of this section is to present examples of derived
module categories of rings such that they possess two
stratifications (by derived module categories of rings) with
different finite lengths. Namely, we consider the following

\smallskip
{\bf Question.} Is there a ring $R$ such that $\D{R}$ has two
stratifications of different finite lengths by derived module
categories of rings ?

Thus we solve the whole problem in \cite{akl} negatively.

Let $k$ be a field. We denote by $k[x]$ and $k[[x]]$ the polynomial
and formal power series algebras over $k$ in one variable $x$,
respectively, and by $k((x))$ the Laurent power series algebra in
one variable $x$, that is, $k((x)):=\{x^{-n}a\mid n\in {\mathbb N},
a\in k[[x]]\}$.

Now, let $k$ be an algebraically closed field, and
let $R$ be the Kronecker algebra $\left(\begin{array}{cc} k & k^2\\
0 & k\end{array}\right)$. It is known that $R$ can be given by the
following quiver
$$Q:\;
\xymatrix{2\ar@<0.4ex>[r]^-{\alpha}\ar@<-0.4ex>[r]_-{\beta}& 1},$$
and that $R\Modcat$ is equivalent to the category of representations
of $Q$ over $k$.

Let $V$ be a simple regular $R$-module. For each $m>0$, we denote by
$V[m]$ the module of regular length $m$ on the ray
$$V=V[1]\subset V[2]\subset\cdots\subset V[m]\subset
V[m+1]\subset\cdots,$$ and let $V[\infty]={\underrightarrow{\lim}}\;
V[m]$ be the corresponding $Pr\ddot{u}fer$ modules. Note that the
only regular submodule of $V[\infty]$  of regular length $m$ is
$V[m]$ with its canonical inclusion in $V[\infty]$, and that each
endomorphism of $V[\infty]$ in $R\Modcat$ restricts to an
endomorphism of $V[m]$ for any $m>0$. Thus, $V[\infty]$ admits a
unique chain of regular submodules. For more details, we refer to
\cite[Section 4.5]{R}.

From now on, we denote by $V$ the simple regular $R$-module:
$\xymatrix{k\ar@<0.4ex>[r]^-{0}\ar@<-0.4ex>[r]_-{1}& k}.$

Let  $e_1 = \left(\begin{array}{ll} 1 &0\\
0 & 0\end{array}\right)$ and  $e_2=\left(\begin{array}{ll} 0 &0\\
0 & 1\end{array}\right)$. Since $\Hom_R(Re_1,Re_2)\simeq k^2$, we
can identify a homomorphism from $Re_1$ to $Re_2$ in $R\Modcat$ with
an element in $k^2$. Fix a minimal projective resolution of $V$:
 $$\begin{CD}0@>>> Re_1@>{\partial:=(1,0)}>> Re_2@>>> V @>>> 0,\end{CD}$$
and denote by $\lambda:R\to R_{V}$ the universal localization of $R$
at the set $\Sigma:=\{\partial\}$.

It follows from \cite[Theorem 4.9, 5.1, and 5.3]{Sch} that $R_{V}$
is hereditary, $\lambda$ is injective, and $R_{V}\oplus R_{V}/R$ is
a tilting $R$-module. Moreover, by \cite[Proposition 1.8]{HJ2}, we
get $R_{V}/R\simeq V[\infty]^2$ as $R$-modules. Note that
$\Hom_R(R_V/R,R_V)=0$ because $R_{V}/R$ is a torsion module and
$R_V$ is a torsion-free module.

For simplicity of notation, we denote by $T$ the tilting module
$R_{V}\oplus V[\infty]^2$. Now, applying Corollary \ref{6.15} to the
module $T$, we can get the following recollement of derived module
categories:
$$(\ast)\quad
\xymatrix@C=1.2cm{\D{R_{V}\sqcup_RS'}\ar[r]&\D{B}\ar[r]
\ar@/^1.2pc/[l]\ar@/_1.2pc/[l]
&\D{R}\ar@/^1.2pc/[l]\ar@/_1.2pc/[l]},\vspace{0.3cm}$$ where
$B:=\End_R(T)$, $S':=M_2\big(\End_R(V[\infty])\big)$ and
$R_{V}\sqcup_RS'$ is the coproduct of $R_{V}$ and $S'$ over $R$.

In the following, we shall describe the rings $B$, $S'$ and
$R_{V}\sqcup_RS'$ explicitly.

First,  by Lemma \ref{lem32}, we can check that $R_{V}=M_2(k[x])$,
the $2\times 2$ matrix algebra over $k[x]$, and the map
$\lambda:R\to R_{V}$ is given by
$\left(\begin{array}{cc} a& (c,d)\\
0 & b\end{array}\right)\mapsto\left(\begin{array}{cc} a
& c+dx\\
0 & b \end{array}\right)$ for any $a,b,c, d\in k$. This means $R_V
e_1\simeq R_V e_2$ as $R_V$-modules, and therefore we have the
following ring isomorphisms:
$$(\ast\ast)\quad
B\simeq M_2\Big(\End_R(R_{V}e_1\oplus V[\infty])\Big) \simeq
M_2\Bigg(\left(\begin{array}{cc}e_1R_V e_1&
\Hom_R(R_V e_1, V[\infty])\\
0 & \End_R(V[\infty])\end{array}\right)\Bigg).$$

Second, we  claim that $\End_R(V[\infty])$ is isomorphic to
$k[[x]]$. In fact, this follows from the following isomorphisms of
abelian groups:
$$
\End_R(V[\infty])\simeq {\underleftarrow{\lim}}\;\Hom_R(V[m],
V[\infty])\simeq {\underleftarrow{\lim}}\; \Hom_R(V[m], V[m]) \simeq
{\underleftarrow{\lim}}\; k[x]/(x^m)\simeq k[[x]],
$$
where the composition of the above isomorphisms gives rise to a ring
isomorphism $\omega:\End_R(V[\infty])\to k[[x]]$. Thus $S'\simeq
M_2(k[[x]])$ as rings. In this sense, we can identify $S'$ with
$M_2(k[[x]])$ under the isomorphism $\omega$.

Third, a direct calculation shows that the ring homomorphism
$\mu:R\to S'$, which appears in the proof of Corollary \ref{6.15},
is given by
$\left(\begin{array}{cc} a'& (c',d')\\
0 & b'\end{array}\right)\mapsto\left(\begin{array}{cc} a'
& d'+c'x\\
0 & b' \end{array}\right)$ for any $a',b',c', d'\in k$.

Finally, we claim  $R_{V}\sqcup_RS'\simeq M_2\big(k((x))\big)$ as
rings.

Recall that $R_{V}$ is the universal localization of $R$ at
$\Sigma:=\{\partial\}$. Define $\varphi:=S'\otimes_{R}\partial:
S'e_1\to S'e_2$. Then it follows from Lemma \ref{epi01} that
$R_V\sqcup_RS'$ is isomorphic to the universal localization
$S'_{\varphi}$ of $S'$ at $\varphi$. Since
$\Hom_{S'}(S'e_1,S'e_2)\simeq e_1S'e_2\simeq k[[x]]$, the map
$\varphi$ corresponds to the matrix element
$\left(\begin{array}{cc} 0 & x\\
0 & 0\end{array}\right)$ in $S'$. Now, let $\rho_x: k[[x]]\to k[[x]]
$ be the right multiplication map by  $x$. Since  $S'$ is Morita
equivalent to $k[[x]]$, we conclude from Lemma \ref{Morita} that
$S'_{\varphi}=M_2\big(k[[x]]_{\rho_x}\big)$, where $k[[x]]_{\rho_x}$
is the universal localization of $k[[x]]$ at  $\rho_x$. Since
$k[[x]]$ is commutative, the ring $k[[x]]_{\rho_x}$ is isomorphic to
the localization $\Theta^{-1}k[[x]]$ of $k[[x]]$ at the
multiplication subset $\Theta:=\{x^m\mid m\in\mathbb{N}\}$. Thus
$\Theta^{-1}k[[x]]$ is the Laurent power series ring $k((x))$.
Therefore, we get the following isomorphisms of rings:
$$R_V\sqcup_RS'\simeq S'_{\varphi}\simeq
M_2\big(k[[x]]_{\rho_x}\big)\simeq
M_2\big(\Theta^{-1}k[[x]]\big)\simeq M_2\big(k((x))\big).$$

On the one hand, by setting $C:=\End_R(R_{V}e_1\oplus V[\infty])$
and using Morita equivalences, the recollement ($*$) can be
rewritten as

$$
\xymatrix@C=1.2cm{\mathscr{D}\big(k((x))\big)\ar[r]&\D{C}\ar[r]
\ar@/^1.2pc/[l]\ar@/_1.2pc/[l]
&\D{R}\ar@/^1.2pc/[l]\ar@/_1.2pc/[l]}.\vspace{0.3cm}
$$

On the other hand, since $e_1R_V e_1\simeq k[x]$ and
$\End_R(V[\infty])\simeq k[[x]]$, it follows from $(\ast\ast)$ that
the ring $C$ admits another recollement
$$
\xymatrix@C=1.2cm{\mathscr{D}\big(k[x]\big)\ar[r]&\D{C}\ar[r]
\ar@/^1.2pc/[l]\ar@/_1.2pc/[l]
&\D{k[[x]]}\ar@/^1.2pc/[l]\ar@/_1.2pc/[l]}.\vspace{0.3cm}
$$
Since derived equivalences preserve the centers of rings,  all rings
$k$, $k[x]$, $k[[x]]$ and $k((x))$ are pairwise not derived
equivalent. Moreover, they are derived simple. Clearly, $\D{R}$ has
a stratification of length 2 with composition factors $\D{k}$ and
$\D{k}$. Thus $C$ admits two stratifications, one of which is of
length $3$ with three composition factors $k((x))$, $k$ and $k$, and
the other is of length $2$  with composition factors $k[x]$ and
$k[[x]]$.  As a result, we have shown that the two stratifications
of $\D{C}$ by derived categories of rings are of different lengths
and without any common composition factors.

\medskip
{\it Remarks.} $(1)$ For any simple regular $R$-module $V'$, we can
choose an automorphism $\sigma:R\to R$, such that the induced
functor $\sigma_*: R\Modcat\to R\Modcat$ by $\sigma$ is an
equivalence and satisfies $\sigma_*(V')\simeq V$. Hence, instead of
$V$, we may use $V'$ to proceed the above procedure, but we will get
the same recollements, up to derived equivalence of each term.

$(2)$ Let $K_0(R)$ be the Grothendieck group of $R$, that is, the
abelian group generated by isomorphism classes $[P]$ of finitely
generated projective $R$-modules $P$ subject to the relation
$[P]+[Q]=[P\oplus Q]$, where $P$ and $Q$ are finitely generated
projective $R$-modules. One can check that
$K_0\big(k((x))\big)\simeq \mathbb{Z}$ and
$K_0(C)\simeq\mathbb{Z}\oplus \mathbb{Z}$. The above example shows
that, even if $\D{A_2}$ is a recollement of $\D{A_1}$ and $\D{A_3}$,
where $A_i$ are  rings for $i=1,2,3,$ we cannot get $K_0(A_2)\simeq
K_0(A_1)\oplus K_0(A_3)$ in general.

\medskip
For a general consideration of stratifications of the endomorphism
algebras of tilting modules over tame hereditary algebras, we shall
discuss it in a forthcoming paper.

\medskip
Finally, we remark that Theorem \ref{th01}(2) can be extended to
$n$-tilting modules. However, since there is not defined any
reasonable torsion theory in module categories for general
$n$-tilting modules, we are not able to extend Theorem \ref{th01}(1)
to $n$-tilting modules. So we mention the following open question.

\smallskip
{\bf Question 1.} Is Theorem \ref{th01}(1) true for $n$-good tilting
modules ?

\smallskip
Another question related to our examples is:

\smallskip
{\bf Question 2.} Is there a ring $R$ such that $\D{R}$ has two
stratifications by derived module categories of rings, one of which
is of finite length, and the other is of infinite length ?




\bigskip
October 10, 2010
\end{document}